\documentclass[10pt]{article}
\usepackage{amsmath, amssymb, amsthm}
\usepackage{mathrsfs}
\newcommand{\no}[1]{#1}
\renewcommand{\no}[1]{}
\no{\usepackage{times}\usepackage[subscriptcorrection, slantedGreek, nofontinfo]{mtpro}
\renewcommand{\Delta}{\upDelta}}
\usepackage{color}
\usepackage{enumerate,comment}
\usepackage{amssymb}
\usepackage{graphicx}
\usepackage{dsfont}
\usepackage{epstopdf}
\graphicspath{{./pics/}}
\usepackage{booktabs}
\usepackage{threeparttable}
\usepackage{url}

\date{\today}
\newcommand{\bel}{\begin{equation} \label}
\newcommand{\ee}{\end{equation}}
\def\beq{\begin{equation}}
\def\eeq{\end{equation}}
\newcommand{\bea}{\begin{eqnarray}}
\newcommand{\eea}{\end{eqnarray}}
\newcommand{\beas}{\begin{eqnarray*}}
\newcommand{\eeas}{\end{eqnarray*}}

\newcommand{\R}{\mathbb{R}}

\newtheorem{theorem}{Theorem}[section]
\newtheorem{lemma}[theorem]{Lemma}

\newtheorem{proposition}[theorem]{Proposition}
\newtheorem{defn}[theorem]{Definition}
\newtheorem{remark}{Remark}[section]
\newtheorem{example}{Example}[section]
\numberwithin{equation}{section}

\renewcommand{\d}{\mathrm{d}}
\allowdisplaybreaks

\providecommand{\norm}[1]{\left\lVert#1\right\rVert}

\def\phi {\varphi}

\title{Point Source Identification in Subdiffusion from A Posteriori Internal Measurement\thanks{The work of K. Huang is supported by a start-up fund from The Chinese University of Hong Kong, and that of B. Jin  by Hong Kong RGC General Research Fund (Project 14306423), and a
start-up fund from The Chinese University of Hong Kong. The work of Z. Zhou is supported by Hong Kong
Research Grants Council (15303122) and an internal grant of Hong Kong Polytechnic University (Project
ID: P0038888, Work Programme: 1-ZVX3).}}

\author{Kuang Huang\thanks{Department of Mathematics, The Chinese University of Hong Kong, Shatin, New Territories, Hong Kong, P.R. China (\texttt{kuanghuang@cuhk.edu.hk, b.jin@cuhk.edu.hk, bangti.jin@gmail.com})} \and Bangti Jin\footnotemark[2]\and
Yavar Kian\thanks{Univ Rouen Normandie, CNRS, Laboratoire de Math\'{e}matiques Rapha\"{e}l Salem, UMR 6085, F-76000 Rouen, France
(\texttt{yavar.kian@univ-rouen.fr}, \texttt{georges.sadaka@univ-rouen.fr})} \and Georges Sadaka\footnotemark[3] \and Zhi Zhou\thanks{Department of Applied Mathematics, The Hong Kong Polytechnic University, Kowloon, Hong Kong, P.R. China (\texttt{zhizhou@polyu.edu.hk})}}

\date{\today}

\begin{document}

\maketitle

\begin{abstract}
In this work we investigate an inverse problem of recovering point sources and their time-dependent
strengths from \textit{a posteriori} partial internal measurements in a subdiffusion model which involves a
Caputo fractional derivative in time and a general second-order elliptic operator in space. We establish the
well-posedness of the direct problem in the sense of transposition and improved local regularity. Using classical unique continuation of the subdiffusion model and improved local solution regularity, we prove the
uniqueness of simultaneously recovering the locations of point sources, time-dependent strengths and initial condition
for both one- and multi-dimensional cases. Moreover, in the one-dimensional case, the elliptic operator can have
time-dependent coefficients. These results extend existing studies on point source identification for parabolic type problems.
Additionally we present several numerical experiments to show the feasibility of numerical reconstruction.\\
\textbf{Keywords}:  subdiffusion, point source identification, uniqueness, a posteriori measurement    
\end{abstract}


\section{Introduction}

This work investigates an inverse source problem for the subdiffusion model. Let
$\Omega\subset \mathbb{R}^d$ ($d= 2,3$) be an open bounded and connected subset with a
$ C^{2}$ boundary $\partial \Omega$. We define an elliptic operator $\mathcal{A}$  by
\begin{equation}\label{A}
\mathcal{A} u(x) :=-\sum_{i,j=1}^d \partial_{x_i} \left( a_{i,j}(x) \partial_{x_j} u(x) \right)+\sum_{k=1}^db_k(x)\partial_{x_k}u+q(x)u(x),\quad  x\in\Omega,
\end{equation}
where $q \in C^1(\overline{\Omega})$ is nonnegative, $b_1,\ldots,b_k\in C^1(\overline{\Omega})$
and $a:=(a_{i,j})_{1 \leq i,j \leq d} \in C^{\infty}
(\R^d;\mathbb{R}^{d\times d})\cap W^{\infty,\infty}
(\R^d;\mathbb{R}^{d\times d})$ is symmetric and fulfills the standard ellipticity condition
\begin{equation}\label{ell}
\exists c>0 :\,  \sum_{i,j=1}^d a_{i,j}(x) \xi_i \xi_j \geq c |\xi|^2,\quad  x \in \R^d,\ \xi=(\xi_1,\ldots,\xi_d) \in \mathbb{R}^d.
\end{equation}
Let the weight function $\rho\in C^{\infty}
(\R^d)\cap W^{\infty,\infty}(\R^d)$ satisfy the following condition
\begin{equation}\label{eqn:rho}
 0<c_0 \leq\rho(x) ,\quad x\in \R^d,
\end{equation}
with $c_0>0$ being a constant. Throughout, we set $\R_+:=(0,\infty)$ and $\mathbb N:=\{1,2,\ldots\}$. The Riemann-Liouville fractional integral $_0I_t^\beta$ and the Riemann-Liouville fractional derivative  $D_t^\beta$, of order $\beta\in(0,1)$, are defined, respectively, by (see, e.g., \cite[p. 21 and p. 28]{Jin:2021book})
\[
_0I_t^\beta h(\cdot,t):=\frac{1}{\Gamma(\beta)}\int_0^t\frac{h(\cdot,s)}{(t-s)^{1-\beta}}\,\d s\quad \mbox{and}\quad D_t^\beta:=
\partial_t\circ {_0I_t^{1-\beta}}.
\]
The Caputo fractional derivative  $\partial_t^\beta$, of order $\beta\in(0,1)$, is defined by \cite[p. 41]{Jin:2021book}
$$\partial_t^\beta h=D_t^\beta (h-h(\cdot,0)),\quad h\in  C([0,\infty);L^2(\Omega)).$$
For all $\beta\in(0,1)$, we have
$\partial_t^\beta h={_0I_t^{1-\beta}}\partial_th$ for $h\in W^{1,1}_{loc}(\mathbb R_+;L^2(\Omega)).$
Fix $T>0$, $N\in\mathbb N$, $\{x_k\}_{k=1}^N\subset \Omega$ (distinct from each other), $\{\lambda_k\}_{k=1}^N\subset L^2(0,T)$,  $u_0\in H^1_0(\Omega)$ and $\alpha \in (0, 1)$.
Consider the following initial boundary value problem (IBVP):
\begin{equation}\label{eq1}
\left\{\begin{aligned}
\rho(x)\partial_t^{\alpha}u +\mathcal{A} u &=  \sum_{k=1}^N\lambda_k(t)\delta_{x_k}, && \mbox{in }\Omega\times(0,T),\\
 u&= 0, && \mbox{on } \partial\Omega\times(0,T), \\
u&=u_0, && \mbox{in } \Omega\times \{0\},
\end{aligned}\right.
\end{equation}
where $\delta_{x_k}(x)$ denotes Dirac delta function concentrated at $x_k$. When $\alpha\to1^-$, the model \eqref{eq1} recovers the standard parabolic problem. The model \eqref{eq1}
can accurately describe anomalously slow diffusion (i.e., subdiffusion) processes. At a microscopic level, these processes can be modelled by continuous time random walk with the waiting time between consecutive jumps following a heavy-tailed distribution, and the subdiffusion model describes the evolution of the probability density function of particles appearing at time $t$ and spatial location $x$. It has found applications in physics, engineering and biology, e.g., diffusion in media with fractal geometry \cite{Nigmatulin:1986}, solute transport in heterogeneous media \cite{AdamsGelhar:1992,HatanoHatano:1998,berkowitz2006modeling} and protein transport in membranes \cite{kou2004generalized}.

The existing solution theory for subdiffusion \cite{KRY,Jin:2021book} does not apply directly to problem \eqref{eq1} due to the presence of point sources. We shall establish the well-posedness of problem \eqref{eq1} in the transposition sense \cite{LM1} and local solution regularity. Now fix $\omega$ to be a nonempty open subset of $\Omega$, and consider the inverse problem of simultaneously determining the source parameters $\{(x_k,\lambda_k)\}_{k=1}^N$ and initial condition $u_0$ from the knowledge of $u(x,t)$, $(x,t)\in\omega\times(T-\epsilon,T)$ with  $\epsilon\in(0,T)$ arbitrarily small. It arises naturally in the context of pollution detection in underground water, when the transport process of the pollutant is described by anomalously slow diffusion. In this context, $u$ denotes the concentration of the pollutant, and $\{x_k\}_{k=1}^N$ and $\{\lambda_k\}_{k=1}^N$ are the locations and release strengths of the pollutant, and the inverse problem is to determine the pollution sources from observations made at monitoring stations after the release of pollutant.

Now we give the unique identifiability of the initial condition, and locations and strengths of point sources.
\begin{theorem}\label{t1} For $j=1,2$, let $N_j\in\mathbb N$, $\{\lambda_k^j\}_{k=1}^{N_j}\subset H_\alpha(0,T)$ be $N_j$ non-uniformly vanishing functions, $\{x_k^j\}_{k=1}^{N_j}$ be $N_j$ distinct points in $\Omega$ and let $u_0^j\in H^1_0(\Omega)$. Fix also $\omega\subset \Omega\setminus\{x_1^1,\ldots,x_{N_1}^1,x_1^2,\ldots,x_{N_2}^2\}$ an open nonempty subset of $\Omega$ and, for $j=1,2$, let $u_j$ be the solution in the transposition sense of problem \eqref{eq1} with $u_0=u_0^j$, $N=N_j$, $\lambda_k=\lambda_k^j$ and $x_k=x_k^j$, $k=1,\ldots,N_j$. Then, for every $\epsilon\in(0,T)$, the condition
\begin{equation}\label{t1a}
u_1(x,t)=u_2(x,t),\quad (x,t)\in\omega\times (T-\epsilon,T)
\end{equation}
implies
\begin{equation}
\label{t1b}N_1=N_2=N\quad \mbox{and}\quad u_0^1=u_0^2,
\end{equation}
and moreover, there exists a one-to-one map $\pi$ of $\{1,\ldots,N\}$ such that
\begin{equation}
\label{t1c}x_k^1=x_{\pi(k)}^2,\quad \lambda_k^1=\lambda_{\pi(k)}^2,\quad k=1,\ldots,N.
\end{equation}
\end{theorem}

In the 1D case, it is impossible to recover point source when $N\geq2$ using the datum over one subinterval; see Example \ref{exam:1d} for a counterexample. Thus we restrict the analysis to $N=1$. Let $\ell>0$, $\rho \in C^{\infty}
(\R;(0,\infty))\cap W^{\infty,\infty}(\R)$, and the operator $\mathcal A$ be defined by
$$\mathcal{A}(t) w(x) :=-\partial_x(a(x,t)\partial_xw)(x)+b(x,t)\partial_xw(x)+q(x,t)w(x),\quad (x,t)\in (0,\ell)\times(0,T),$$
with $a\in C^{\infty}
(\R^2;(0,\infty))\cap W^{\infty,\infty}(\R^2)$ and $b, q\in C^1([0,\ell]\times[0,T])$. Consider  the following IBVP:
\begin{equation}\label{eqq11}
\left\{\begin{aligned}
\rho(x)\partial_t^{\alpha}u +\mathcal{A}(t) u &=  \lambda_0(t)\delta_{x_0}, && \mbox{in }(0,\ell)\times(0,T),\\
 u(0,t)=u(\ell,t)&= 0, && t\in(0,T), \\
u&=u_0, && \mbox{in } (0,\ell)\times \{0\}.
\end{aligned}\right.
\end{equation}

The result in the 1D case is stated below, which notably allows also time-dependent coefficients.
\begin{theorem}\label{t1D} For $j=1,2$, let  $\lambda^j\in L^2(0,T)$ be non-uniformly vanishing, $x_0^j\in(0,\ell)$ and let $u_0^j\in H^1_0(0,\ell)$. Fix  $I_1\subset (0,\min(x_0^1,x_0^2))$ and $I_2\subset (\max(x_0^1,x_0^2),\ell)$ be two open intervals and, for $j=1,2$, let $u_j$ be the solution  of \eqref{eqq11} with $u_0=u_0^j$,  $\lambda=\lambda^j$ and $x_0=x_0^j$. Further, assume  that
\begin{equation}
\label{ini} u_0^1(x)=u_0^2(x),\quad x\in[\min(x_0^1,x_0^2),\max(x_0^1,x_0^2)].
\end{equation}
Then, for every $\epsilon\in(0,T)$, the condition
\begin{equation}
\label{t1da} u_1(x,t)=u_2(x,t),\quad (x,t)\in(I_1\cup I_2)\times (T-\epsilon,T)
\end{equation}
implies
\begin{equation}
\label{t1db} x_0^1=x_0^2,\quad \lambda^1=\lambda^2,\quad u_0^1=u_0^2.
\end{equation}
\end{theorem}

Inverse source problems for subdiffusion have been intensively studied \cite{JinRundell:2015,LLY2}. However, the study on the recovery of point sources remains scarce: the presence of point sources greatly complicates the analysis of the direct problem; see Section \ref{sec:direct} for a new solution theory in the sense of transposition \cite{LM1}. Hrizi et al \cite{HHN} investigated the case of constant intensities for all point sources, i.e., $F(x,t)=\sum_{k=1}^N \lambda_k\delta_{x_k}(x)$, with $\lambda_k\in \mathbb{R}$, in the subdiffusion model with the Laplace operator, from the measurement taken on a subdomain $\omega$ over the whole time horizon $(0,T)$.
In contrast, Theorem \ref{t1} allows multiple point sources and time-dependent intensities, and the elliptic operator $\mathcal{A}$ can have spatially variable coefficients and is not necessarily self-adjoint. Thus it substantially generalizes the setting of the existing study \cite{HHN}. Moreover, we utilize only \textit{a posteriori} measurement that is made over a time interval away from $t=0$. This is expected to be of practical relevance, since often one can only make the measurement after the occurrence of  excitation, e.g., in pollution detection. This relies heavily on the memory effect of subdiffusion. See Remark \ref{rmk:memory} for more details. In addition, we provide several numerical experiments confirming the feasibility of the numerical reconstruction. The improved regularity in Theorem \ref{t2} gives a strong solution away from the point sources, which enables applying classical unique continuation principle \cite{LN}; see Lemma \ref{l2} for the precise form. To the best of our knowledge, we give a first application of unique continuation principle in \cite{LN} to inverse source problems. Meanwhile, several uniqueness results of recovering one point source
in the standard parabolic problem have been obtained, either from boundary measurement or from observations at multiple internal points; see, e.g., \cite{ElBadia:2002jiip,ElBadiaHamdi:2005,KusiakWeatherwax:2007}. In contrast, we can determine multiple point sources and the initial condition and utilize \textit{a posteriori} measurement due to the memory effect of the mathematical model, which has not been investigated in the standard parabolic case.

The rest of the paper is organized as follows. In Section \ref{sec:direct}, we develop a solution theory for the direct problem, including unique existence in the sense of transposition and improved local regularity. In Section \ref{sec:proof}, we give detailed proofs of Theorems \ref{t1} and \ref{t1D}. Last, in Section \ref{sec:numer}, we present  numerical experiments to complement the analysis. Throughout $c$ denotes a generic constant which may change at each occurrence.

\section{Analysis of the direct problem}\label{sec:direct}
Existing solution theory for subdiffusion \cite{KRY,Jin:2021book} does not directly cover point sources. We shall present a solution theory in the sense of transposition \cite{LM1} and establish improved local regularity. For any Banach space $X$ and $T\in(0,\infty]$, we denote by $H_1(0,T;X)$ the space of functions $u\in H^1(0,T;X)$ satisfying  $u(0)=0$. In view of \cite{LM1} (see also \cite[p. 12]{KRY}), for all $\beta\in (0,1)$, we denote by $H_\beta(0,T;X)$ the interpolation space of order $1-\beta$ between $H_1(0,T;X)$ and $L^2(0,T;X)$ which is a subspace of $H^\beta(0,T;X)$. We set also $H_\beta(0,T):=H_\beta(0,T;\R)$. Throughout $Q=\Omega\times (0,T)$ denotes the space-time cylinder.
\subsection{Solution in the transposition sense}\label{ssec:transposition}
 Consider the following IBVP:
\begin{equation}\label{eq2}
\left\{\begin{aligned}
\rho\partial_t^{\alpha}v +\mathcal{A} v &=  \sum_{k=1}^N\lambda_k(t)\delta_{x_k}, && \mbox{in }\Omega\times(0,T),\\
 v&= 0, && \mbox{on } \partial\Omega\times(0,T), \\
v&=0, && \mbox{in } \Omega\times \{0\}.
\end{aligned}\right.
\end{equation}

Now we study solutions of problem \eqref{eq2} in the transposition sense \cite[Chapter 3, Section 2]{LM1}. By Sobolev embedding theorem \cite[Theorem 4.12, p. 85]{Adams:2023}, $H^2(\Omega)$ embeds continuously into $C(\overline{\Omega})$ for $d=1,2,3$. Thus, the maps $\{\delta_{x_k}\}_{k=1}^N$ are continuous linear forms on $H^2(\Omega)$ and for any $v\in H^2(\Omega)$, we have $\delta_{x_k}(v)=v(x_k)$, $k=1,\ldots,N$. This motivates the following definition of solutions in the transposition sense.
\begin{defn}\label{def:transposition}
A function $v\in L^2(0,T;L^2(\Omega))$ is said to be a solution of problem \eqref{eq2} in the transposition sense if for all $w\in H_\alpha(0,T;L^2(\Omega))\cap L^2(0,T;H^2(\Omega)\cap H^1_0(\Omega))$, we have
\begin{equation*} \int_0^T\int_\Omega v(x,t)[\rho(x) D_t^\alpha w(x,T-t)+\mathcal A^* w(x,T-t)]\,\d x\d t=\sum_{k=1}^N\int_0^T\lambda_k(t)w(x_k,T-t)\,\d t,
\end{equation*}
where
$\mathcal{A}^* v :=-\sum_{i,j=1}^d \partial_{x_i} ( a_{i,j} \partial_{x_j} v-\sum_{k=1}^db_k\partial_{x_k}v+(q-\sum_{k=1}^d\partial_{x_k}b_k)v$ is the adjoint of $\mathcal{A}$. 
\end{defn}

The first result gives the unique existence of a solution to  \eqref{eq2} in the sense of transposition.
\begin{proposition}\label{p1} Problem \eqref{eq2} has a unique solution $v\in L^2(0,T;L^2(\Omega))$ in the transposition sense. Moreover, there exists $c>0$ depending only on $\rho$, $\mathcal A$, $\Omega$ and $T$ such that
\begin{equation}\label{p1a}
\norm{v}_{L^2(0,T;L^2(\Omega))}\leq c\sum_{k=1}^N\norm{\lambda_k}_{L^2(0,T)},
\end{equation}
and $v$ fulfills the condition
\begin{equation}\label{p1aa}\left\langle \rho(x)D_t^\alpha v+\mathcal A v,\psi\right\rangle_{D'(Q),C^\infty_0(Q)}=\sum_{k=1}^N\int_0^T\lambda_k(t)\psi(x_k,t)\,\d t,\quad \forall \psi\in C^\infty_0(Q).
\end{equation}
\end{proposition}
\begin{proof} Let $\zeta = \sum_{k=1}^N\norm{\lambda_k}_{L^2(0,T)}$.
By \cite[Theorem 2.1]{KRY}, for any $F\in L^2(0,T;L^2(\Omega))$, there exists a unique $w_F\in H_\alpha(0,T;L^2(\Omega))\cap L^2(0,T;H^2(\Omega)\cap H^1_0(\Omega))$ such that
\begin{equation}\label{eqn:w-F} 
\rho(x) D_t^\alpha w_F(x,t)+\mathcal A^* w_F(x,t)=F(x,T-t),\quad (x,t)\in \Omega\times(0,T).
\end{equation}
Moreover, the map $w_F$ satisfies the estimate
\begin{equation}
\label{p1b}\norm{w_F}_{L^2(0,T;H^2(\Omega))}\leq c\norm{F}_{L^2(0,T;L^2(\Omega))},
\end{equation}
with the constant $c>0$ depending only on $\rho$, $\mathcal A$, $\Omega$, and $T$. Next we define the map
\begin{equation*}
\mathcal G:L^2(0,T;L^2(\Omega))\ni F\mapsto \sum_{k=1}^N\int_0^T\lambda_k(t)w_F(x_k,T-t)\,\d t.
\end{equation*}
By the Cauchy-Schwarz inequality, the embedding $H^2(\Omega)\hookrightarrow C(\overline{\Omega})$ \cite[Theorem 4.12]{Adams:2023} and the estimate \eqref{p1b}, for any $F\in L^2(0,T;L^2(\Omega))$, we have
\begin{equation*}
|\mathcal{G}(F)|
\leq \zeta \norm{w_F}_{L^2(0,T;C(\overline{\Omega}))}
\leq c\zeta\norm{F}_{L^2(0,T;L^2(\Omega))},
\end{equation*}
with $c>0$ depending only on $\rho$, $\mathcal A$, $\Omega$ and $T$. Thus, $\mathcal G$ is a continuous linear form on $L^2(0,T;L^2(\Omega))$, and further, there holds
$\norm{\mathcal G}_{L^2(0,T;L^2(\Omega))}\leq
c\zeta .$
By Riesz representation theorem, there exists a unique $v\in L^2(0,T;L^2(\Omega))$ such that
\begin{equation*}
\left\langle v,F\right\rangle_{L^2(0,T;L^2(\Omega))}=\mathcal G(F)= \sum_{k=1}^N\int_0^T\lambda_k(t)w_F(x_k,T-t)\,\d t,\quad F\in L^2(0,T;L^2(\Omega)).
\end{equation*}
Then $v$ is the unique solution in the transposition sense of problem \eqref{eq2}. Moreover, the estimate \eqref{p1a} follows directly. To prove \eqref{p1aa}, we fix $\psi\in C^\infty_0(Q)$ and  define $\tilde{\psi}$ by $\tilde{\psi}(x,t)=\psi(x,T-t)$, $(x,t)\in \Omega\times (0,T)$. By the definition of solution $v$, we have
\begin{equation}\label{p1d}
\int_0^T\int_\Omega v(x,t)[\rho(x) D_t^\alpha \tilde{\psi}(x,T-t)+\mathcal A^* \tilde{\psi}(x,T-t)]\,\d x\d t=\sum_{k=1}^N\int_0^T \lambda_k(t)\psi(x_k,t)\,\d t.
\end{equation}
From the condition $\tilde{\psi}(x,0)=0$, $x\in\Omega$, we obtain
$$D_t^\alpha \tilde{\psi}(x,T-t)=\int_0^{T-t}\frac{(T-t-s)^{-\alpha}}{\Gamma(1-\alpha)}\partial_s \tilde{\psi}(x,s)\,\d s.$$
It follows directly from this identity and Fubini's theorem that
\begin{align*}
\int_0^T\int_\Omega v(x,t) D_t^\alpha \tilde{\psi}(x,T-t)\,\d x\d t&=\int_{\Omega}\int_0^Tv(x,t)\int_0^{T-t}\frac{(T-t-s)^{-\alpha}}{\Gamma(1-\alpha)}\partial_s \tilde{\psi}(x,s)\,\d s \d t \d x\\
=&\int_{\Omega}\int_0^Tv(x,t)\int_t^{T}\frac{(s-t)^{-\alpha}}{\Gamma(1-\alpha)}\partial_s \tilde{\psi}(x,T-s)\,\d s \d t \d x\\
=&-\int_{\Omega}\int_0^T\left(\int_0^{s}\frac{(s-t)^{-\alpha}}{\Gamma(1-\alpha)}v(x,t)\,\d t\right)\partial_s [\tilde{\psi}(x,T-s)]\d s  \d x\\
=&-\left\langle {_0I_t^{1-\alpha}}v,\partial_s \psi\right\rangle_{D'(Q),C^\infty_0(Q)}=\left\langle D^\alpha_tv, \psi\right\rangle_{D'(Q),C^\infty_0(Q)}.
\end{align*}
This and the identity \eqref{p1d} imply
\begin{align*} 
\left\langle \rho(x) D_t^\alpha v+\mathcal A v,\psi\right\rangle_{D'(Q),C^\infty_0(Q)}&= \int_0^T\int_\Omega v(x,t)[\rho(x) D_t^\alpha \tilde{\psi}(x,T-t)+\mathcal A^* \tilde{\psi}(x,T-t)]\,\d x\d t\\
&=\sum_{k=1}^N\int_0^T\lambda_k(t)\psi(x_k,t)\,\d t,\end{align*}
which proves the identity \eqref{p1d}, and also the theorem.
\end{proof}

By \cite[Theorem 2.1]{KRY}, there exists a unique function $v_1\in H^\alpha(0,T;L^2(\Omega))\cap L^2(0,T;H^2(\Omega)\cap H^1_0(\Omega))$ such that $v_1-u_0\in H_\alpha(0,T;L^2(\Omega))$ and
$$\rho(x) D_t^\alpha [v_1-u_0](x,t)+\mathcal A v_1(x,t)=0,\quad (x,t)\in \Omega\times(0,T).$$
Then, we can define the unique solution $u$ in the transposition sense of problem \eqref{eq1} in such a way that $u-v_1$ is the unique solution in the transposition sense of problem \eqref{eq2}.

The definition of solutions in the transposition sense extends easily to problem \eqref{eqq11}. In addition, since the space $H^1_0(0,\ell)$ embeds continuously into $C[0,\ell]$, we have $\delta_{x_0}\in H^{-1}(0,\ell)$. By \cite[Theorem 4.1]{KRY}, for $u_0\in L^2(0,\ell)$ and $\lambda\in L^2(0,T)$, there exists a unique $u_\star\in L^2(0,T;H^1_0(0,\ell))$ satisfying $u_\star-u_0\in H_\alpha(0,T;H^{-1}(0,\ell))$ such that for a.e. $t\in(0,T)$ and all $\psi\in H^1_0(0,\ell)$, we have
\begin{equation*}
\left\langle \rho D_t^\alpha(u_\star-u_0)(\cdot,t)+\mathcal A(t) u_\star(\cdot,t),\psi\right\rangle_{H^{-1}(0,\ell),H^1_0(0,\ell)}=\lambda(t)\psi(x_0).
\end{equation*}
Similar to Proposition \ref{p1}, $u_\star$ coincides with the solution $u$ in the transposition sense  of \eqref{eqq11} and $u\in L^2(0,T;H^1_0(0,\ell))$ satisfies $u-u_0\in H_\alpha(0,T;H^{-1}(0,\ell))$.

\subsection{Improved regularity}

In Section \ref{ssec:transposition}, we have shown that problem \eqref{eq1} has a unique solution $u\in L^2(0,T;L^2(\Omega))$. Since the singularities of the source $\sum_{k=1}^N\lambda_k(t)\delta_{x_k}$ are located at $N$ points $\{x_k\}_{k=1}^N\subset \Omega$, we can improve the regularity on $\Omega\setminus \{x_k\}_{k=1}^N$. Below $B(y,r):=\{x\in\R^d:\ |x-y|<r\}$ denotes the ball centered at $y$ with a radius $r$.
\begin{theorem}\label{t2}
Let $\{\lambda_k\}_{k=1}^N\subset H_\alpha(0,T)$. Then for the unique solution $u$  of problem \eqref{eq1} in the transposition sense, we have $u-u_0\in H_\alpha(0,T;L^2(\Omega))$ and, for all $r>0$, by fixing
$\Omega_r:=\Omega\setminus\cup_{k=1}^N \overline{B(x_k,r)}$,
we have $u\in L^2\left(0,T; H^2(\Omega_r)\right)$ and
\bel{t2aa}\rho(x) D_t^\alpha(u-u_0)(x,t)+\mathcal A u(x,t)=0,\quad (x,t)\in \Omega_r\times(0,T).\ee
\end{theorem}

To prove this result, consider $R>0$ such that $\overline{\Omega}\subset B_R:=\{x\in\R^d:\ |x|<R\}$, fix $R_1>R$ and define the operator $\mathcal A$ acting on $B_{R_1}$ by extending the coefficients $b_1,\ldots,b_k$, $q$ into functions lying in $C^1(\overline{B_{R_1}})$. By the theory of fundamental solutions for elliptic equations \cite[Theorem 3]{Ka}, there exists a function $P$ defined in $B_{R_1}^2\setminus D$, with $D:=\{(x,y)\in B_{R_1}^2:\ x=y\}$, such that
\begin{itemize}
\item[(1)] For all $y\in B_{R_1}$, the map $x\mapsto P(x,y)\in C^2(B_{R_1}\setminus \{y\})$ and
\begin{equation}\label{fond1}
|P(x,y)|\leq c\max\left(|x-y|^{2-d},\ln(|x-y|)\right),\quad \mathcal A P(x,y)=0,\quad x\in B_{R_1}\setminus \{y\};
\end{equation}
\item[(2)] For all $\psi \in C^\infty_0(B_{R_1})$ we have
\begin{equation}\label{fond2}
\int_{B_{R_1}}P(x,y)\mathcal A^*\psi(x)\,\d x=\psi(y).
\end{equation}
\end{itemize}

Note that since $d\leq3$, the estimate \eqref{fond1} implies that for any $y\in B_{R}$, $x\mapsto P(x,y)\in L^2(B_{R})$.
Moreover, we have the following result.
\begin{lemma}\label{l1} For any $v\in H^2(\Omega)\cap H^1_0(\Omega)$ and $y\in\Omega$, we have
\begin{equation}\label{l1a}
\int_\Omega P(x,y)\mathcal A^* v(x)\,\d x=v(y)-\int_{\partial\Omega}P(x,y)\partial_{\nu_a}v(x)\,\d\sigma(x),
\end{equation}
where $\partial_{\nu_a}v(x)=\sum_{i,j=1}^da_{ij}(x)\partial_{x_j}v(x)\nu_i(x)$, for $x\in \partial\Omega$, is the conormal derivative associated with the coefficient $a$.
\end{lemma}
\begin{proof}
First, by density, the identity \eqref{fond2} holds also for $\psi\in H^2_0(B_{R_1})$, the closure of $C^\infty_0(B_{R_1})$ in $H^2(B_{R_1})$. To prove \eqref{l1a}, we first extend $v$ into an element of $H^2_0(B_{R_1})$ satisfying $v=0$ on $B_{R_1}\setminus B_R$, still denoted by $v$. Then the identity \eqref{fond2} implies
\begin{align*}
  \int_{B_{R}}P(x,y)\mathcal A^*v(x)\,\d x=\int_{B_{R_1}}P(x,y)\mathcal A^*v(x)\,\d x=v(y).
\end{align*}
For any fixed $y\in\Omega$, we have $x\mapsto P(x,y)\in C^2(\overline{B_{R}}\setminus \Omega)$ and, by integration by parts and using \eqref{fond1}, we get
\begin{align*}
\int_{B_{R}\setminus\Omega}P(x,y)\mathcal A^*v(x)\,\d x&=\int_{\partial\Omega}P(x,y)\partial_{\nu_a}v(x)\,\d\sigma(x)+\int_{B_{R}\setminus\Omega}\mathcal A P(x,y)v(x)\,\d x\\
&=\int_{\partial\Omega}P(x,y)\partial_{\nu_a}v(x)\,\d\sigma(x).\end{align*}
Then it follows that
\begin{align*}\int_{\Omega}P(x,y)\mathcal A^*v(x)\,\d x&=\int_{B_{R}}P(x,y)\mathcal A^*v(x)\,\d x-\int_{B_{R}\setminus\Omega}P(x,y)\mathcal A^*v(x)\,\d x\\
&=v(y)-\int_{\partial\Omega}P(x,y)\partial_{\nu_a}v(x)\,\d\sigma(x),
\end{align*}
which proves the identity \eqref{l1a}, and also completes the proof of the lemma.
\end{proof}

For any fixed $y\in \Omega$, we have $\partial\Omega\ni x\mapsto P(x,y)\in C^2(\partial\Omega)\subset H^{\frac{3}{2}}(\partial\Omega)$, and by \cite[Theorems 8.3 and 8.12]{GT}, we can consider the $H^2(\Omega)$-solution $v_y$ of the boundary value problem
\begin{equation}\label{eq3}
\left\{\begin{aligned}
\mathcal{A} v_y &= 0, && \mbox{in }\Omega,\\
 v_y(x)&= -P(x,y), && x\in\partial\Omega.
\end{aligned}\right.
\end{equation}
Then Lemma \ref{l1} and integration by parts yield that $\mathcal W_y\in L^2(\Omega)$ given by
\begin{align}\label{w}
\mathcal W_y(x)=P(x,y)+v_y(x),\quad x\in\Omega
\end{align}
satisfies the condition
$\int_{\Omega}\mathcal W_y \mathcal A^*h\,\d x=h(y)$ for $h\in H^2(\Omega)\cap H^1_0(\Omega)$.

In light of \cite[Theorem 2.3]{KRY}, for each $k=1,\ldots,N$, since $\lambda_k\in H_\alpha(0,T)$, we have $_0I_t^{1-\alpha}\lambda_k\in H_1(0,T)$ and $D_t^\alpha \lambda_k=\partial_t{_0I_t^{1-\alpha}}\lambda_k\in L^2(0,T)$. Thus, we  can define $\mathcal R$ to be the unique element of $H_\alpha(0,T;L^2(\Omega))\cap L^2(0,T;H^2(\Omega)\cap H^1_0(\Omega))$ satisfying
\begin{align}\label{R}
D_t^\alpha \mathcal R+\mathcal A \mathcal R=-\rho(x)\sum_{k=1}^ND_t^\alpha \lambda_k(t) \mathcal W_{x_k}(x).
\end{align}

Using these tools, we can prove Theorem \ref{t2}.
\begin{proof}[Proof of Theorem \ref{t2}] Without loss of generality, we may assume $u_0\equiv0$. Consider the function $u_1\in L^2(0,T;L^2(\Omega))$ defined by
\begin{equation*}
u_1(x,t)=\sum_{k=1}^N \lambda_k(t) \mathcal W_{x_k}(x)+\mathcal R(x,t),\quad (x,t)\in \Omega\times(0,T),
\end{equation*}
with $\mathcal W_{x_k}$ given by \eqref{w} with $y=x_k$. The preceding discussion shows $ u_1\in H_\alpha(0,T;L^2(\Omega))\cap L^2\left(0,T; H^2(\Omega_r)\right)$. It suffices to show that $u_1$ is a solution in the transposition sense of \eqref{eq1}, and by density, that the identity
\begin{align}\label{t2d}
\left\langle u_1,F\right\rangle_{L^2(0,T;L^2(\Omega))}=\sum_{k=1}^N\int_0^T\lambda_k(t)w_F(x_k,T-t)\,\d t
\end{align}
holds for any $F\in C^\infty_0(Q)$, with $w_F\in H_\alpha(0,T;L^2(\Omega))\cap L^2(0,T;H^2(\Omega)\cap H^1_0(\Omega))$ satisfying
\begin{align*}
\rho(x) D_t^\alpha w_F(x,t)+\mathcal A^* w_F(x,t)=F(x,T-t),\quad (x,t)\in \Omega\times(0,T).
\end{align*}
Since $F\in C^\infty_0(Q)$, we have $\partial_tF\in L^2(0,T;L^2(\Omega))$ and can define $w_{\partial_tF}\in H_\alpha(0,T;L^2(\Omega))\cap L^2(0,T;H^2(\Omega)\cap H^1_0(\Omega))$, cf. \eqref{eqn:w-F}. Since $\rho$ and $\mathcal A^*$ are independent of $t$, direct computation gives
$w_F(x,t)=-\int_0^tw_{\partial_tF}(x,s)\,\d s$ for $(x,t)\in \Omega\times(0,T)$, which implies
$w_F\in H_1(0,T;L^2(\Omega))\cap L^2(0,T;H^2(\Omega)\cap H^1_0(\Omega))$. Then, similar to the proof of Proposition \ref{p1}, we can prove
\begin{align*}
\int_0^T\int_\Omega u_1(x,t) D_t^\alpha w_F(x,T-t)\,\d x\d t=-\int_{\Omega}\int_0^TI^{1-\alpha}u_1(x,t)\partial_t [w_F(x,T-t)]\,\d t  \d x.
\end{align*}
Further, the condition $u_1\in H_\alpha(0,T;L^2(\Omega))$ implies $_0I_t^{1-\alpha}u_1\in H_1(0,T;L^2(\Omega))$ \cite[Theorem 2.3]{KRY}. Then by integration by parts, we get
\begin{equation}\label{t2e}\begin{aligned}
&\quad \int_0^T\int_\Omega u_1(x,t) D_t^\alpha w_F(x,T-t)\,\d x\d t=\int_{\Omega}\int_0^T\partial_t I^{1-\alpha}u_1(x,t) w_F(x,T-t)\,\d t \d x\\
&=\int_0^T\int_\Omega D_t^\alpha u_1(x,t)  w_F(x,T-t)\,\d x\d t\\
&=\int_0^T\int_\Omega \Big[D_t^\alpha \mathcal R+\sum_{k=1}^ND_t^\alpha \lambda_k(t) \mathcal W_{x_k}(x)\Big]w_F(x,T-t)\,\d x \d t.
\end{aligned}
\end{equation}
Moreover, by Lemma \ref{l1} and integration by parts, we get
\begin{align*}
&\int_0^T\int_\Omega\Big[\sum_{k=1}^N \lambda_k(t) \mathcal W_{x_k}(x)\Big]\mathcal A^* w_F(x,T-t)\,\d x \d t\\
=&\sum_{k=1}^N\int_0^T\lambda_k(t)\left(\int_\Omega [P(x,x_k)+v_{x_k}(x)]\mathcal A^* w_F(x,T-t)\,\d x\right)\,\d t\\
=&\sum_{k=1}^N\int_0^T\lambda_k(t) \left(w_F(x_k,T-t) -\int_{\partial\Omega}[P(x,x_k)+v_{x_k}(x)]\partial_{\nu_a}w_F(x,T-t)\,\d\sigma(x)\right) \,\d t\\
=&\sum_{k=1}^N\int_0^T\lambda_k(t) w_F(x_k,T-t)\,\d t,
\end{align*}
since $P(x,x_k)+v_{x_k}(x)$ vanishes on the boundary $\partial\Omega$.
Thus, it follows from the last identity that
\begin{align*}\left\langle u_1,F\right\rangle_{L^2(0,T;L^2(\Omega))}
=&\int_0^T\int_\Omega u_1(x,t)[\rho(x) D_t^\alpha w_F(x,T-t)+\mathcal A^* w_F(x,T-t)]\,\d x \d t\\
=&\int_0^T\int_\Omega\Big[\rho D_t^\alpha \mathcal R+\mathcal A\mathcal R+ \rho(x)\sum_{k=1}^ND_t^\alpha \lambda_k(t) \mathcal W_{x_k}(x)\Big]w_F(x,T-t)\,\d x \d t\\ &+\sum_{k=1}^N\int_0^T\lambda_k(t)w_F(x_k,T-t)\,\d t
=\sum_{k=1}^N\int_0^T\lambda_k(t)w_F(x_k,T-t)\,\d t.
\end{align*}
Therefore, the identity \eqref{t2d} holds and  $u=u_1$ is the unique solution in the transposition sense of problem \eqref{eq1}. The definition of $u_1$ indicates that $u=u_1$ fulfills \eqref{t2aa}.
\end{proof}

We have the following improved regularity in the 1D case. For any $x_0\in(0,\ell)$ and $r\in (0,\min(x_0,\ell-x_0))$, we denote $I_{x_0,r}^-=(0,x_0-r)$ and $I_{x_0,r}^+=(x_0+r,\ell)$.

\begin{theorem}\label{tt2} Let  $u_0\in H^1_0(0,\ell)$. Then for the unique solution $u\in L^2(0,T;H^1_0(0,\ell))$ of problem \eqref{eqq11} and all $r\in(0,\min(x_0,\ell-x_0))$,   we have
\begin{align*}
u|_{I_{x_0,r}^-\times(0,T)}-u_0|_{I_{x_0,r}^-}&\in H_\alpha(0,T;L^2(I_{x_0,r}^-)),\quad  u|_{I_{x_0,r}^-\times(0,T)}\in L^2(0,T;H^2(I_{x_0,r}^-)),\\
u|_{I_{x_0,r}^+\times(0,T)}-u_0|_{I_{x_0,r}^+}&\in H_\alpha(0,T;L^2(I_{x_0,r}^+)), \quad  u|_{I_{x_0,r}^+\times(0,T)}\in L^2(0,T;H^2(I_{x_0,r}^+)).
\end{align*}
\end{theorem}
\begin{proof}
Fix $\chi\in C^\infty[0,\ell]$ satisfying $\chi(x_0)=0$. Then, by choosing $v(x,t)=\chi(x)u(x,t)$ for $(x,t)\in (0,\ell)\times(0,T)$, we derive that $v\in L^2(0,T;H^1_0(0,\ell))$ satisfies $v-\chi u_0\in H_\alpha(0,T;H^{-1}(0,\ell))$ and, for all $\psi\in H^1_0(0,\ell)$ and a.e. $t\in(0,T)$,  the condition
\begin{equation}\label{eqq12}
\left\langle \rho(x)D_t^{\alpha}(v-\chi u_0)(\cdot,t) +\mathcal{A}(t) v(\cdot,t),\psi\right\rangle_{H^{-1}(0,\ell),H^1_0(0,\ell)}=\left\langle G(\cdot,t),\psi\right\rangle_{H^{-1}(0,\ell),H^1_0(0,\ell)},
\end{equation}
where the function $G$ is defined by
\begin{align*}
G(x,t)&=\chi(x)(\rho(x)\partial_t^{\alpha}u +\mathcal{A}(t) u)-2a(x,t)\chi'(x)\partial_xu(x,t)+b(x,t)\chi'(x)u(x,t)\\
&=\underbrace{\chi(x_0)}_{=0}\lambda(t)\delta_{x_0}(x)-2a(x,t)\chi'(x)\partial_xu(x,t)+b(x,t)\chi'(x)u(x,t).
\end{align*}
Since $u\in L^2(0,T;H^1_0(0,\ell))$, we have $G\in L^2(0,T;L^2(0,\ell))$. This, the fact $\chi u_0\in H^1_0(0,\ell)$ and \cite[Theorem 2]{KuY} yield $v\in L^2(0,T;H^2(0,\ell))$ and $v-\chi(x)u_0(x)\in H_\alpha(0,T;L^2(0,\ell))$. Thus, we have $\chi (x) u\in L^2(0,T;H^2(0,\ell))$, $\chi(x)(u-u_0)\in H_\alpha(0,T;L^2(0,\ell))$ and, since $\chi\in C^\infty[0,\ell]$ is only subject to the condition $\chi(x_0)=0$, we deduce the desired result.
\end{proof}

\section{The proofs of the main theorems}\label{sec:proof}
\subsection{The proof of Theorem \ref{t1}}
One key ingredient of the proof of Theorem \ref{t1} is the classical unique continuation principle of \cite[Theorem 1.2]{LN} that we adapt slightly to solutions of problem \eqref{l2a}. For any Banach space $X$, we denote by $D'_+(\R;X)$ the set of distributions taking values in $X$ and supported on $[0,\infty)$.
\begin{lemma}\label{l2} Let $\mathcal O$ be an open connected subset of $\mathbb{R}^d$ {\rm(}with $d=2,3${\rm)} with a $C^2$ boundary and let $w_0\in H^1(\mathcal O)$, let $w\in  H^\alpha(0,T;L^2(\mathcal O))\cap L^2(0,T;H^2(\mathcal O))$ be such that $w-w_0\in H_\alpha(0,T;L^2(\mathcal O))$ and
\begin{equation}\label{l2a}
\rho(x) D_t^\alpha [w-w_0(x)](x,t)+\mathcal A w(x,t)=0,\quad (x,t)\in \mathcal O\times(0,T).
\end{equation}
Suppose also that there exists an open nonempty subset $U$ of $\mathcal O$ such that
$w(x,t)=0$ for $(x,t)\in U\times(0,T).$
Then, we have $w\equiv 0$ in $\mathcal O\times(0,T)$.
\end{lemma}
\begin{proof}
In view of \cite{KS} (see also \cite[Chap. 6, Section 5, (15)]{S}), the map $_0I_t^{1-\alpha}$ can be extended to a continuous linear map from $D'_+(\R;L^2(\mathcal O))$ to $D'_+(\R;L^2(\mathcal O))$ by setting, for all $v \in D_+'(\R,L^2(\Omega))$ and all $\psi\in D_{\mathrm a}(\R,L^2(\Omega)):=\{\psi \in C^\infty(\R,L^2(\Omega));\ \exists R>0,\ \textrm{supp} (\psi) \subset(-\infty,R)\}$,
\begin{equation*}
\langle {_0I_t^{1-\alpha}}v,\psi \rangle_{D'_+(\R,L^2(\Omega)), D_{\mathrm a}(\R,L^2(\Omega))} := \left \langle v(s,\cdot),\int_0^{+\infty}\!\!\!\frac{t^{-\alpha}}{\Gamma(1-\alpha)}\psi(t+s,\cdot)\,\d t \right \rangle_{D'_+(\R,L^2(\Omega), D_{\mathrm a}(\R,L^2(\Omega))}.
\end{equation*}
For any $v\in D'_+(\R;L^2(\mathcal O))$, we have
$\partial_t{_0I_t^{1-\alpha}v}={_0I_t^{1-\alpha}\partial_tv}$, and can study the Caputo derivative in the sense of distribution, i.e., $\partial_t^\alpha v={_0I_t^{1-\alpha}}\partial_tv$ for $v\in D'_+(\R;L^2(\mathcal O))$. This definition uses convolution of elements of $D'_+(\R)$ and coincides with that in \cite[Theorem 1.2]{LN}. By symmetry, we can extend $w$ to $w_*\in H^\alpha(\R_+;L^2(\mathcal O))\cap L^2(\R_+;H^2(\mathcal O))$ and then extend $w_*$ to $\tilde{w}$ defined on $\mathcal O\times\R$ by $\tilde{w}(x,t)=w_0(x)$, $(x,t)\in\mathcal O\times(-\infty,0]$. Since for all $T_1>0$, $w_*-w_0\in H_\alpha(0,T_1;L^2(\Omega))$, by \cite[Theorems 11.4 and 5.1, Chap. 1]{LM1}, the map $\tilde{w}-w_0$ belongs to $H^\alpha_{loc}(\R;L^2(\mathcal O))\cap D'_+(\R;L^2(\mathcal O))$. Thus \cite[Theorem 2.3]{KRY} implies $_0I_t^{1-\alpha}(\tilde{w}-w_0)\in H^1_{loc}(\R;L^2(\mathcal O))$ and
\begin{align*}
_0I_t^{1-\alpha}\partial_t(\tilde{w}-w_0)(x,t)&=\partial_t{_0I_t^{1-\alpha}}(\tilde{w}-w_0)(x,t)\\
&=\partial_t\left(\int_0^t\frac{(t-s)^{-\alpha}}{\Gamma(1-\alpha)}(\tilde{w}(x,s)-w_0)\,\d s\right),\quad (x,t)\in \mathcal O\times\R.
\end{align*}
Meanwhile, using the identity $\partial_t(\tilde{w}-w_0)=\partial_t\tilde{w}$ in $ \mathcal O\times\R$, the condition $\tilde{w}=w$ in $\mathcal O\times(0,T)$ (from the construction of $\tilde w$) and the identity \eqref{l2a} yield
\begin{equation*}
\partial_t^\alpha\tilde{w}(x,t)+\rho^{-1}(x)\mathcal A \tilde{w}(x,t)=D_t^\alpha [w-w_0](x,t)+\rho^{-1}(x)\mathcal A w(x,t)=0,\quad (x,t)\in \mathcal O\times(0,T).
\end{equation*}
This, the condition
$\tilde{w}(x,t)=w(x,t)=0$ for $ (x,t)\in U\times(0,T)$
and \cite[Theorem 1.2]{LN} yield
$$w(x,t)=\tilde{w}(x,t)=0,\quad (x,t)\in\mathcal O\times(0,T),$$
which completes the proof of the lemma.
\end{proof}

Using Theorem \ref{t2} and Lemma \ref{l2}, we can now prove Theorem \ref{t1}.
\begin{proof}
The proof of the theorem is lengthy and divided into four steps below.

\smallskip
\noindent\textbf{Step 1:} This step proves that condition \eqref{t1a} implies
\begin{equation}\label{t1c0} u_1(x,t)=u_2(x,t),\quad (x,t)\in\omega\times (0,T).
\end{equation}
Let $u=u_1-u_2$. Then $u$ solves in the sense of transposition
\begin{equation}\label{eq7}
\left\{\begin{aligned}
\rho(x)\partial_t^{\alpha}u +\mathcal{A} u &=  F(x,t), && \mbox{in }\Omega\times(0,T),\\
 u&= 0, && \mbox{on } \partial\Omega\times(0,T), \\
u&=u_0, && \mbox{in } \Omega\times \{0\},
\end{aligned}\right.
\end{equation}
with $u_0=u_0^1-u_0^2$ and the source $F$ given by
$$F(x,t)=\sum_{k=1}^{N_1}\lambda_k^1(t)\delta_{x_k^1}(x)-\sum_{k=1}^{N_2}\lambda_k^2(t)\delta_{x_k^2}(x),\quad (x,t)\in Q.$$
For $r>0$, let $\Omega_r=\Omega\setminus(\cup_{j=1}^2\cup_{k=1}^{N_j} \overline{B(x_k^j,r)})$,
and choose $r>0$ sufficiently small such that $\omega\subset\Omega_r$. By Theorem \ref{t2}, $u\in H_\alpha(0,T;L^2(\Omega))\cap L^2(0,T;H^2(\Omega_r))$ and satisfies
\begin{align}\label{t1e}
\rho(x)D_t^\alpha(u-u_0)(x,t)+\mathcal A u(x,t)=0,\quad (x,t)\in \Omega_r\times (0,T).
\end{align}
Now fix any $\psi\in C^\infty_0(\omega)$, and let $k_\psi(t)=\left\langle u(\cdot,t),\psi\right\rangle_{L^2(\Omega)}$ and  $\mu_\psi=\left\langle u_0,\psi\right\rangle_{L^2(\Omega)}$. By repeating the preceding argument, we have $k_\psi-\mu_\psi\in H_\alpha(0,T)$ with
\begin{align*}
D_t^\alpha(k_\psi-\mu_\psi)(t)+\left\langle \rho^{-1}\mathcal A u(\cdot,t),\psi\right\rangle_{L^2(\Omega)}=0.
\end{align*}
By condition \eqref{t1a}, we have $u=0$ in $\omega\times(T-\epsilon,T)$, and also  $\rho^{-1}\mathcal A u=0$ in $ \omega\times(T-\epsilon,T)$. Thus, 
\begin{align}\label{t1f}
k_\psi(t)=D_t^\alpha(k_\psi-\mu_\psi)(t)=0,\quad t\in (T-\epsilon,T).
\end{align}
Now we extend $k_\psi$ by zero to $\mathbb{R}_+$, still denoted by $k_\psi$. In view of \eqref{t1f}, we have that, for all $T_1>0$, $k_\psi-\mu_\psi\in H_\alpha(0,T_1)$. Thus, we can also define $h_\psi$ on $\R_+$  by $ h_\psi(t)=D_t^\alpha(k_\psi-\mu_\psi)(t)$,
and note that for all $T_1>0$, $h_\psi\in L^2(0,T_1)$. Since $u\in H^\alpha(0,T;L^2(\Omega))$, we deduce  $k_\psi\in L^1(\R_+)\cap L^2(\R_+)$ and, since $k_\psi=0$ on $(T-\epsilon,+\infty)$, we obtain
$$_0I_t^{1-\alpha}k_\psi(t)=\int_0^{T-\epsilon}\frac{(t-s)^{-\alpha}}{\Gamma(1-\alpha)}k_\psi(s)\,\d s,\quad t\in(T-\epsilon,T).$$
It follows that
\begin{align*}
h_\psi(t)=\int_0^{T-\epsilon}\frac{-\alpha(t-s)^{-1-\alpha}}{\Gamma(1-\alpha)}k_\psi(s)\,\d s-\frac{t^{-\alpha}}{\Gamma(1-\alpha)}\mu_\psi,\quad t\in(T-\epsilon,\infty),
\end{align*}
and hence the map $h_\psi(t)$ is analytic with respect to $t\in(T-\epsilon,+\infty)$. Moreover, condition \eqref{t1f} implies  $h_\psi=0$ on $(T-\epsilon,T)$ and by the isolated zero theorem for analytic functions, we have $h_\psi=0$ on $(T-\epsilon,\infty)$. Thus, the Laplace transform $\widehat{h}_\psi$ of $h_\psi$ is well defined and holomorphic in the whole complex plane $\mathbb{C}$. Moreover, since $k_\psi\in L^1(\R_+)\cap L^2(\R_+)$, we deduce
\begin{align*}
\widehat{h}_\psi(p)=p^\alpha(\widehat{k}_\psi(p)-p^{-1}\mu_\psi),\quad p\in\mathbb C_+:=\{z\in\mathbb C:\ \Re(z)>0\}.
\end{align*}
Since the maps $\widehat{h}_\psi(p)$ and $p^\alpha(\widehat{k}_\psi(p)-p^{-1}\mu_\psi)$ are holomorphic with respect to $p\in\mathbb C\setminus(-\infty,0]$, we can extend the identity as
\begin{align*}
\widehat{h}_\psi(p)=p^\alpha(\widehat{k}_\psi(p)-p^{-1}\mu_\psi),\quad p\in\mathbb C\setminus(-\infty,0].
\end{align*}
By fixing $r>0$, $\theta\in(0,\pi)$ and setting $p=re^{\pm i\theta}$ in the identity, we get
\begin{align*}
\widehat{h}_\psi(re^{ i\theta})-\widehat{h}_\psi(re^{ -i\theta})=r^\alpha e^{ i\alpha\theta}(\widehat{k}_\psi(re^{ i\theta})-r^{-1}e^{-i\theta}\mu_\psi)-r^\alpha e^{ -i\alpha\theta}(\widehat{k}_\psi(re^{ -i\theta})-r^{-1}e^{i\theta}\mu_\psi),\quad r>0.
\end{align*}
By sending $\theta\to\pi$ and the continuity of $\widehat{h}_\psi$ and $\widehat{k}_\psi$ in $\mathbb C$, we get
\begin{align*}
0=\widehat{h}_\psi(-r)-\widehat{h}_\psi(-r)= 2i \sin(\alpha\pi) r^\alpha(\widehat{k}_\psi(-r)+r^{-1}\mu_\psi),\quad r>0.
\end{align*}
Thus,
$\widehat{k}_\psi(-r)=-r^{-1}\mu_\psi$ for $r>0$,
and the unique continuation of holomorphic functions implies
$\widehat{k}_\psi(z)=z^{-1}\mu_\psi$ for $z\in\mathbb C\setminus\{0\}$. Since $\widehat{k}_\psi$ is holomorphic in $\mathbb{C}$, this implies  $\mu_\psi=0$ and similarly also $\widehat{k}_\psi\equiv0$. By the injectivity of Laplace transform, we deduce
$\left\langle u(\cdot,t),\psi\right\rangle_{L^2(\Omega)}=k_\psi(t)=0$ for $t\in(0,T)$. Since $\psi\in C^\infty_0(\omega)$ is arbitrarily chosen, we have $u=0$ on $\omega\times(0,T)$, which implies \eqref{t1c0}.

\medskip\noindent
\textbf{Step 2:} This step proves that condition \eqref{t1c0} implies
\bel{t1h} u_1(x,t)=u_2(x,t),\quad (x,t)\in\Omega\times (0,T).\ee
By fixing $r>0$ sufficiently small, $\Omega_r$ is a connected open set with a $C^2$ boundary and $\omega\subset\Omega_r$. Moreover, by Theorem \ref{t2}, $u:=u_1-u_2\in H^\alpha(0,T;L^2(\Omega_r))\cap L^2(0,T;H^2(\Omega_r))$, $u-u_0\in H_\alpha(0,T;L^2(\Omega_r))$ and
\bel{t1i}\rho(x) D_t^\alpha [u-u_0(x)](x,t)+\mathcal A u(x,t)=0,\quad (x,t)\in \Omega_r\times(0,T).\ee
This, \eqref{t1c0} and Lemma \ref{l2} yield
$u=0$ in $ \Omega_r\times(0,T).$
Since it holds for any $r>0$ arbitrarily small and $u\in L^2(Q)$, $u\equiv0$ as a function lying in $L^2(0,T;L^2(\Omega))$ and \eqref{t1h} holds.

\medskip\noindent\textbf{Step 3:} This step proves  $u_0^1=u_0^2$. Combining \eqref{t1h} and \eqref{t1i} yields that for  $r>0$ small enough
$$-\rho(x)\frac{t^{-\alpha}}{\Gamma(1-\alpha)}u_0(x)=\rho(x) D_t^\alpha [u-u_0(x)](x,t)+\mathcal A u(x,t)=0,\quad (x,t)\in \Omega_r\times(0,T),$$
which implies $u_0=0$ in $\Omega_r$ for $r>0$ arbitrarily small. Thus $u_0\equiv0$ and $u_0^1=u_0^2$.

\medskip
\noindent\textbf{Step 4:} This step completes the proof. From Step 3, $u$ solves problem \eqref{eq7} in the sense of transposition with $u_0\equiv0$. That is, for any $\psi\in C^\infty_0(Q)$, by letting $\tilde{\psi}(x,t)=\psi(x,T-t)$, $(x,t)\in Q$, we get
$$\int_0^T\int_\Omega u(x,t)[\rho(x) D_t^\alpha \tilde{\psi}(x,T-t)+\mathcal A^* \tilde{\psi}(x,T-t)]\,\d x\d t=\left\langle F,\psi\right\rangle_{D'(Q),C^\infty_0(Q)}.$$
This and \eqref{t1h} give
$\left\langle F,\psi\right\rangle_{D'(Q),C^\infty_0(Q)}=0$ for any $\psi\in C^\infty_0(Q).$
That is, the distribution $F$ is equal to zero and
\bel{t1ab}\sum_{k=1}^{N_1}\lambda_k^1(t)\delta_{x_k^1}(x)=\sum_{k=1}^{N_2}\lambda_k^2(t)\delta_{x_k^2}(x),\quad (x,t)\in \Omega\times(0,T).\ee
Since $\lambda_1^j,\ldots,\lambda_{N_j}^j$, $j=1,2$, are non-uniformly vanishing, condition \eqref{t1ab}
implies  $N_1=N_2=N$ and that \eqref{t1c} holds. This completes the proof of the theorem.
\end{proof}

\begin{remark}\label{rmk:memory}
\textit{A posteriori} measurement was first studied in \cite{JK} for an inverse source problem in subdiffusion; see \cite{KLY,Y1} for related works. In the ODE case, the result asserts: If $\partial_t^\alpha v (t) = v(t) = 0$ for $t>t_0$, with $t_0>0$, then $v(t)=0$ for $t>0$ \cite[Theorem 1]{JK}. Indeed, for $v\in L^1(0,t_0)$, the fractional integral ${_0I_t^{1-\alpha}} v(t)$ is analytic over $(t_0,\infty)$ \cite[Lemma 1]{JK}. This work exploits the singularity of Laplace transform, and the result shows the beneficial memory effect of subdiffusion.
\end{remark}

\subsection{The proof of Theorem \ref{t1D}}

Now we discuss the one-dimensional case. First we give a counterexample for uniqueness.
\begin{example}\label{exam:1d}
There is one point source at $x_0\in\Omega=(0,1)$ with $\lambda(t)\equiv\lambda_0\in \mathbb{R}$,
$\mathcal{A}u=-u''$, with $u=0$ at $x=0,1$ and suitable  $u_0$ such that the
solution $u$ to problem \eqref{eqq11} is given by
\[ u(x,t) = \bar{u}(x)  = \begin{cases}
\lambda_0 (1-x_0) x, & x<x_0, \\
\lambda_0 x_0 (1-x), & x>x_0,
\end{cases} \quad (x,t) \in \Omega\times(0,T).
\]
The observation is given on the subdomain $\omega=(0,x_\omega)$ with $x_\omega < x_0$ over the time interval $(0,T)$, i.e.,
$\bar{u}(x) = \lambda_0 (1-x_0) x$ for $(x,t) \in \omega\times(0,T)$.
Note that the source at $x_1\in(x_\omega,1)$ with an intensity $\lambda(t)\equiv\lambda_1=\frac{\lambda_0 (1-x_0)}{1-x_1}$ and suitable $u_0$ leads also to $\bar{u}(x) = \lambda_1 (1-x_1) x = \lambda_0 (1-x_0) x$ in $\omega\times(0,T)$. Thus the inverse problem does not have uniqueness. Similarly, the observations at three internal points cannot uniquely determine two point sources \cite[Section 4]{ABE}.
\end{example}

Now we prove Theorem \ref{t1D}. Similar to Lemma \ref{l2}, we can extend the unique continuation result of \cite{LN}  to the given class of solutions. The proof is similar to Lemma \ref{l2} and hence omitted.
\begin{lemma}\label{l5} Fix $I\subset(0,\ell)$ an open interval and $w_0\in H^1(I)$, let $w\in  H^\alpha(0,T;L^2(I))\cap L^2(0,T;H^2(I))$ be such that $w-w_0\in H_\alpha(0,T;L^2(I))$ and
\begin{equation}\label{l5a}
\rho(x) D_t^\alpha [w-w_0](x,t)+\mathcal A(t) w(x,t)=0,\quad (x,t)\in I\times(0,T).
\end{equation}
Suppose that there exists a  nonempty open subset $U$ of $I$ such that
$w(x,t)=0$ for $(x,t)\in U\times(0,T).$
Then we have $w=0$ in $I\times(0,T)$.
\end{lemma}

Using Lemma \ref{l5}, we can now prove Theorem \ref{t1D}.
\begin{proof}
Suppose that condition \eqref{t1da} holds for some $\epsilon\in(0,T)$. 
Let $u_0=u_0^1-u_0^2$ and $u=u_1-u_2$. Note that $u$ is the solution in the transposition sense of
$$\left\{\begin{aligned}
\rho(x)\partial_t^{\alpha}u +\mathcal{A}(t) u &=  \lambda_0^1(t)\delta_{x_0^1}-\lambda_0^2(t)\delta_{x_0^2}, && \mbox{in }(0,\ell)\times(0,T),\\
 u(0,t)=u(\ell,t)&= 0, && t\in(0,T), \\
u&=u_0, && \mbox{in } (0,\ell)\times \{0\}.
\end{aligned}\right.
$$
In addition, $u\in L^2(0,T; H^1_0(0,\ell))$ and $u-u_0\in H_\alpha(0,T;H^{-1}(0,\ell))$. By Theorem \ref{tt2},
and letting $x_1=\min(x_0^1,x_0^2)$ and $x_2=\max(x_0^1,x_0^2)$, we deduce that, for all $r\in(0,\min(x_1,\ell-x_2))$,
\begin{align*}
u|_{I_{x_1.r}^-\times(0,T)}-u_0|_{I_{x_1,r}^-}\in H_\alpha(0,T;L^2(I_{x_1,r}^-)), \quad u|_{I_{x_1,r}^-\times(0,T)}\in L^2(0,T;H^2(I_{x_1,r}^-)),\\
u|_{I_{x_2,r}^+\times(0,T)}-u_0|_{I_{x_2,r}^+}\in H_\alpha(0,T;L^2(I_{x_2,r}^+)),\quad
 u|_{(I_{x_2,r}^+\times(0,T)}\in L^2(0,T;H^2(I_{x_2,r}^+)).
\end{align*}
Moreover, we have
\begin{equation}\label{t1dd}
\rho(x)D_t^{\alpha}(u-u_0)(x,t) +\mathcal{A}(t) u(x,t)=0,\quad (x,t)\in(I_{x_1,r}^-\cup I_{x_2,r}^+)\times (0,T).
\end{equation}
Meanwhile, we have $I_1\subset(0,x_1)$, $I_2\subset(x_2,\ell)$ and, for $r\in(0,\min(x_1,\ell-x_2))$ sufficiently small, we have  $I_{1,r}=I_1\cap I_{x_1,r}^- \neq\emptyset$, $I_{2,r}=I_2\cap I_{x_2,r}^+\neq\emptyset$. This and condition \eqref{t1da} imply
$$u(x,t)=D_t^{\alpha}(u-u_0)(x,t)=0,\quad (x,t)\in(I_{1,r}\cup I_{2,r})\times (T-\epsilon,T).$$
Thus, repeating the argument in Step 1 of Theorem \ref{t1} gives
$u=0$ in $(I_{1,r}\cup I_{2,r})\times (0,T),$
and Lemma \ref{l5} implies
$u=0$ in $( I_{x_1,r}^-\cup I_{x_2,r}^+)\times (0,T).$
Since $r\in(0,\min(x_1,\ell-x_2))$ can be arbitrarily small, we obtain
\begin{equation}\label{t1df}
u(x,t)=0,\quad (x,t)\in((0,x_1)\cup(x_2,\ell))\times (0,T).
\end{equation}
This and \eqref{t1dd}, with $r\in(0,\min(x_1,\ell-x_2))$, yield
$$-\rho(x)\frac{t^{-\alpha}}{\Gamma(1-\alpha)}u_0(x)=\rho(x) D_t^\alpha [u-u_0](x,t)+\mathcal A u(x,t)=0,\quad (x,t)\in (I_{x_1,r}^-\cup I_{x_2,r}^+)\times(0,T),$$
which implies $u_0=0$ on $(0,x_1)\cup(x_2,\ell)$. Then, from condition \eqref{ini}, we deduce $u_0\equiv0$.
Next we prove $x_0^1=x_0^2$ by contradiction. By assuming the contrary, we may assume $x_0^1<x_0^2$. Then, \eqref{t1df} holds with $x_1=x_0^1$, $x_2=x_0^2$. Since $u_0\equiv0$, condition \eqref{t1df} implies  that $v=u|_{(x_1,x_2)\times(0,T)}\in H_\alpha(0,T;H^{-1}(x_1,x_2))\cap L^2(0,T;H^1_0(x_1,x_2))$ satisfies for any $\phi\in H^1_0(x_1,x_2),\ t\in(0,T)$,
\begin{equation}\label{t1ddf}
\begin{aligned}&\left\langle \rho D_t^\alpha v(\cdot,t)+\mathcal A(t) v(\cdot,t),\phi\right\rangle_{H^{-1}(x_1,x_2),H^1_0(x_1,x_2)}\\
=&\left\langle \rho D_t^\alpha u(\cdot,t)+\mathcal A(t) u(\cdot,t),\phi_\star\right\rangle_{H^{-1}(0,\ell),H^1_0(0,\ell)}\\
=&\lambda_0^1(t)\phi(x_0^1)-\lambda_0^2(t)\phi(x_0^2)=0,
\end{aligned}
\end{equation}
with $\phi_\star$ being the zero extension of $\phi$ to $(0,\ell)$.
\cite[Theorem 4.1]{KRY} and  \eqref{t1ddf} imply $v\equiv0$. This and \eqref{t1df} give $u\equiv0$ and,  similar to Step 4 of the proof of Theorem \ref{t1}, it contradicts the condition $x_0^1\neq x_0^2$. Thus, $x_1=\min(x_0^1,x_0^2)=\max(x_0^1,x_0^2)=x_2$ and, then the identity \eqref{t1df} implies \eqref{t1db}.
\end{proof}

\section{Numerical results and discussions}\label{sec:numer}

Now we present numerical reconstructions using the Levenberg-Marquardt algorithm \cite{Levenberg:1944,Marquardt:1963}. For example, to recover one point source $x^*$ and its strength $\lambda^*$, we define a nonlinear operator $F:(x,\lambda)\in \Omega \times L^2(0,T) \rightarrow u_{x,\lambda}\in L^2(T-\epsilon,T;L^2(\omega))$, with $u$ solving problem \eqref{eq1} with parameters $x$ and $\lambda$. Fix an initial guess $(x^0,\lambda^0)$. Now given $(x^k,\lambda^k)$, the next approximation $(x^{k+1},\lambda^{k+1})$ is given by
\begin{equation*}
    (x^{k+1},\lambda^{k+1})={\arg\min}_{x,\lambda} J_k(x,\lambda),
\end{equation*}
with the functional $J_k(x,\lambda)$ based at $(x^k,\lambda^k)$ given by
\begin{align*}
    J_k(x,\lambda)=&\|F(x^{k},\lambda^k)-z^{\delta}+\partial_xF(x^k,\lambda^k)(x-x^k)+\partial_{\lambda} F(x^k,\lambda^k)(\lambda-\lambda^k) \|^2_{L^2(T-\epsilon,T;L^2(\omega))}\\
    &+\beta_x^k|x-x^k|^2+{\beta_{\lambda}^k}\|\lambda-\lambda^k\|_{H^1(0,T)}^2,
\end{align*}
where $z^\delta$ is the noisy data, $(T-\epsilon,T)$ is the time horizon, $\beta_x^k,\beta_\lambda^k,\mu>0$ are scalars, and $\partial_x F (x^k,\lambda^k)$ and $\partial_\lambda F (x^k,\lambda^k)$ are the Jacobians of the map $F$ in $x$ and $\lambda$, respectively. The parameters $\beta_x^k$ and $\beta_\lambda^k$ often decrease geometrically with factors $\gamma_x,\gamma_\lambda\in(0,1)$: $\beta_x^{k+1}=\gamma_x\beta_x^k $ and $\beta_\lambda^{k+1}=\gamma_\lambda \beta_\lambda^k $.
To prevent over-fitting, we use early stopping and obtain $x^K$ and $\lambda^K$ at the $K$-th iteration.

The direct problem is solved using the Galerkin FEM with continuous piecewise linear elements in space, and backward Euler convolution quadrature in time \cite{JinZhou:2023book}. Unless otherwise stated, the order $\alpha$ is fixed at $0.5$. The noisy data $z^\delta$ is generated by adding i.i.d. Gaussian noise of mean zero and standard deviation $\delta \norm{u_{x^*,\lambda^*}}_{L^\infty}$ to the exact data $u_{x^*,\lambda^*}$, with $\delta$ denoting the noise level.

\subsection{1D numerical examples}

First we present numerical results for 1D examples with $\Omega=(0,1)$, $T=1$, $\rho\equiv1$, and $\mathcal{A}=-\Delta$. The direct problem is solved on a fine grid with a mesh size $h=2\times10^{-3}$ and time step size $\tau =10^{-3}$, and the inverse problem on a coarse grid with $h=10^{-2}$ and $\tau=5\times10^{-3}$. 

\begin{example}\label{exam1}
Consider one point source at $x=0.5$ with an intensity $\lambda(t)=0.2e^t$ for $t\in[0,T]$, and $u_0(x)=x(1-x)$. The data is given over the subdomain $\omega=(0,\frac14)\cup(\frac34,1)$ for the time interval $(T-\epsilon,T)$ with $\epsilon=\frac34T$ or $\epsilon=\frac12T$.
\end{example}

\begin{figure}[htb!]
\centering\setlength{\tabcolsep}{0pt}
\begin{tabular}{cccc}
	\includegraphics[width=.24\textwidth]{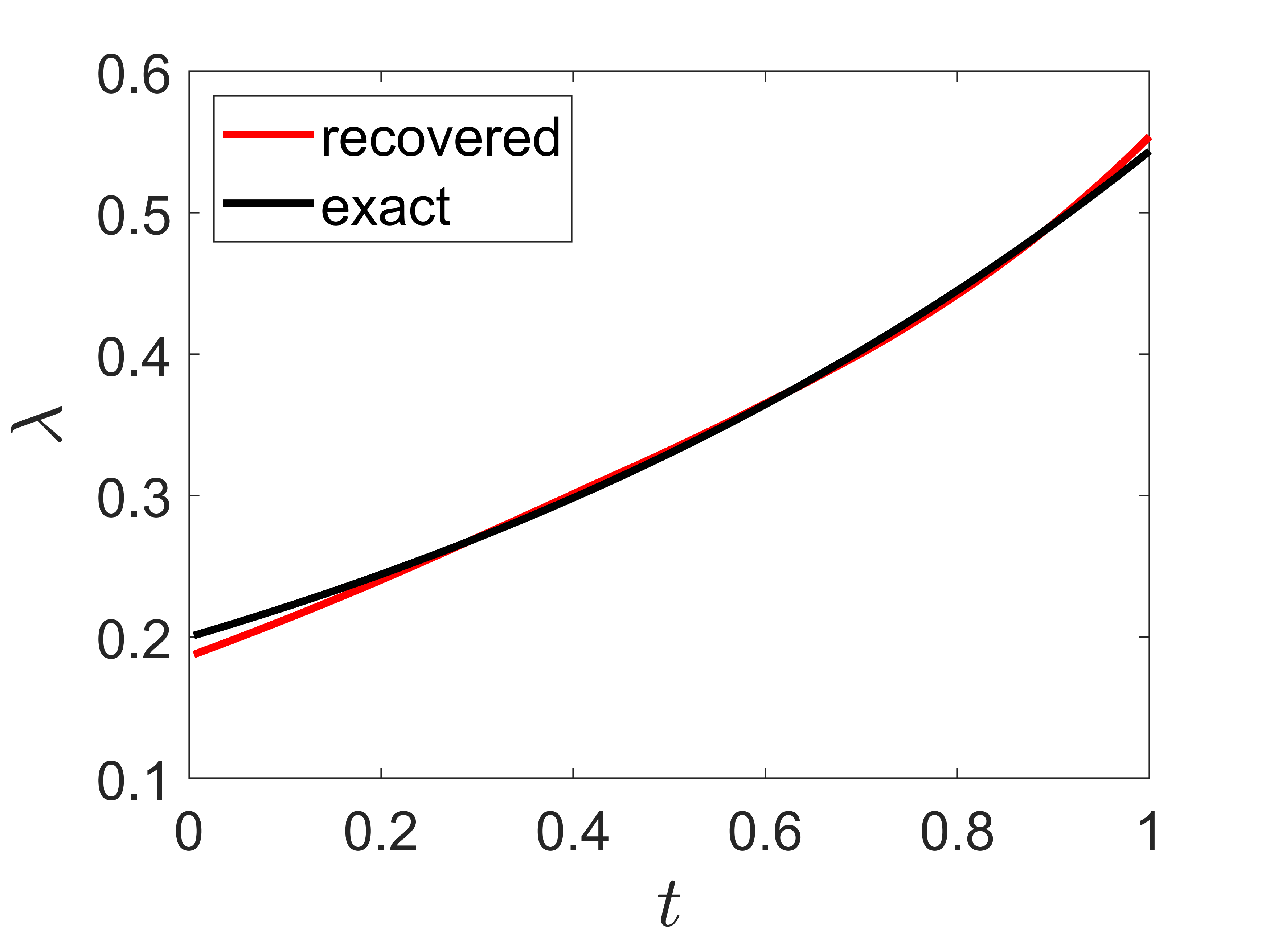} & \includegraphics[width=.24\textwidth]{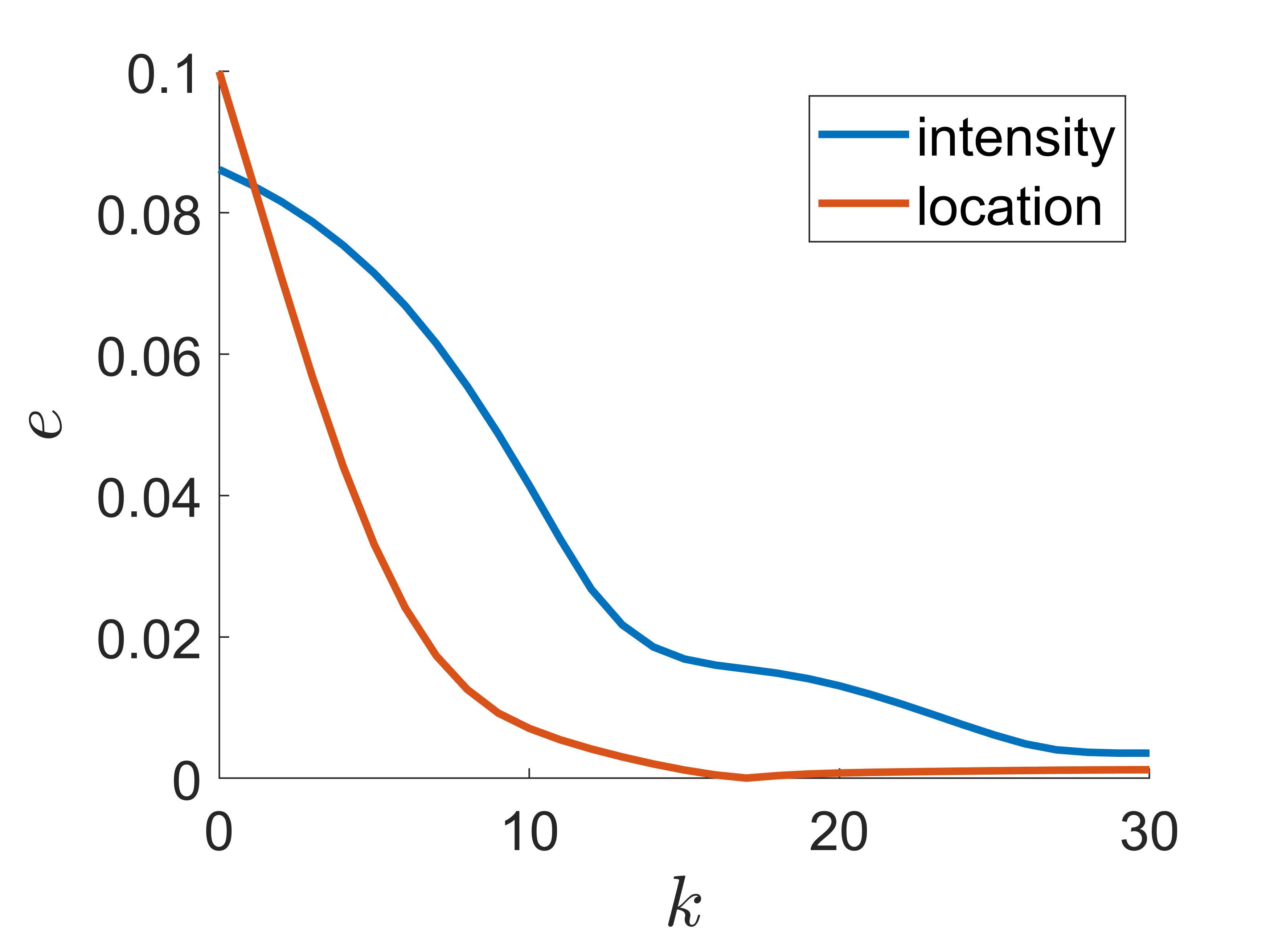} & \includegraphics[width=.24\textwidth]{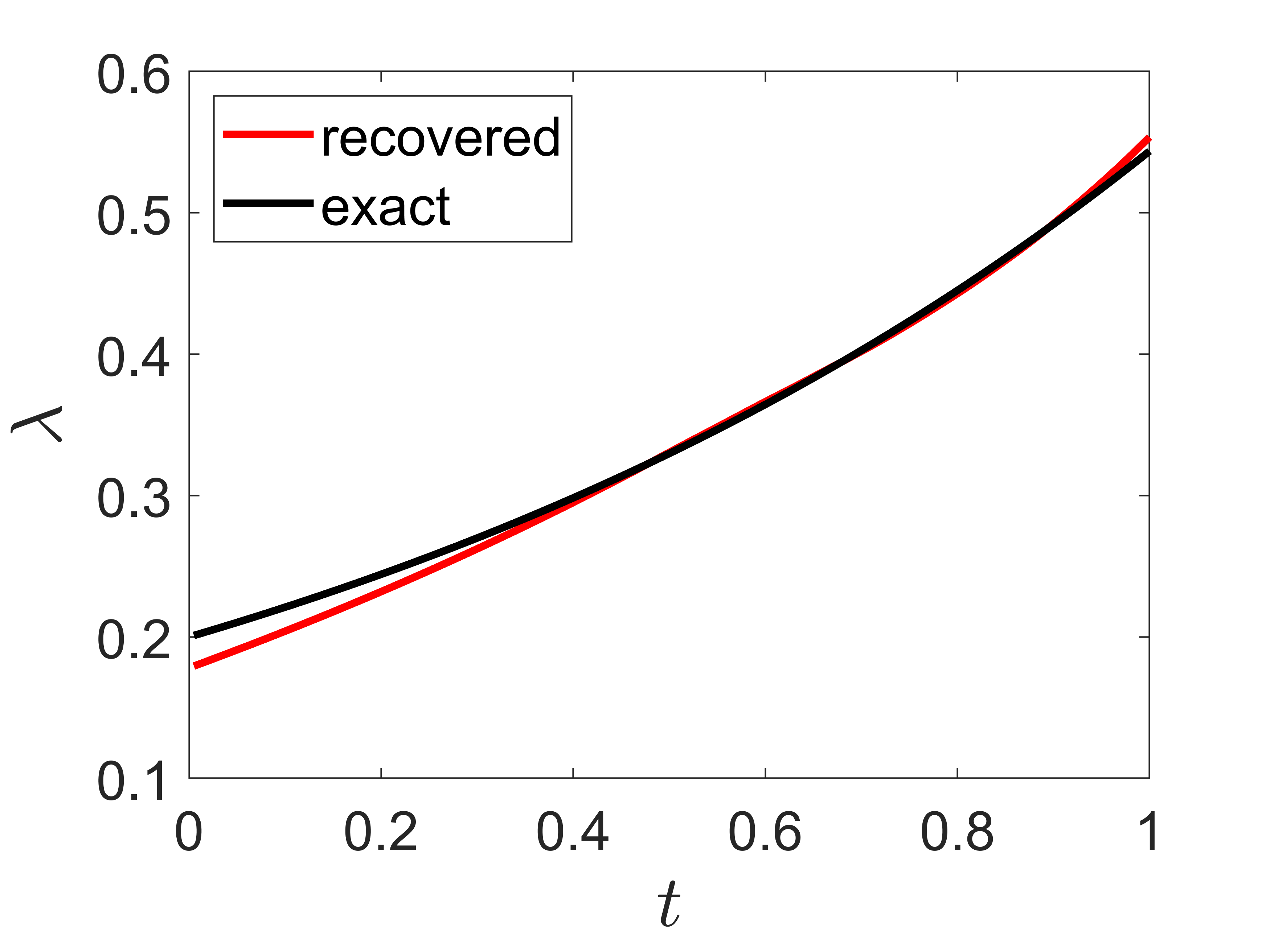} & \includegraphics[width=.24\textwidth]{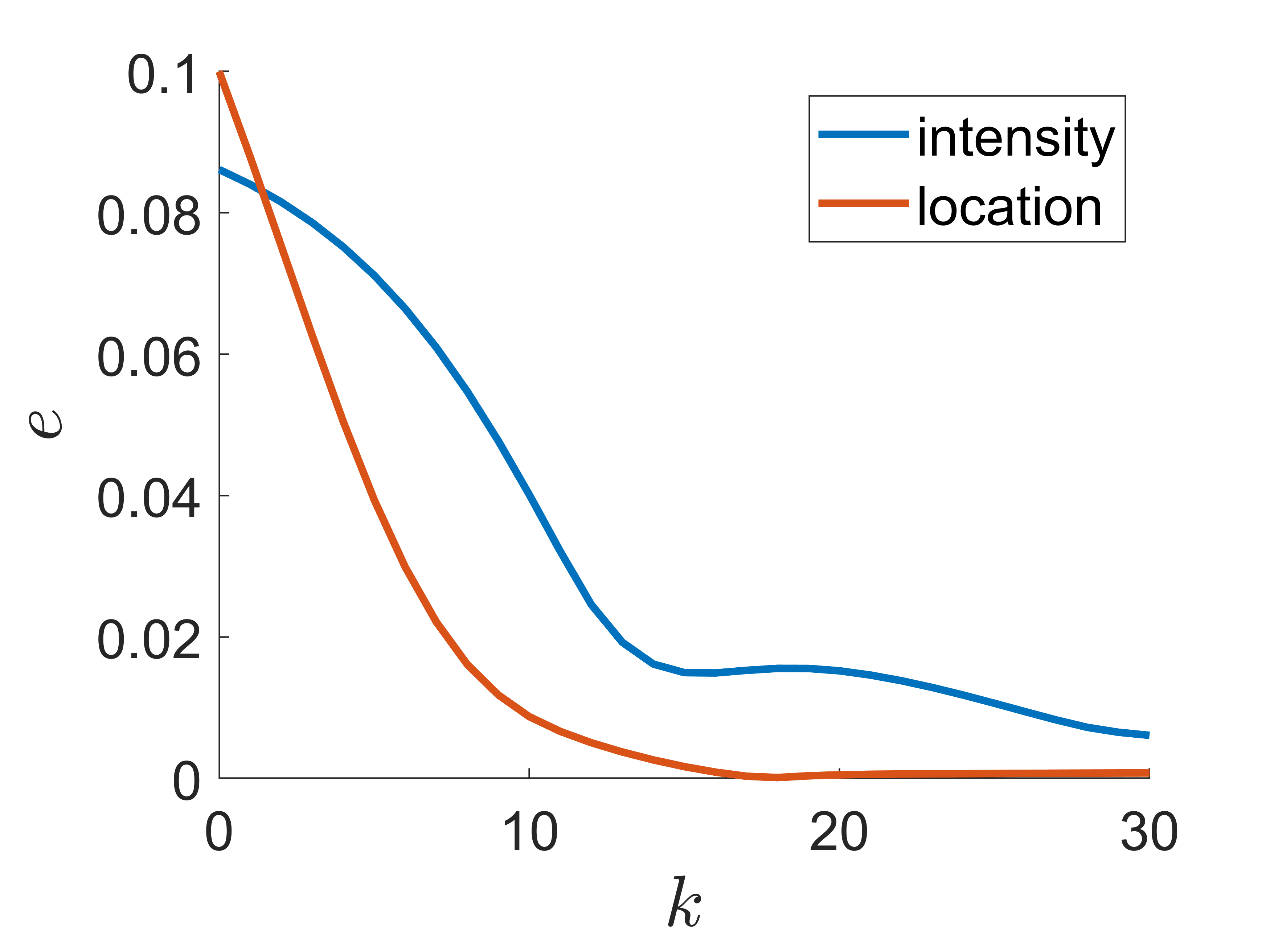} \\
 \includegraphics[width=.24\textwidth]{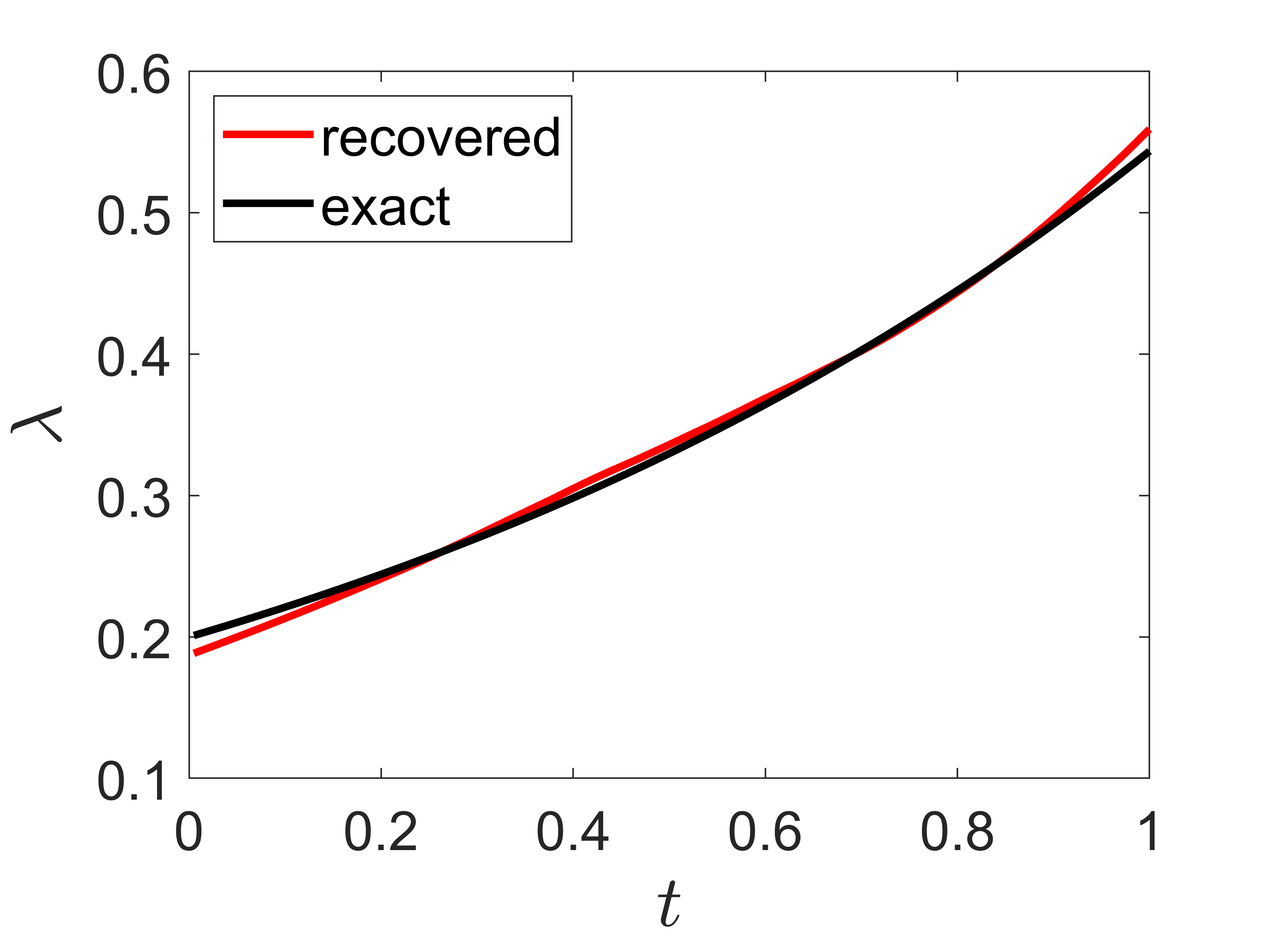} & \includegraphics[width=.24\textwidth]{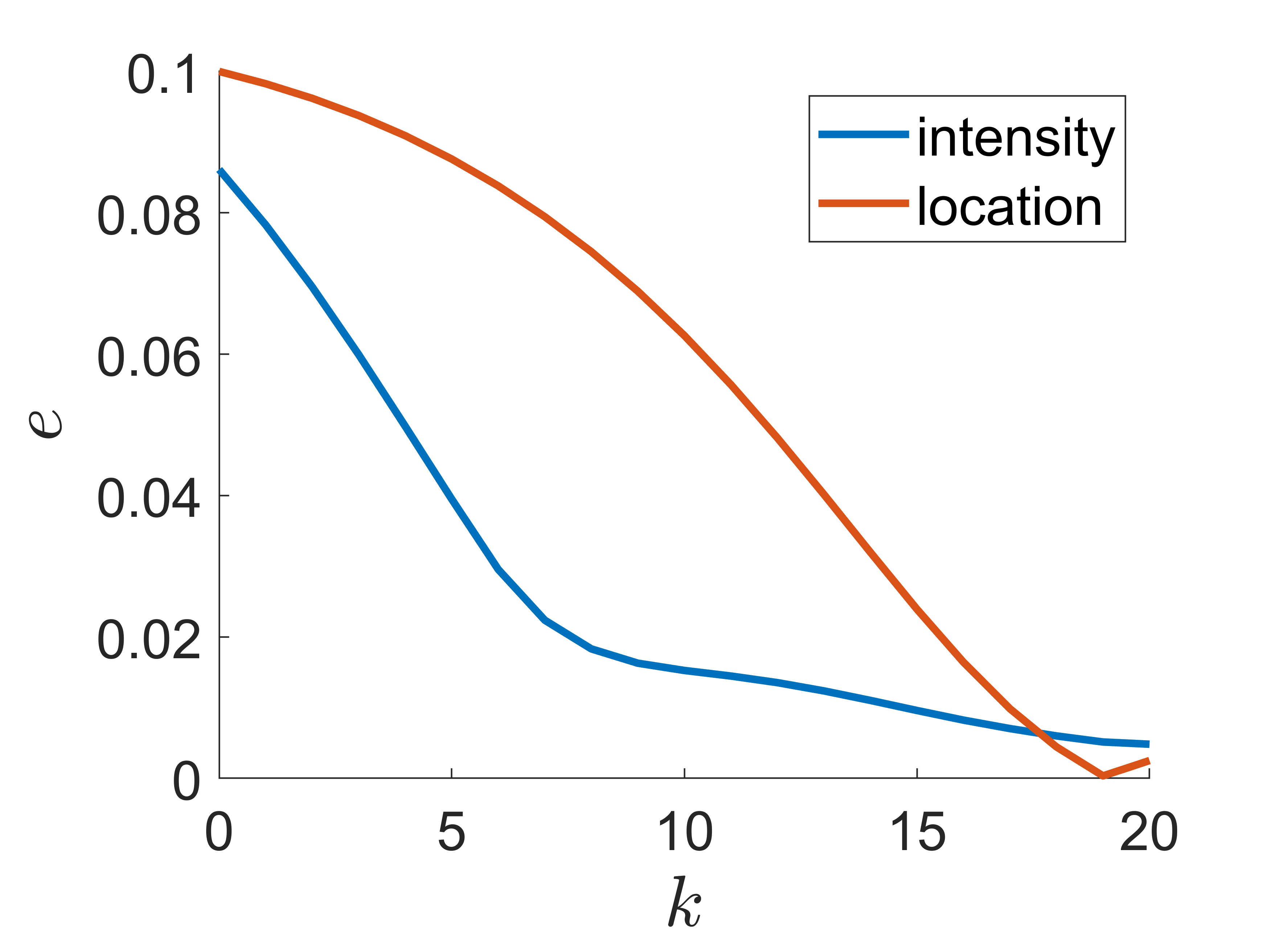} & \includegraphics[width=.24\textwidth]{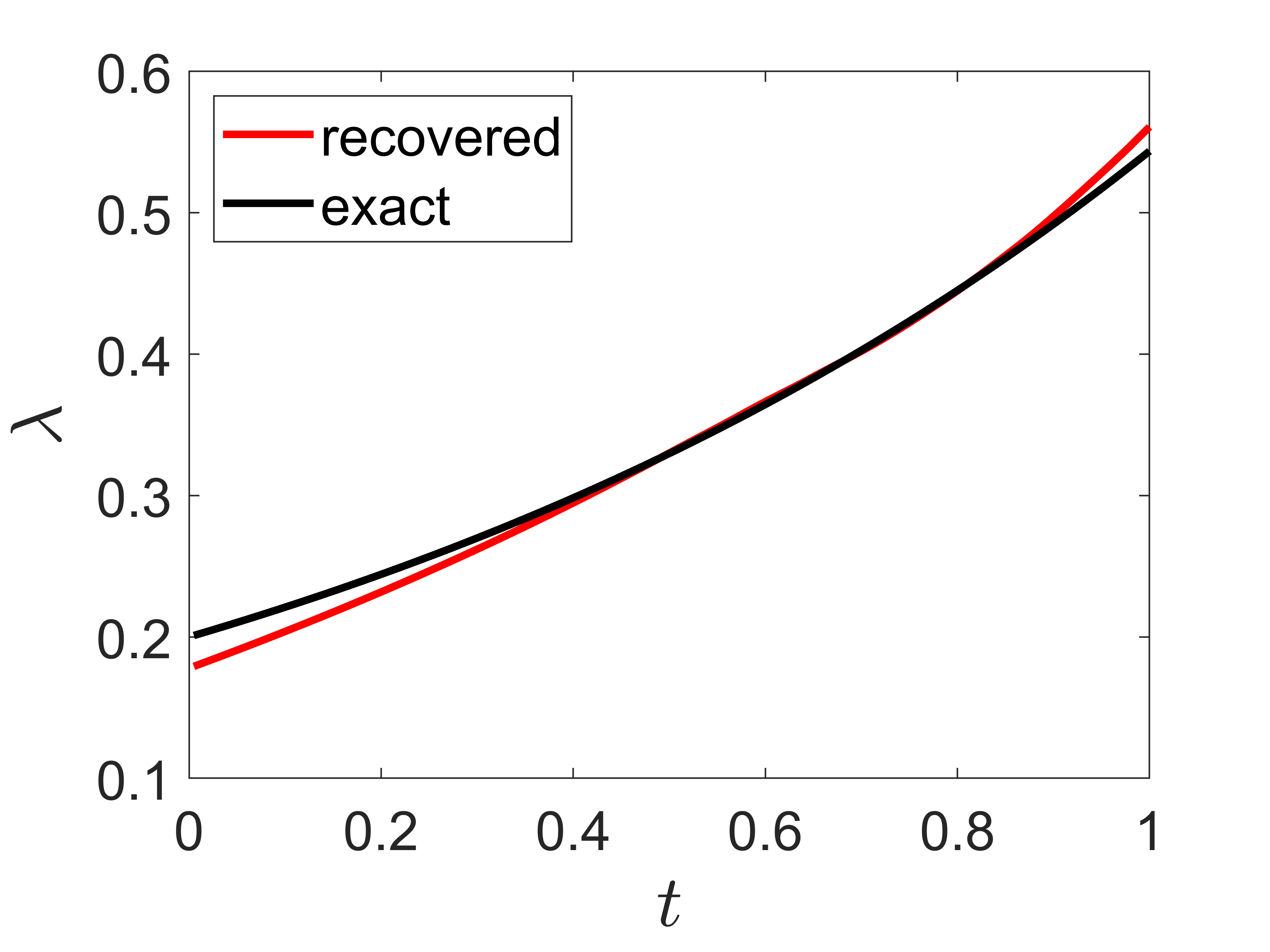} & \includegraphics[width=.24\textwidth]{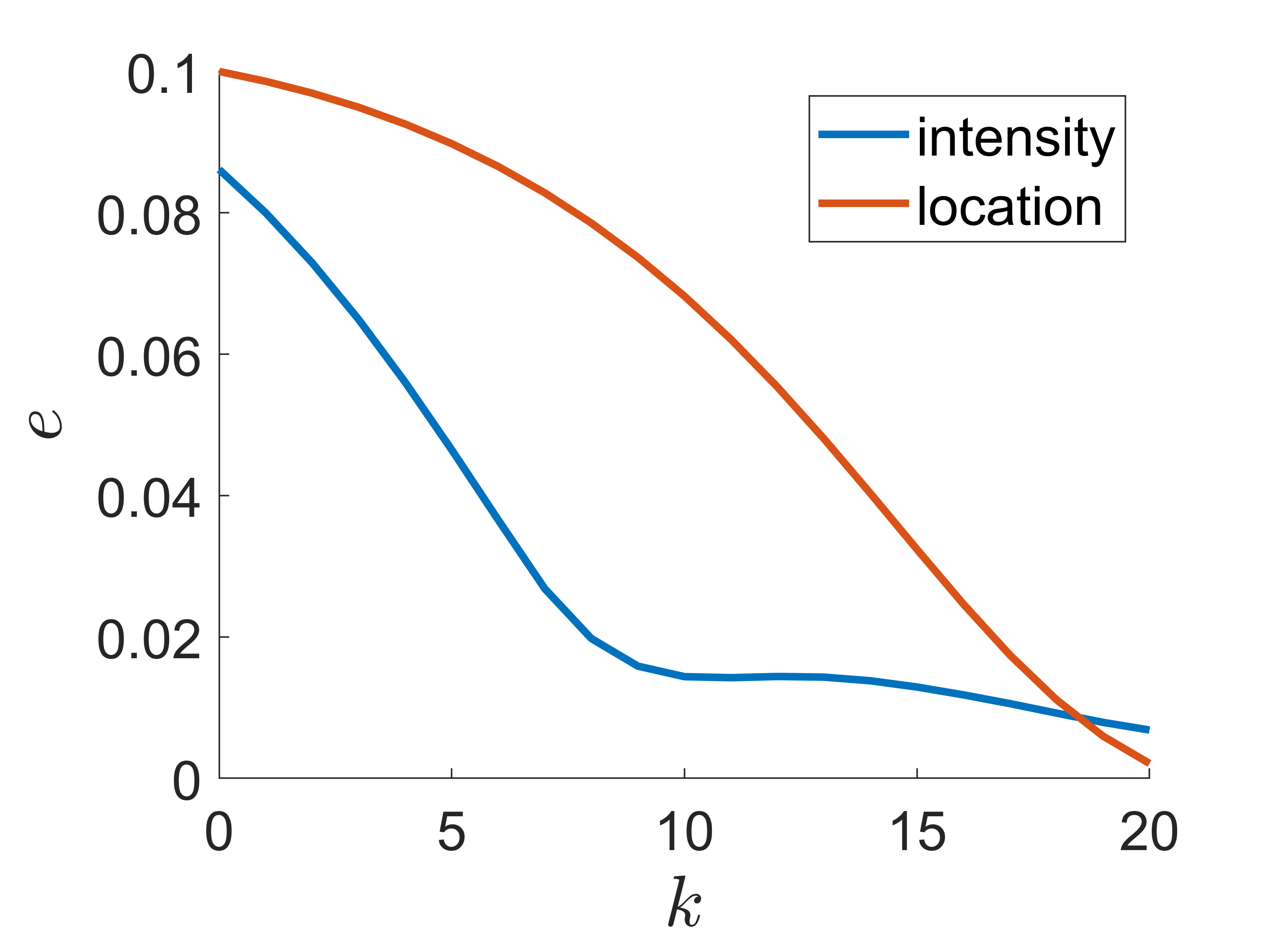}\\
 (a) $\lambda$ & (b) convergence & (c) $\lambda$ & (d) convergence
\end{tabular}
\caption{The recovered $\lambda$ and convergence for Example~\ref{exam1} with $\delta=2\%$ (top) and $\delta=10\%$ (bottom). (a)--(b) and (c)--(d) are for $\epsilon=\frac34T$ and $\epsilon=\frac12T$, respectively.}
\label{fig:exp_1}
\end{figure}

In Table~\ref{tab:exp_1}, we present the recovered location for different $\epsilon$ and $\delta$, with the initial guess  $x^0=0.4$ and $\lambda^0(t)=0.25e^t$, and in Fig.~\ref{fig:exp_1} the recovered intensity. The results show that the recovered location and intensity are accurate for $\delta$ up to 10\%, and also the choice of the time domain (i.e., $\epsilon=\frac34T$ or $\epsilon=\frac12T$) has little influence on the convergence of the algorithm. Fig.~\ref{fig:exp_1}(b) and (d) show the steady decay of errors for both location and intensity during the iteration, and the convergence is reached within tens of iterations, showing the effectiveness of the algorithm.

\begin{table}[htb!]
    \centering
    \renewcommand{\arraystretch}{1.5}
    \begin{tabular}{c|cccc}
    \toprule
    $(\epsilon,\delta)$ & $(\frac34T,2\%)$ & $(\frac34T,10\%)$ & $(\frac12T,2\%)$ & $(\frac12T,10\%)$ \\
    \midrule
     location $x^K$ & $0.5012$ & $0.5025$ & $0.5008$ & $0.4979$ \\
    \bottomrule
    \end{tabular}
    \caption{The recovered source location for Example~\ref{exam1}}
    \label{tab:exp_1}
\end{table}

\begin{example}\label{exam2}
Consider two point sources at $x_1=0.3$ and $x_2=0.7$, with intensities $\lambda_1(t)=0.2e^t$ and $\lambda_2(t)=0.2\sin(2\pi t)+0.4$, and $u_0(x)=x(1-x)$.
\begin{itemize}
    \item[{\rm(i)}] Recover point sources from the data in the time interval $(\frac{1}{4}T,T)$ over the subdomain $\omega=(0,0.2)\cup(0.4,0.6)\cup(0.8,1)$.
    \item[{\rm(ii)}] Recover point sources and initial data from the data in the time interval $(\frac12T,T)$ over the subdomain $\omega=(0,0.25)\cup(0.35,0.65)\cup(0.75,1)$.
\end{itemize}
\end{example}

\begin{figure}[htbp!]
\centering\setlength{\tabcolsep}{0pt}
\begin{tabular}{ccc}
    \includegraphics[width=.32\textwidth]{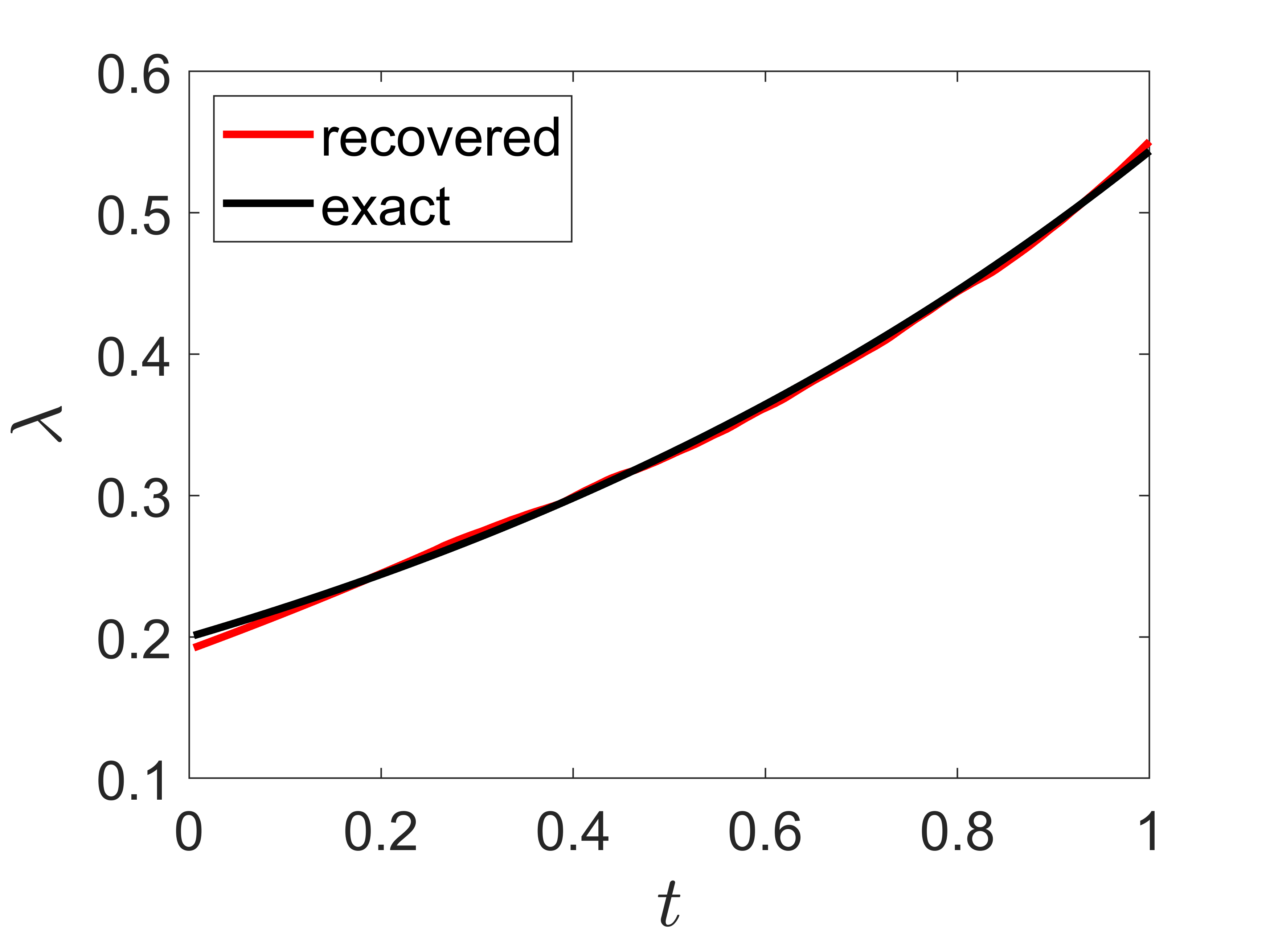} & \includegraphics[width=.32\textwidth]{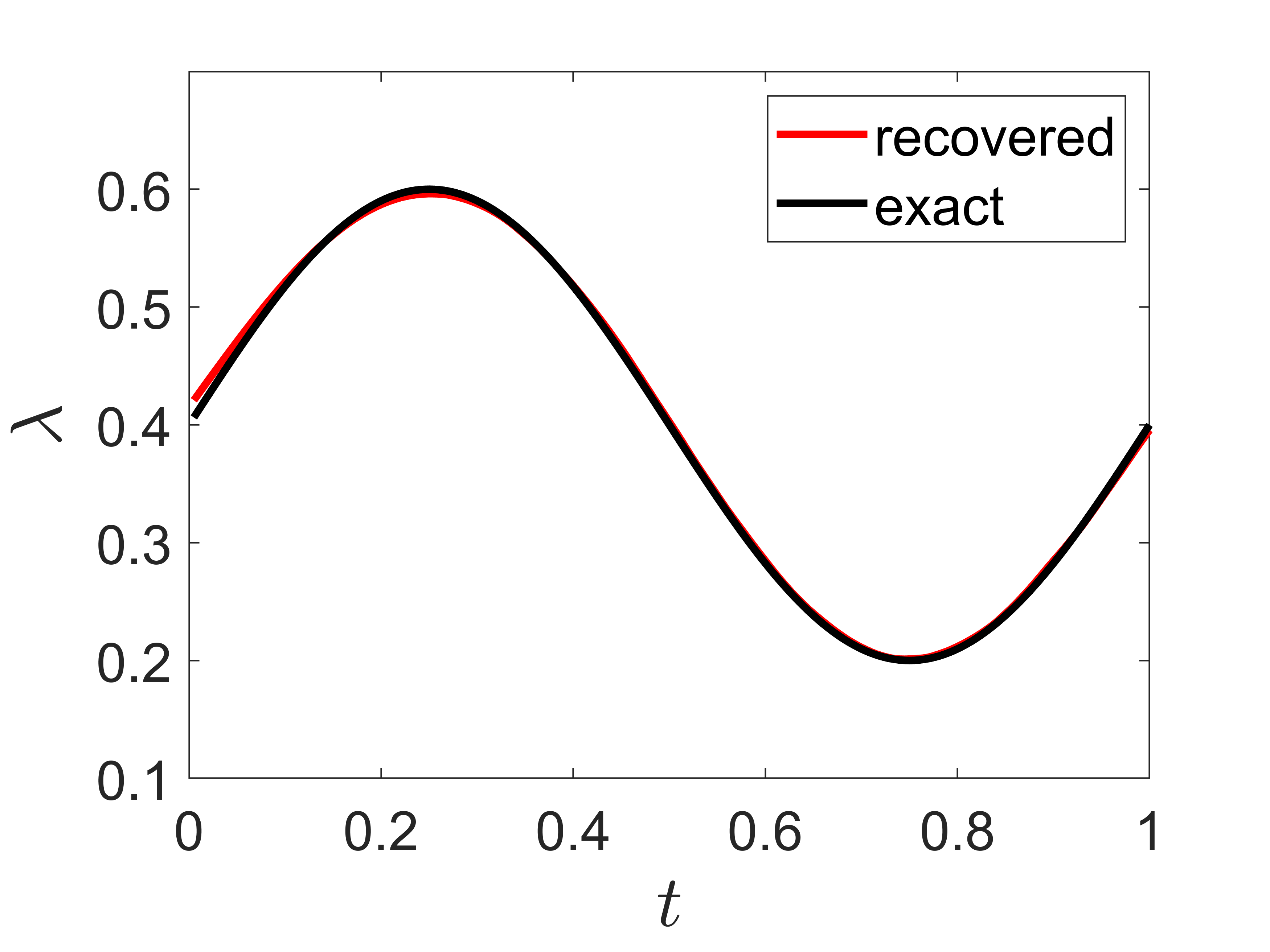} & \includegraphics[width=.32\textwidth]{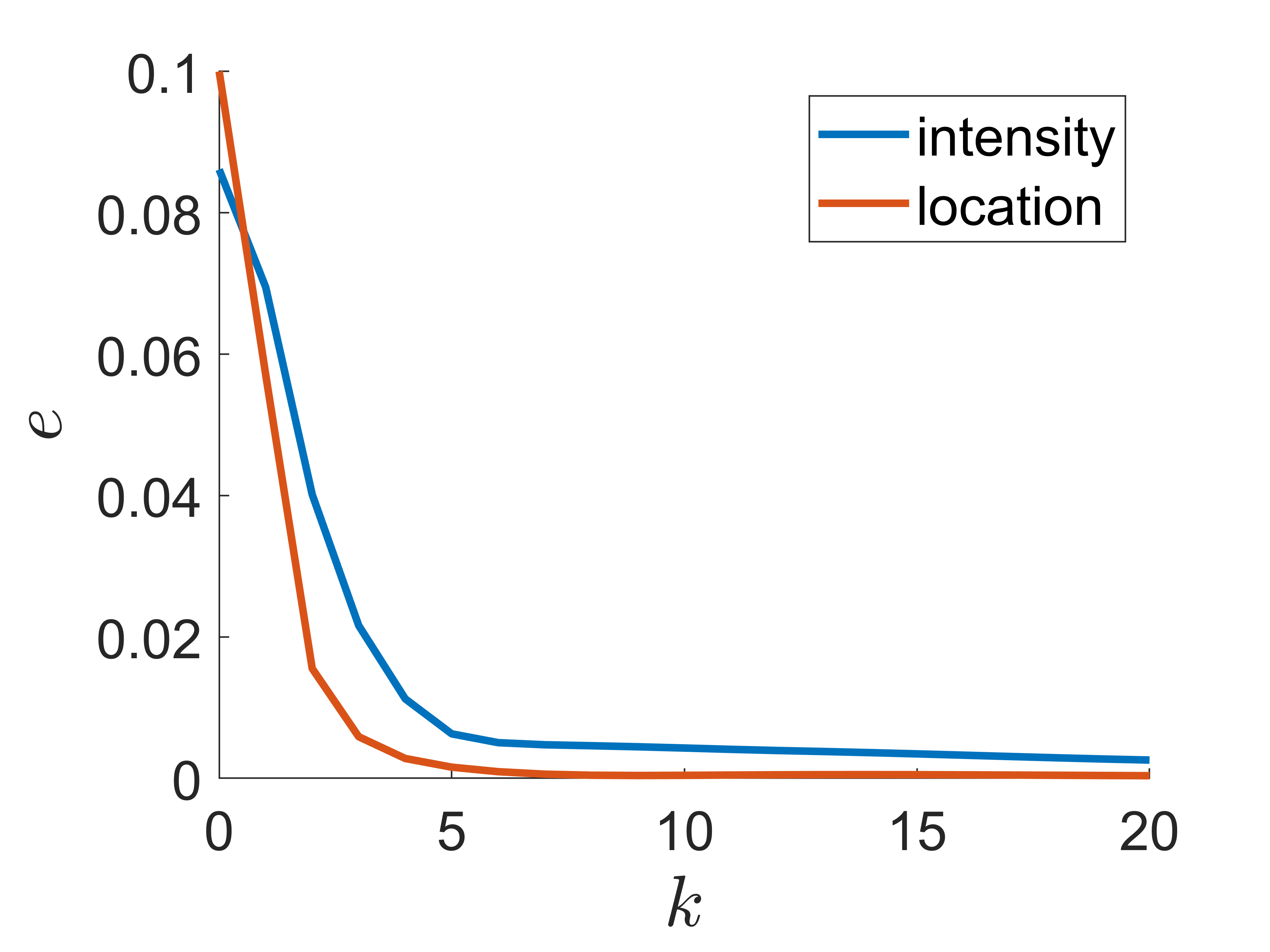} \\
    \includegraphics[width=.32\textwidth]{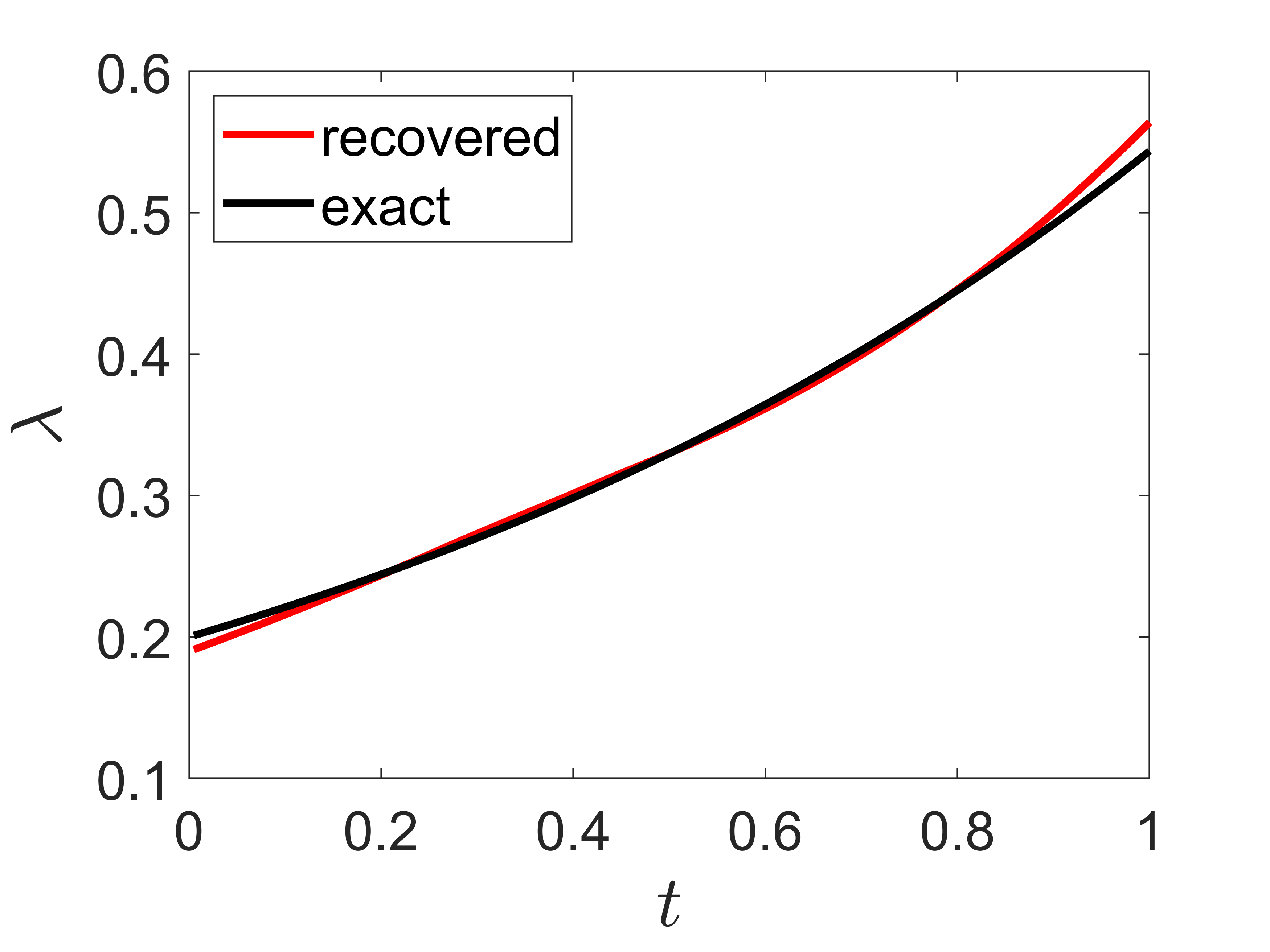} & \includegraphics[width=.32\textwidth]{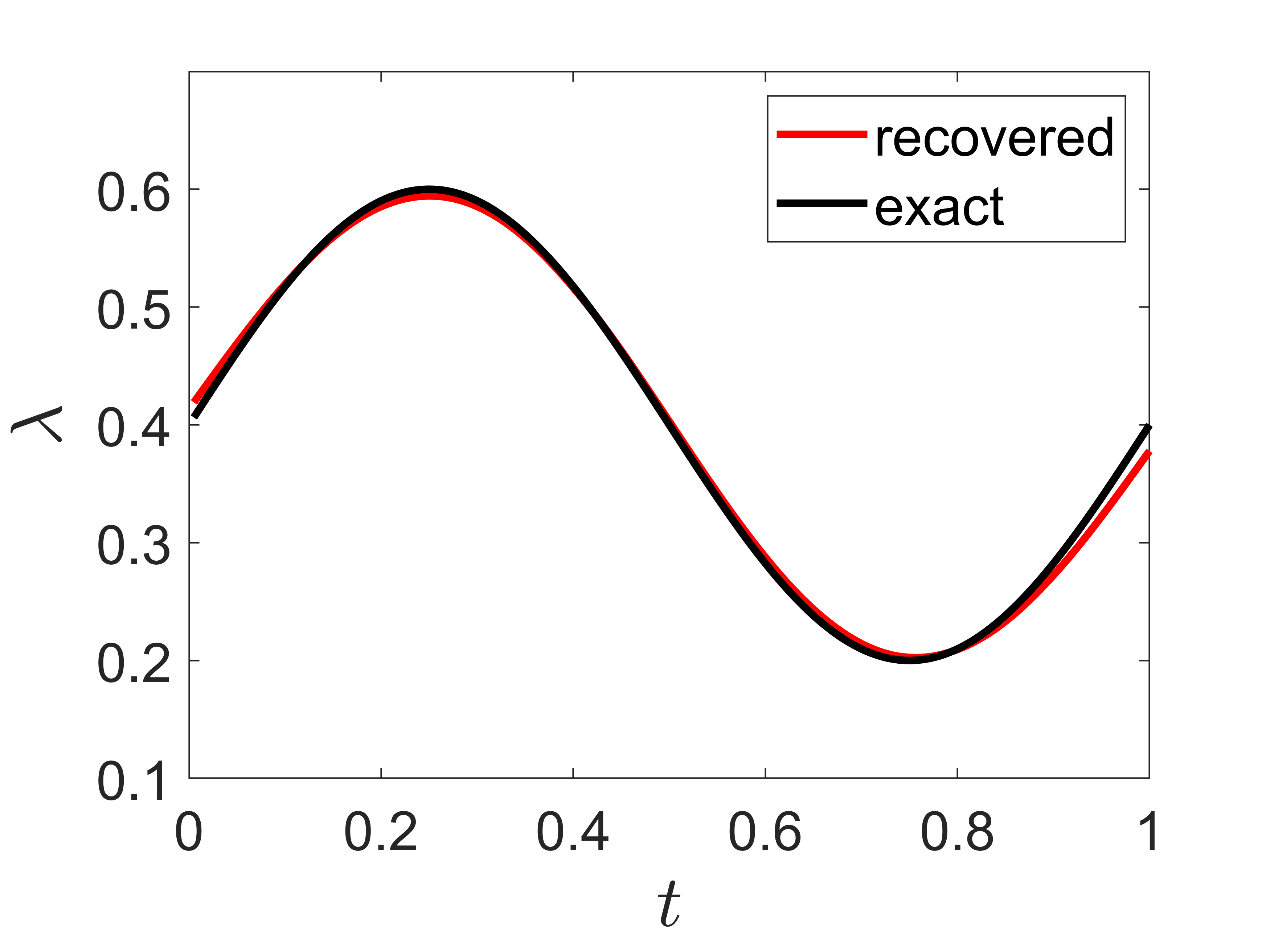} & \includegraphics[width=.32\textwidth]{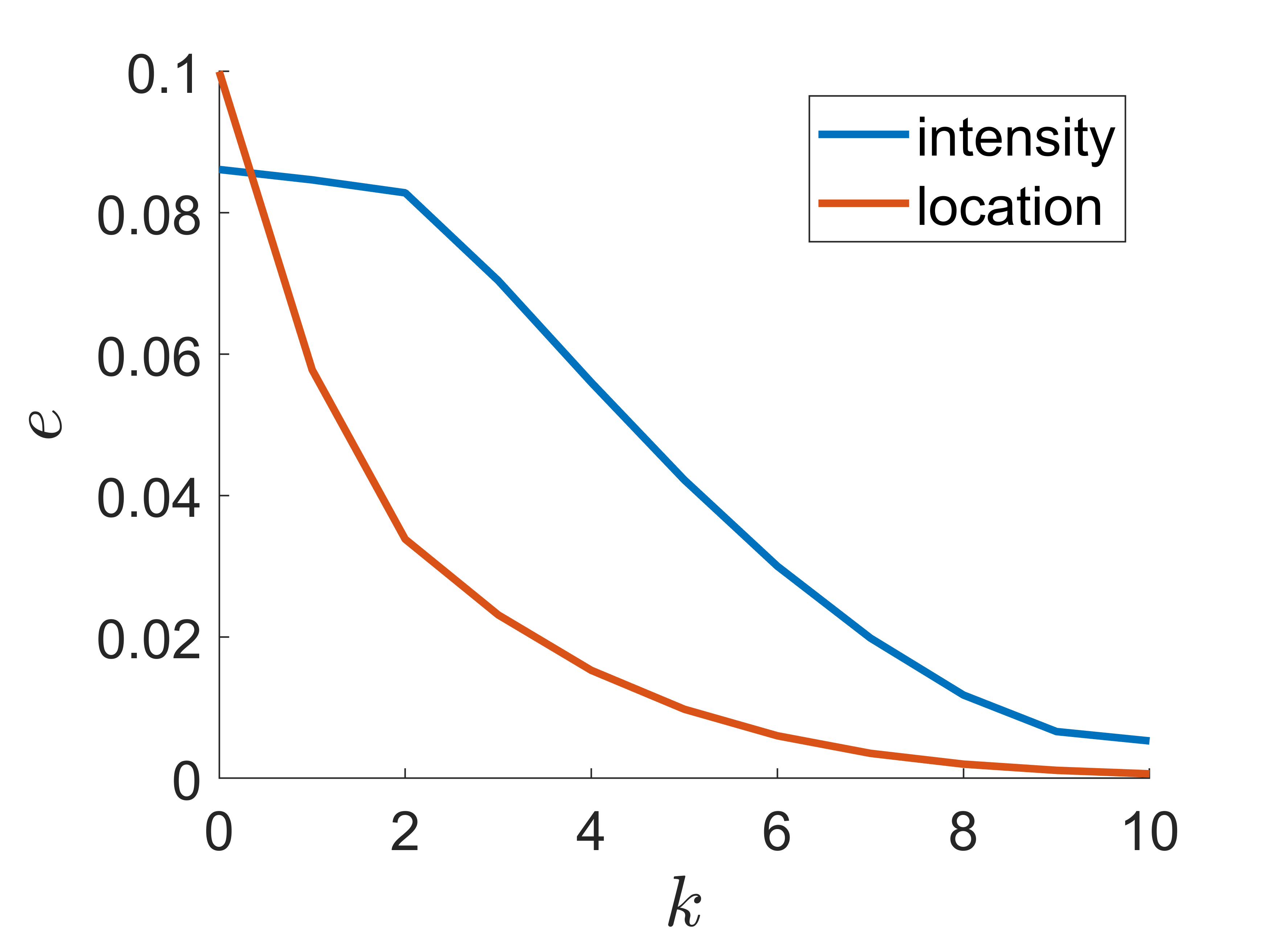} \\
    (a) $\lambda_1$ & (b) $\lambda_2$  & (c) convergence
\end{tabular}
\caption{The recovered $\lambda_i$ and convergence for Example~\ref{exam2}(i) with $\delta=2\%$ (top) and $\delta=10\%$ (bottom).}
\label{fig:exp_2}
\end{figure}

In case (i), we initialize the algorithm with $x_1^0=0.4$, $x_2^0=0.6$ and $\lambda_1^0(t)=0.25e^t, \lambda_2^0(t)=0.18\sin(2\pi t)+0.36$. The recovered source locations are 0.3000 and 0.7004
for $2\%$ noise, and 0.3006 and 0.7005 for 10\% noise, which are fairly accurate. In Fig.~\ref{fig:exp_2}, we show the recovered intensities, which have reasonable accuracy, with the errors concentrating around the end points, and that the convergence of the recovered locations and intensities is very steady.

\begin{figure}[htb!]
\centering\setlength{\tabcolsep}{0pt}
\begin{tabular}{ccc}
\includegraphics[width=.32\textwidth]{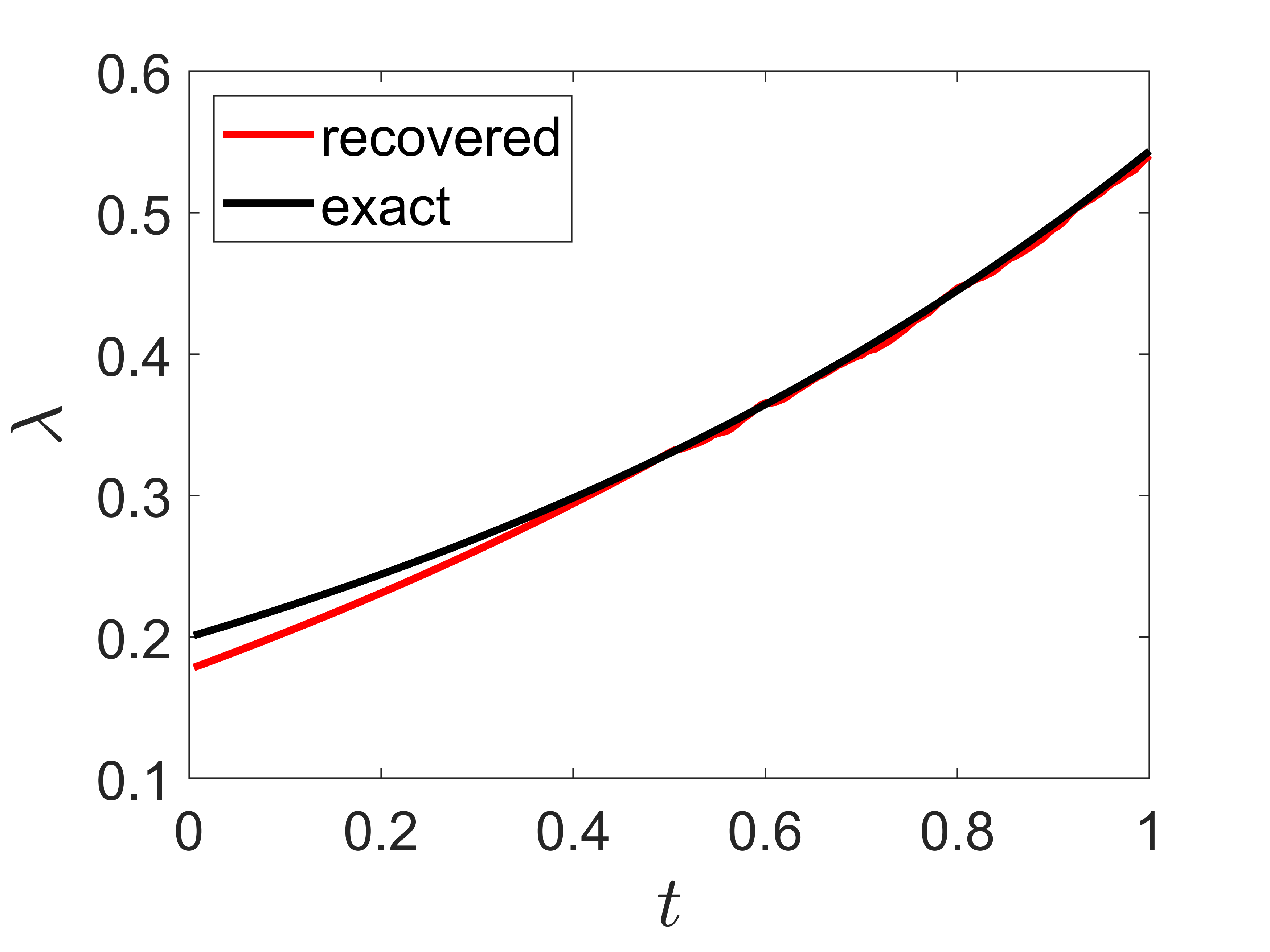} & \includegraphics[width=.32\textwidth]{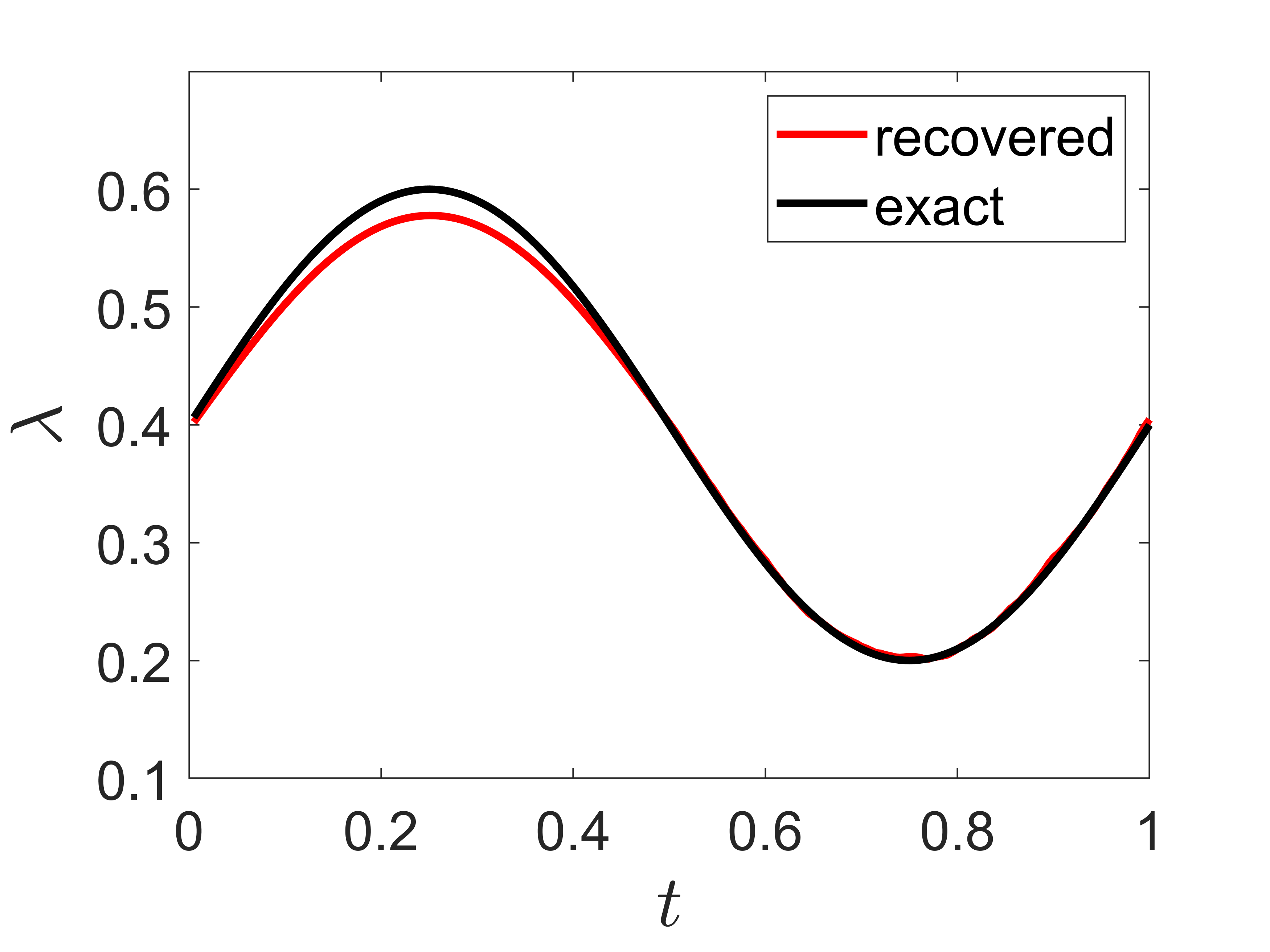} & \includegraphics[width=.32\textwidth]{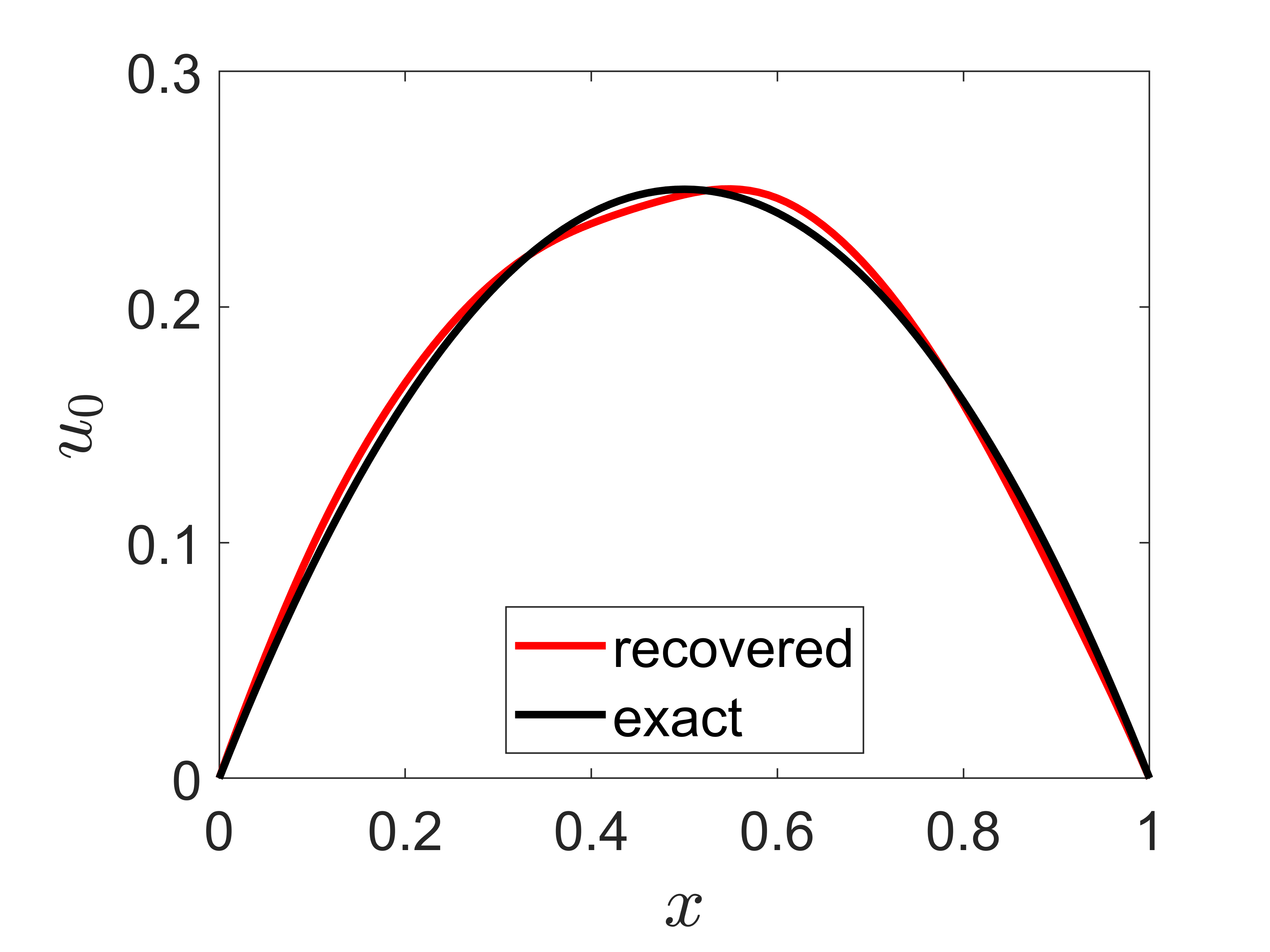} \\
(a) $\lambda_1$ & (b) $\lambda_2$ & (c) $u_0$\\
\includegraphics[width=.32\textwidth]{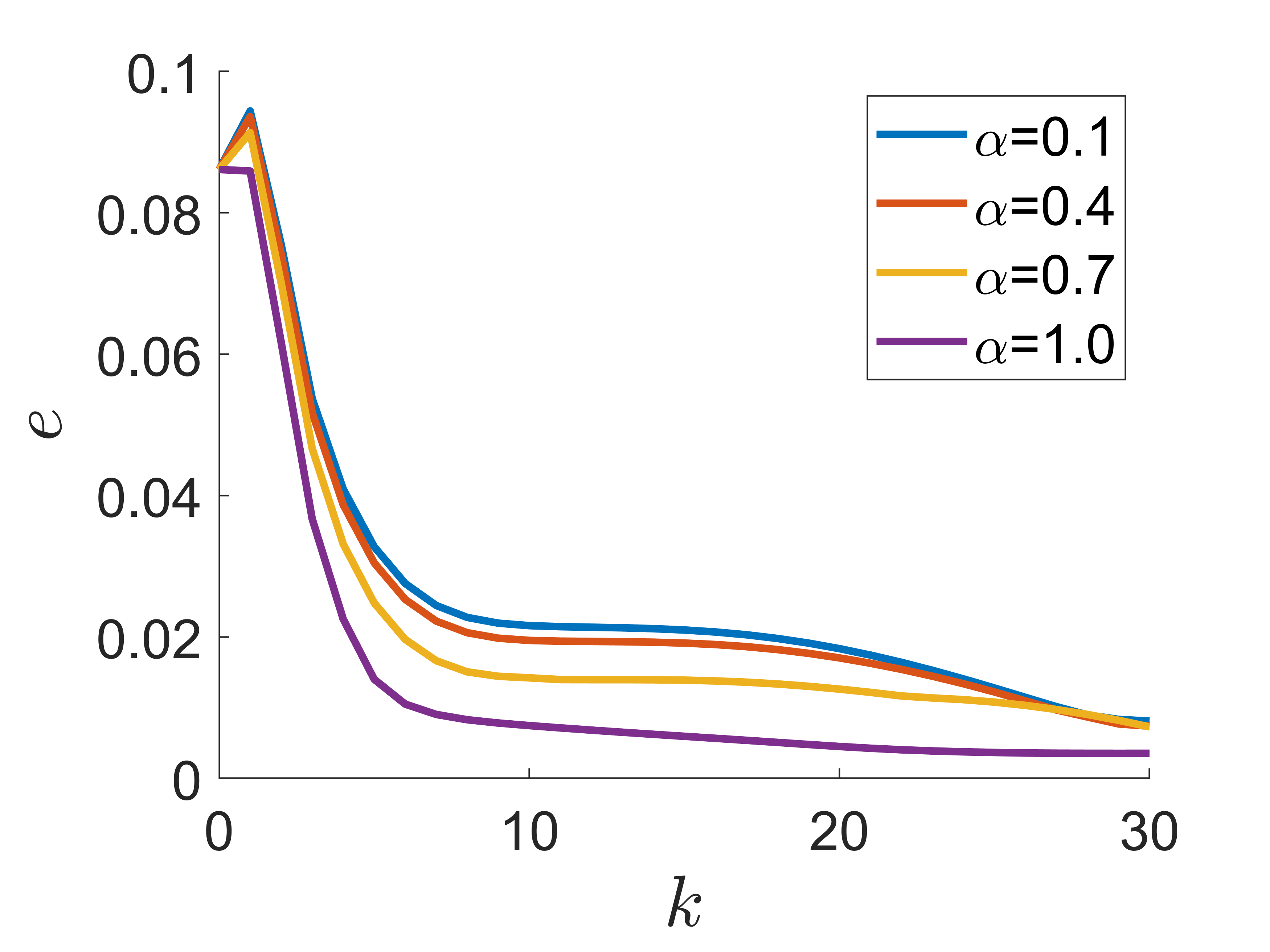} & \includegraphics[width=.32\textwidth]{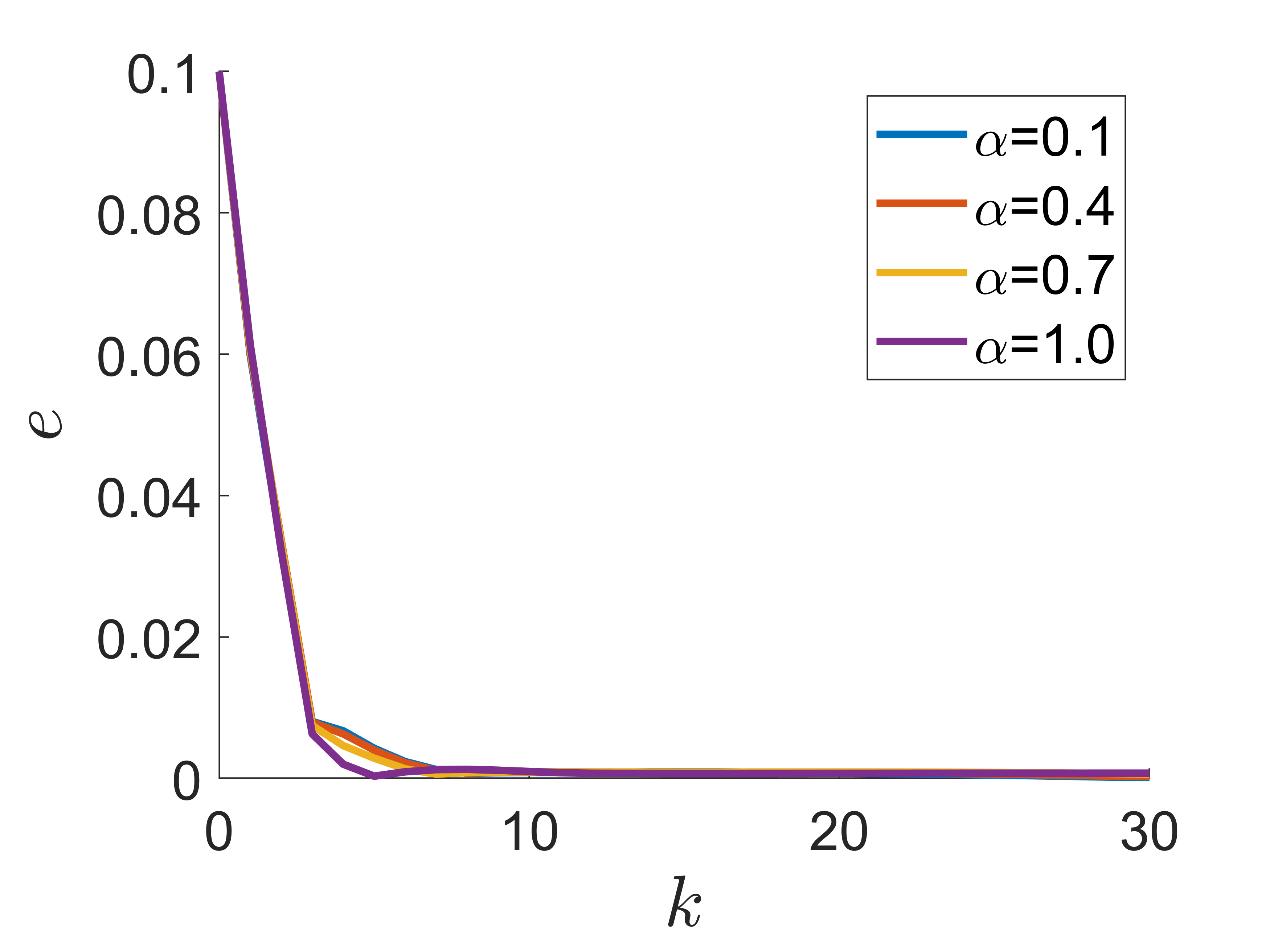} & \includegraphics[width=.32\textwidth]{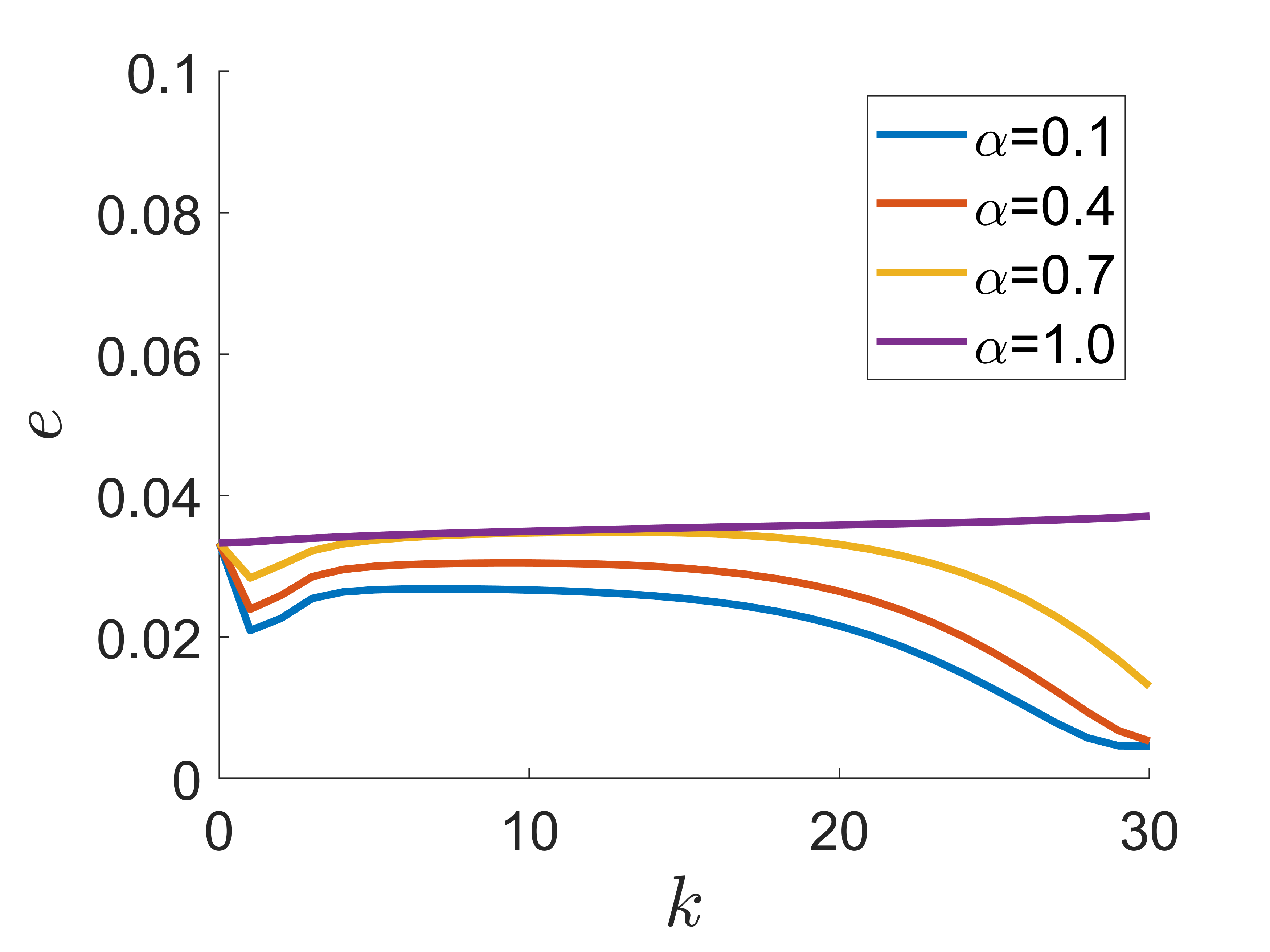} \\
 (d) error of $x$  & (e) error of $\lambda$ & (f)  error of $u_0$
\end{tabular}
\caption{Numerical results for Example \ref{exam2}(ii): the recovered $\lambda_i$ and $u_0$ with $\alpha=0.1$ (top) and the error of the recovered $\lambda_i$ and $u_0$ at different $\alpha$ (bottom).}
\label{fig:exp_3}
\end{figure}

In case (ii), the initial guesses for the source parameters are identical with that for (i), and that
for $u_0$ is $1.2x(1-x)$. The time interval is fixed to be $(\frac12T,T)$ and the
noise level $\delta=2\%$. To recover $u_0$, we add a penalty $\norm{u_0-u_0^k}_{H^1(\Omega)}^2$ to the functional $J_k$. The numerical results are shown in Fig.~\ref{fig:exp_3} and Table~\ref{tab:exp_3}.
The recovered intensities $\lambda_i(t)$ and the initial data are very accurate. The accuracy in cases (i) and (ii) is largely comparable.
For all $\alpha$ values, the source locations and intensities are well recovered. The convergence behavior of
the algorithm does not depend much on $\alpha$: the larger is $\alpha$, the faster is the convergence of the source
location. However, the behavior is reversed for $u_0$, and moreover, for $\alpha$ close to 1, the
algorithm hardly converges for $u_0$. Indeed for $\alpha=1$, the problem is local, and the recovery of $u_0$ is severely ill-posed. This might have contributed to the slow convergence of the algorithm for recovering $u_0$.

\begin{table}[htb!]
    \centering
    \renewcommand{\arraystretch}{1.5}
    \begin{tabular}{c|cccc}
    \toprule
    $\alpha$ & $0.1$ & $0.4$ & $0.7$ & $1$ \\
    \midrule
    location & \begin{tabular}{@{}c@{}} $0.2999$ \\  $0.6999$ \end{tabular} & \begin{tabular}{@{}c@{}} $0.2998$ \\  $0.6997$ \end{tabular} & \begin{tabular}{@{}c@{}} $0.2994$ \\  $0.6998$ \end{tabular} & \begin{tabular}{@{}c@{}} $0.2993$ \\  $0.6996$ \end{tabular} \\
    \bottomrule
    \end{tabular}
    \caption{The recovered source locations for Example~\ref{exam2}(ii).}
    \label{tab:exp_3}
\end{table}

\subsection{2D numerical examples}
Now we show 2D numerical examples with $\Omega=(0,1)^2$, $T=1$, $\rho\equiv1$, and $\mathcal{A}=-\Delta$. The direct problem is solved on a fine triangular mesh derived from a uniform grid over $\Omega$ with a mesh size $h=5\times10^{-3}$, and the time step size $\tau=2.5\times10^{-3}$.
The inverse problem is solved on a coarse mesh with $h=2\times10^{-2}$ and $\tau=10^{-2}$ at a noise level $\delta=2\%$.
The initial data is set to $u_0(x_1,x_2)=x_1x_2(1-x_1)(1-x_2)$. 

\begin{example}\label{exam4} Consider one point source at $x=(0.4, 0.4)$ with an intensity $\lambda(t)=0.5e^t$. The  data is given over the region $\omega=(0.5,1)^2$ for the time interval $(T-\epsilon,T)$ with $\epsilon=\frac34T$ or $\epsilon=\frac12T$.
\end{example}

\begin{figure}[htb!]
\centering\setlength{\tabcolsep}{0pt}
\begin{tabular}{cccc}
	\includegraphics[width=.24\textwidth]{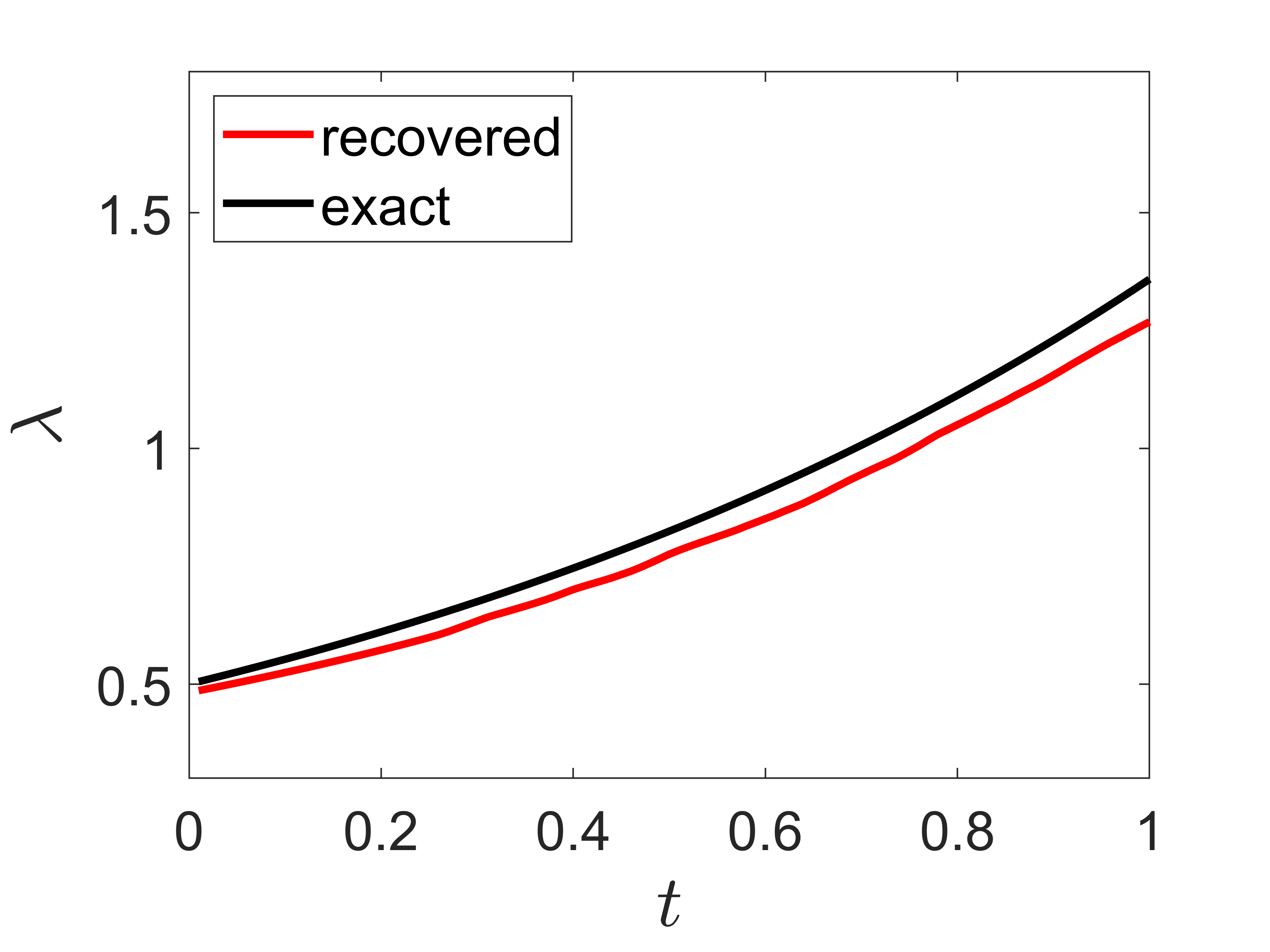} & \includegraphics[width=.24\textwidth]{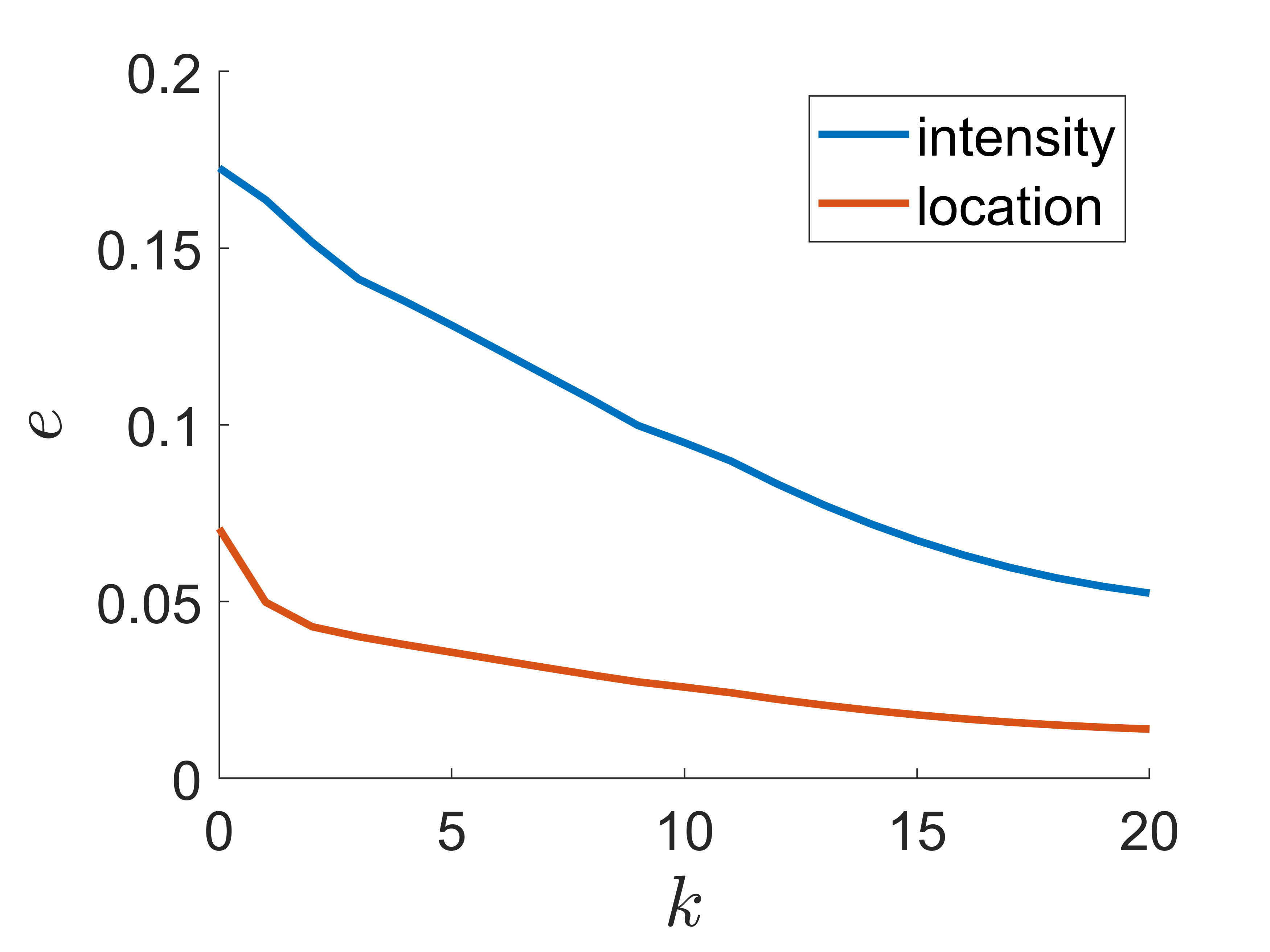} & \includegraphics[width=.24\textwidth]{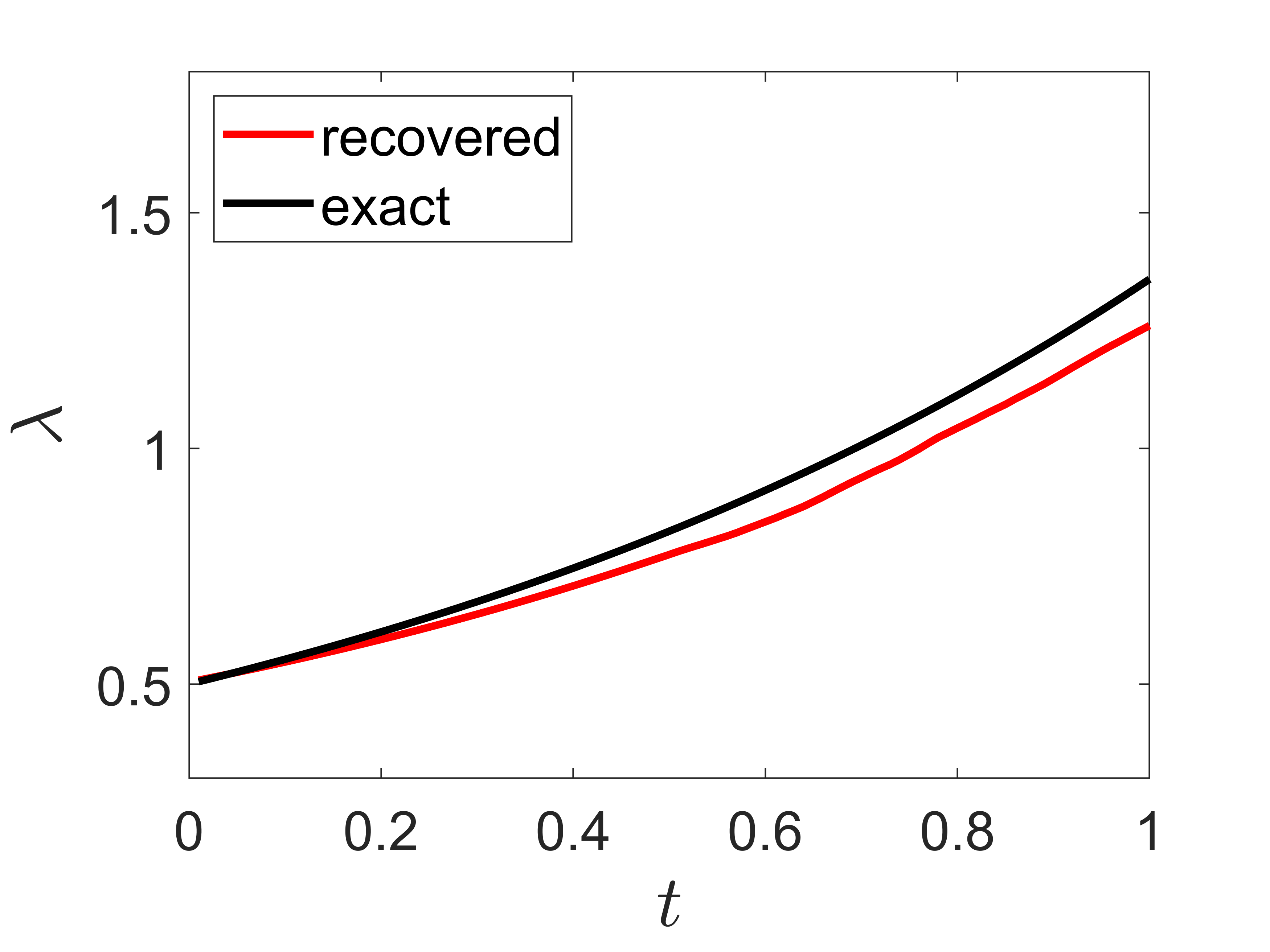} & \includegraphics[width=.24\textwidth]{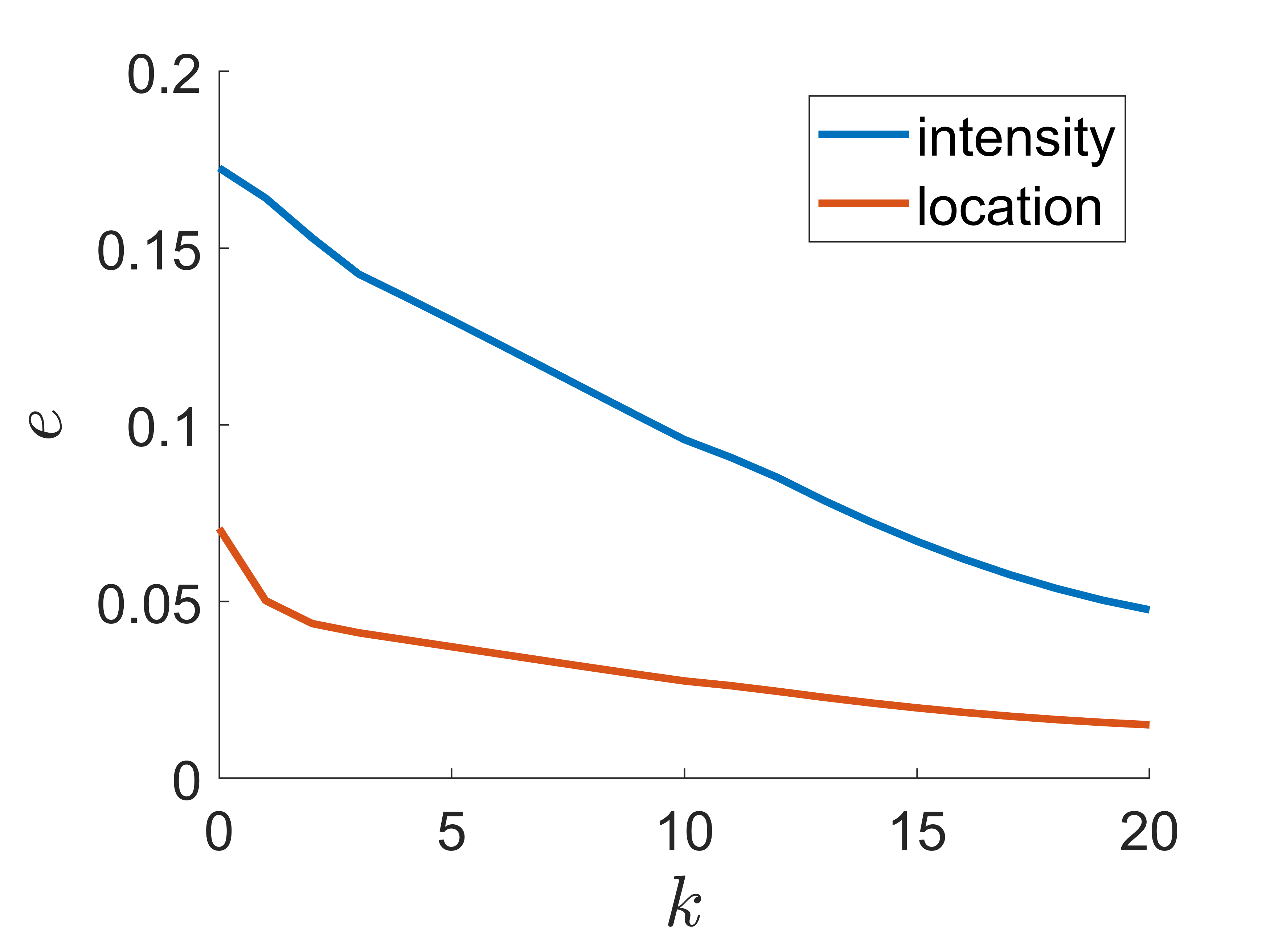}\\
    (a) $\lambda$ &  (b) convergence & (c) $\lambda$ & (d) convergence
\end{tabular}
\caption{The recovered $\lambda_i$ and convergence  for Example~\ref{exam4}. (a)-(b) and (c)-(d) are for $\epsilon= \frac34T$ and  $\epsilon=\frac12T$, respectively.}
\label{fig:exp_4}
\end{figure}

The algorithm is initialized with $x^0=(0.45, 0.45)$ and $\lambda^0(t)=0.4e^t$. The recovered source location and intensity
are fairly reasonable for both time intervals. See Table~\ref{tab:exp_4} and Fig.~\ref{fig:exp_4} for the recovered
location and intensity. The error of the recovered intensity is larger than that for the recovered location, and also less accurate than that for the 1D case.

\begin{table}[htb!]
    \centering
    \renewcommand{\arraystretch}{1.5}
    \begin{tabular}{c|cc}
    \toprule
    example  $\backslash \epsilon$ & $\frac34T$ & $\frac12T$\\
    \midrule
    \ref{exam4} & $(0.4111,0.4084)$ & $(0.4113,0.4100)$ \\
    \midrule
    \ref{exam5} & \begin{tabular}{@{}c@{}} $(0.5015,0.8946)$ \\  $(0.4995,0.1056)$ \end{tabular} & \begin{tabular}{@{}c@{}} $(0.4993,0.8931)$ \\  $(0.5050,0.1092)$ \end{tabular} \\
    \bottomrule
    \end{tabular}
    \caption{The recovered source locations for Examples~\ref{exam4} and \ref{exam5} at two $\epsilon$.}
    \label{tab:exp_4}
\end{table}

\begin{example} \label{exam5}
Consider two point sources at $x_1=(0.5, 0.9)$ and $x_2=(0.5, 0.1)$ with intensities $\lambda_1(t)=0.5e^t$ and $\lambda_2(t)=0.5\sin(2\pi t)+1$ for $t\in[0,T]$. The data is given over the disk centered at $(0.5,0.5)$ with a radius $0.3$ for the time interval $(T-\epsilon,T)$ with $\epsilon=\frac34T$ or $\epsilon=\frac12T$.
\end{example}

\begin{figure}[htbp!]
\centering\setlength{\tabcolsep}{0pt}
\begin{tabular}{ccc}
    \includegraphics[width=.32\textwidth]{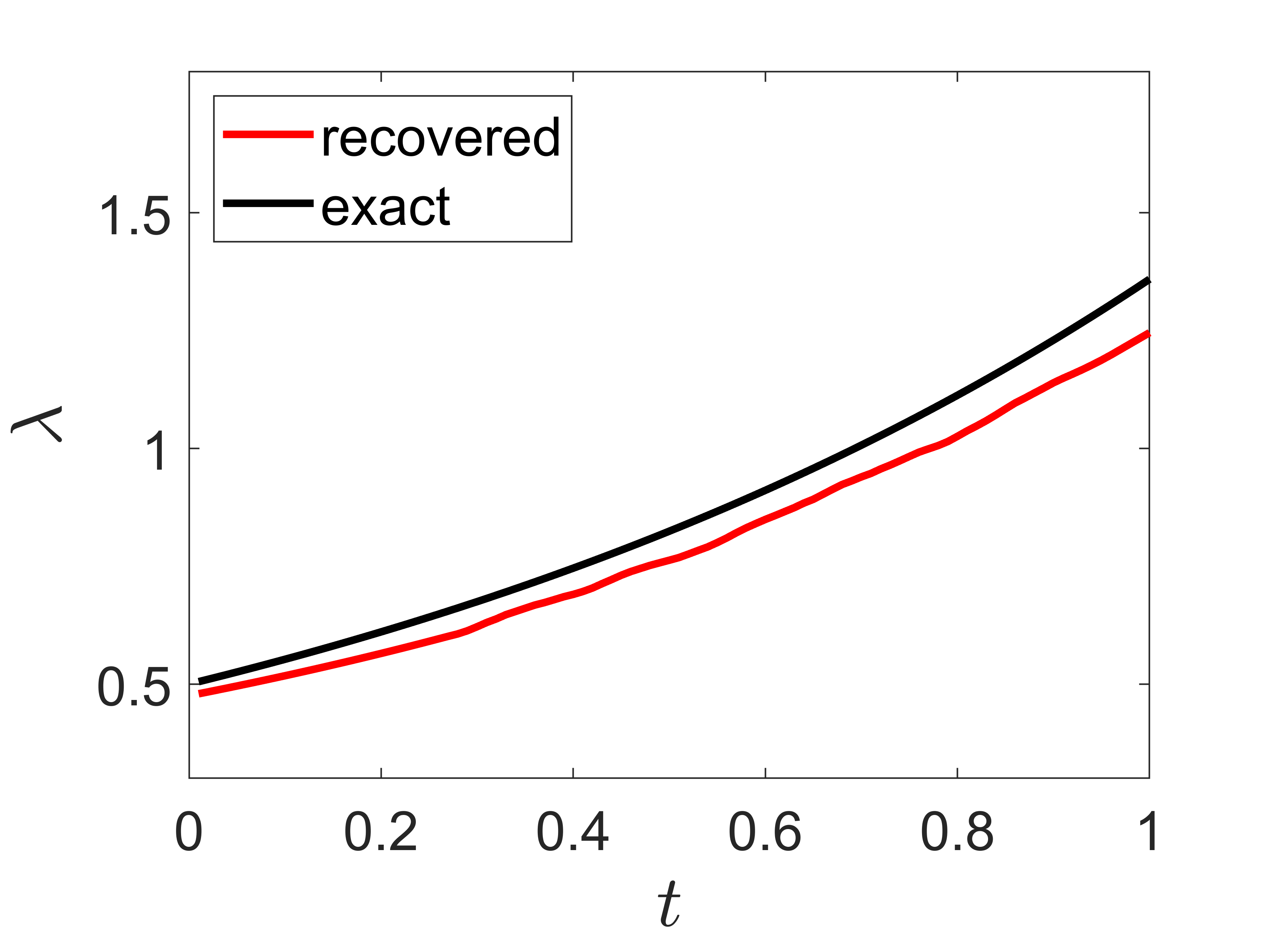} & \includegraphics[width=.32\textwidth]{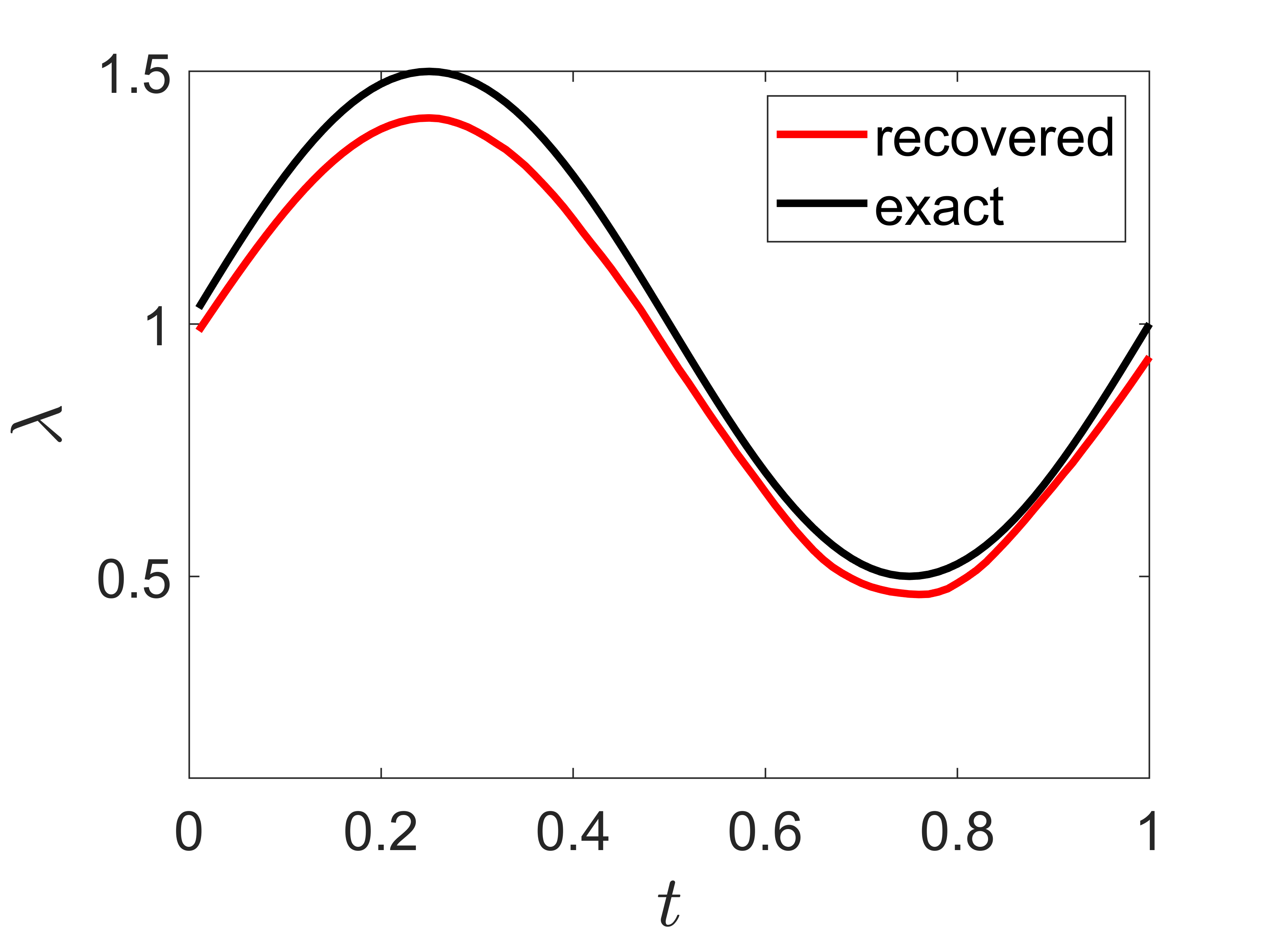} & \includegraphics[width=.32\textwidth]{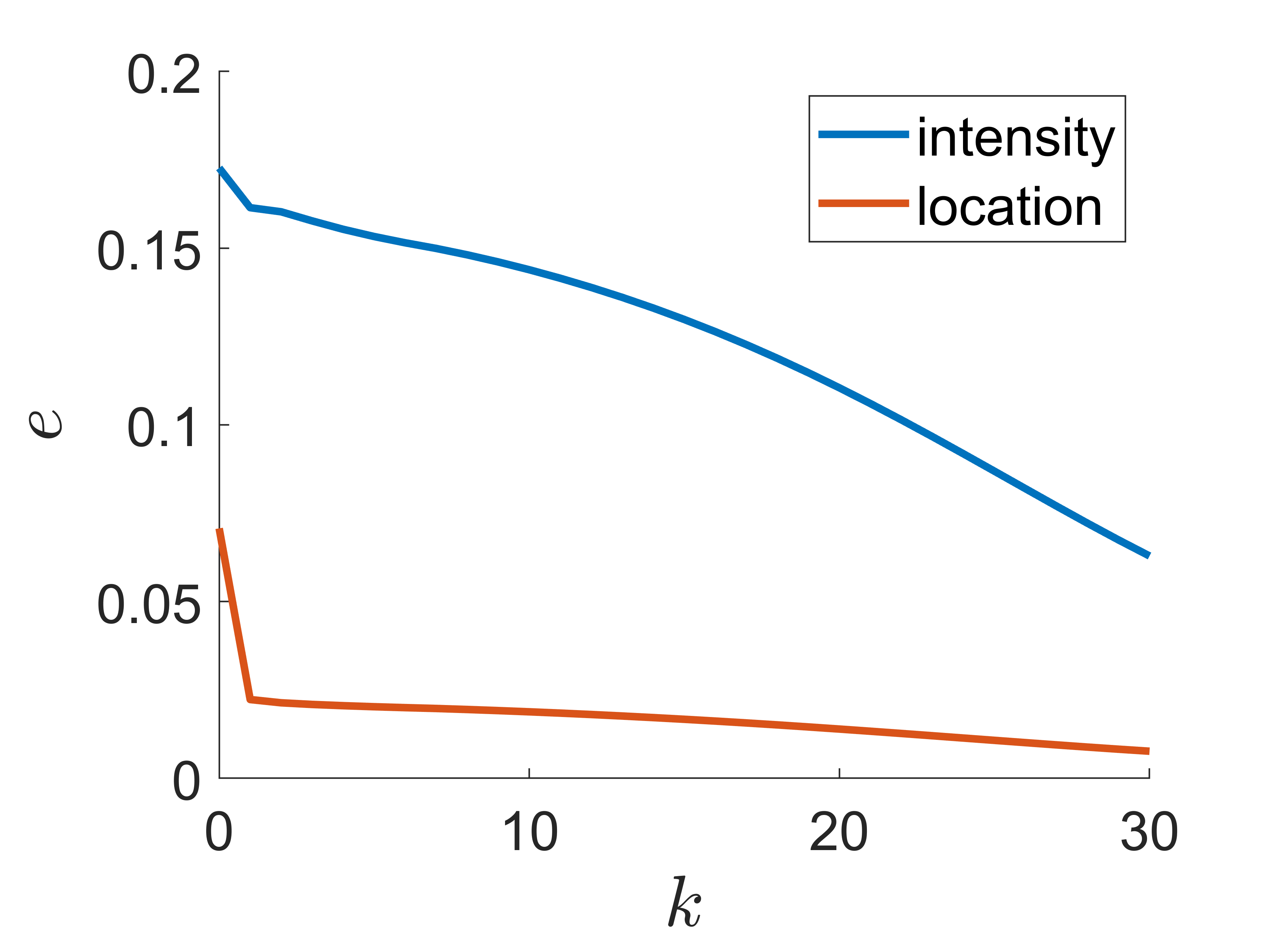} \\
    \includegraphics[width=.32\textwidth]{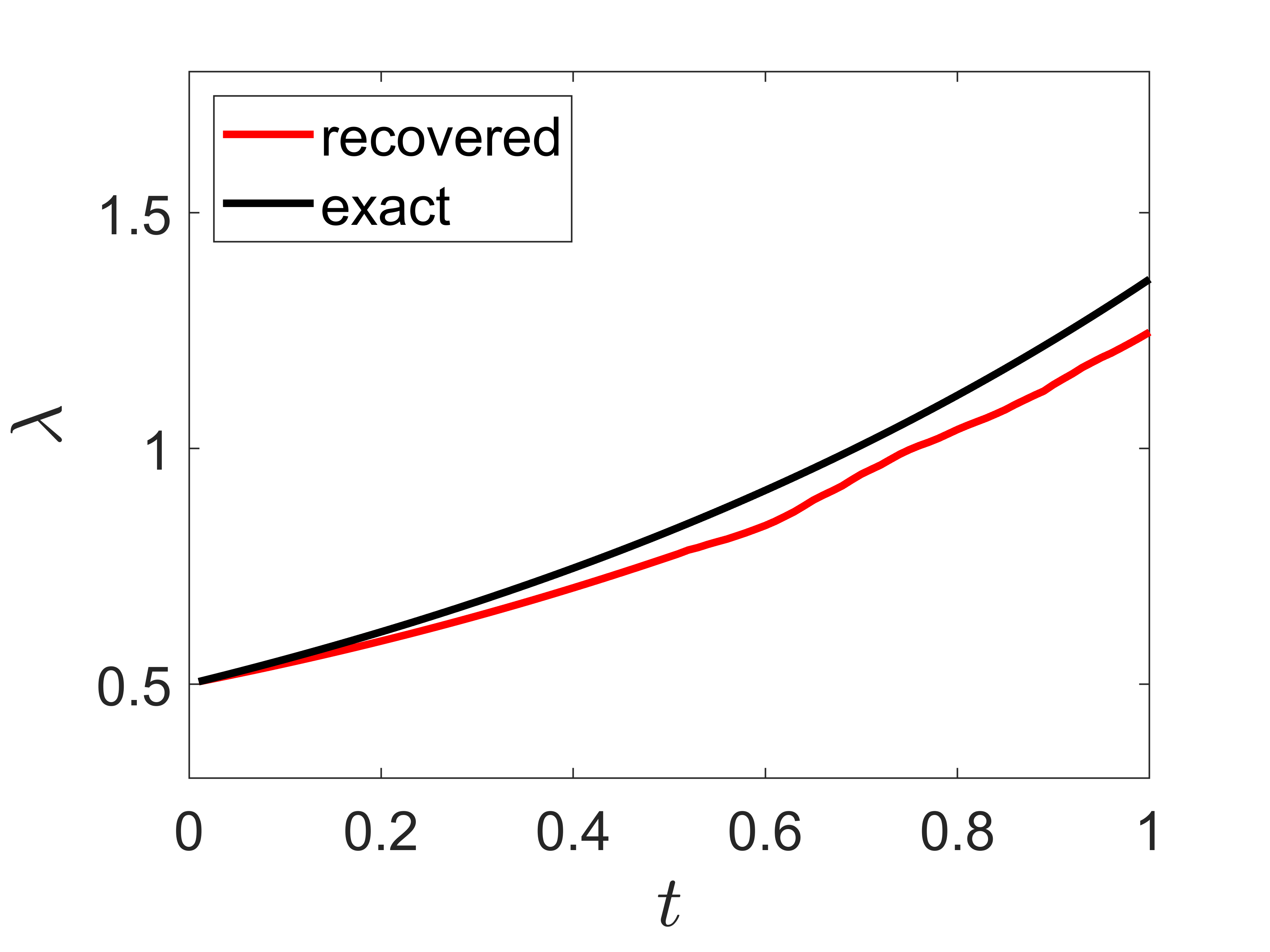} & \includegraphics[width=.32\textwidth]{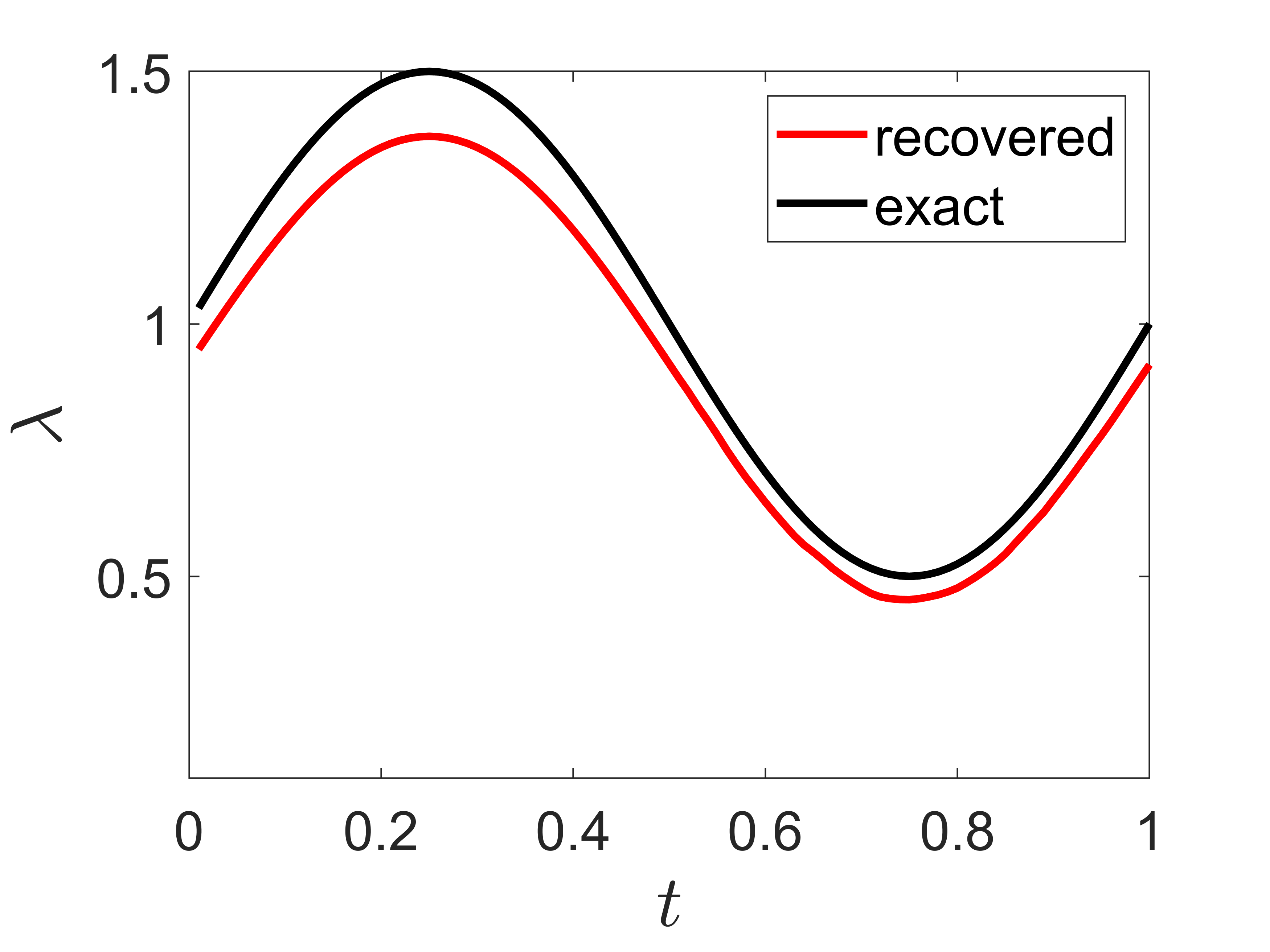} & \includegraphics[width=.32\textwidth]{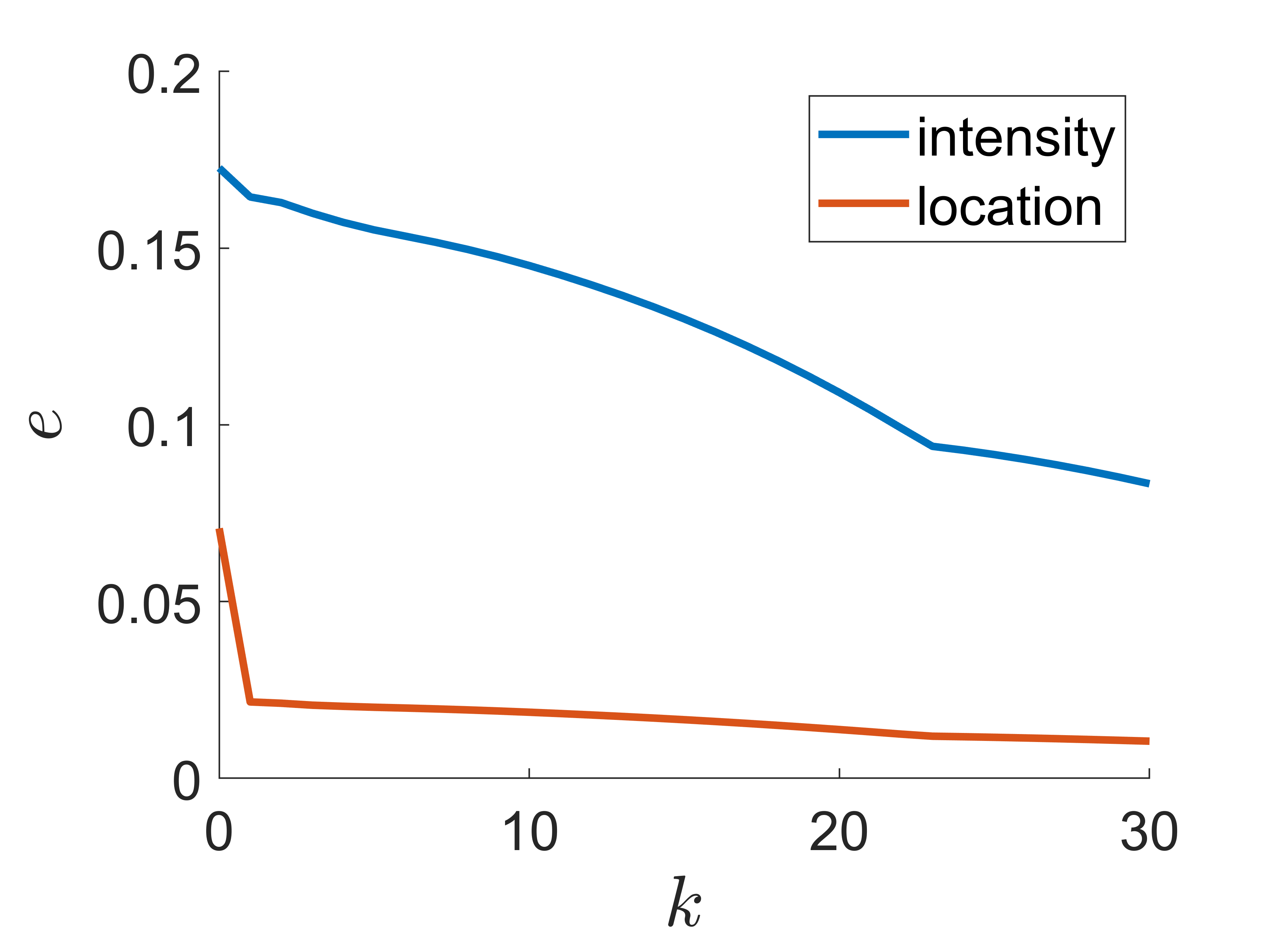} \\
    (a) $\lambda_1$ & (b) $\lambda_2$  & (c) convergence
\end{tabular}
\caption{The recovered $\lambda_i$ and convergence for Example~\ref{exam5} with $\epsilon=\frac34T$ (top) and $\epsilon=\frac12T$ (bottom).}
\label{fig:exp_5}
\end{figure}

The algorithm is initialized with $x_1^0=(0.45,0.85)$, $x_2^0=(0.45,0.15)$, $\lambda_1^0(t)=0.4e^t$ and $ \lambda_2^0(t)=0.45\sin(2\pi t)+0.9$. The recovered source locations and intensities are reasonably accurate, and the choice of the time interval does not influence much the reconstruction accuracy; see Fig.~\ref{fig:exp_5} and Table~\ref{tab:exp_4} for the recovered intensities and locations.

\subsection{3D numerical examples}
Now we use the open source finite element software package \texttt{FreeFem++} (\url{https://freefem.org}) \cite{hecht:2012} to solve PDEs. We show 3D numerical examples with $\Omega=(0,1)^3$, $T=1$, $\rho\equiv1$, and $\mathcal{A}=-\Delta$. To resolve the strong singularity of solutions, the direct problem is solved on an adaptive tetrahedral mesh used in \cite{Volpert:2024} with the minimum and maximum tetrahedron edge sizes $h_{\mathrm{min}}=3\times10^{-3}$ and $h_{\mathrm{max}}=0.12$, and a time step size $\tau=5\times10^{-3}$, and
the inverse problem on a coarse grid with $h_{\mathrm{min}}=7\times10^{-3},h_{\mathrm{max}}=0.15$ and $\tau =2\times10^{-2}$.
Moreover, due to the strong singularity of solution, the noisy data $z^\delta$ is generated by adding i.i.d. noise from the uniform distribution $\mathcal{U}[-\delta,\delta]$ to the exact data, with $\delta$ being fixed at $2\%$.

\begin{example}\label{exam6}  
The initial data $u_0(x_1,x_2,x_3)=\prod_{i=1}^3x_i(1-x_i)$. The data is given over the sphere centered at $(0.5,0.5,0.5)$ with a radius $0.3$ for the time interval $(\frac14T,T)$ with $\epsilon=\frac34T$.
\begin{itemize}
    \item[{\rm(i)}] One point source at $x=(0.5, 0.5, 0.9)$ with intensity $\lambda(t)=0.5e^t$.
    \item[{\rm(ii)}] Two point sources at $x_1=(0.5, 0.5, 0.9)$ and $x_2=(0.5, 0.5, 0.1)$ with intensities $\lambda_1(t)=0.5e^t$ and $\lambda_2(t)=0.5\sin(2\pi t)+1$.
\end{itemize}
\end{example}

\begin{figure}[htbp!]
\centering\setlength{\tabcolsep}{0pt}
\begin{tabular}{cc}    \includegraphics[width=.32\textwidth]{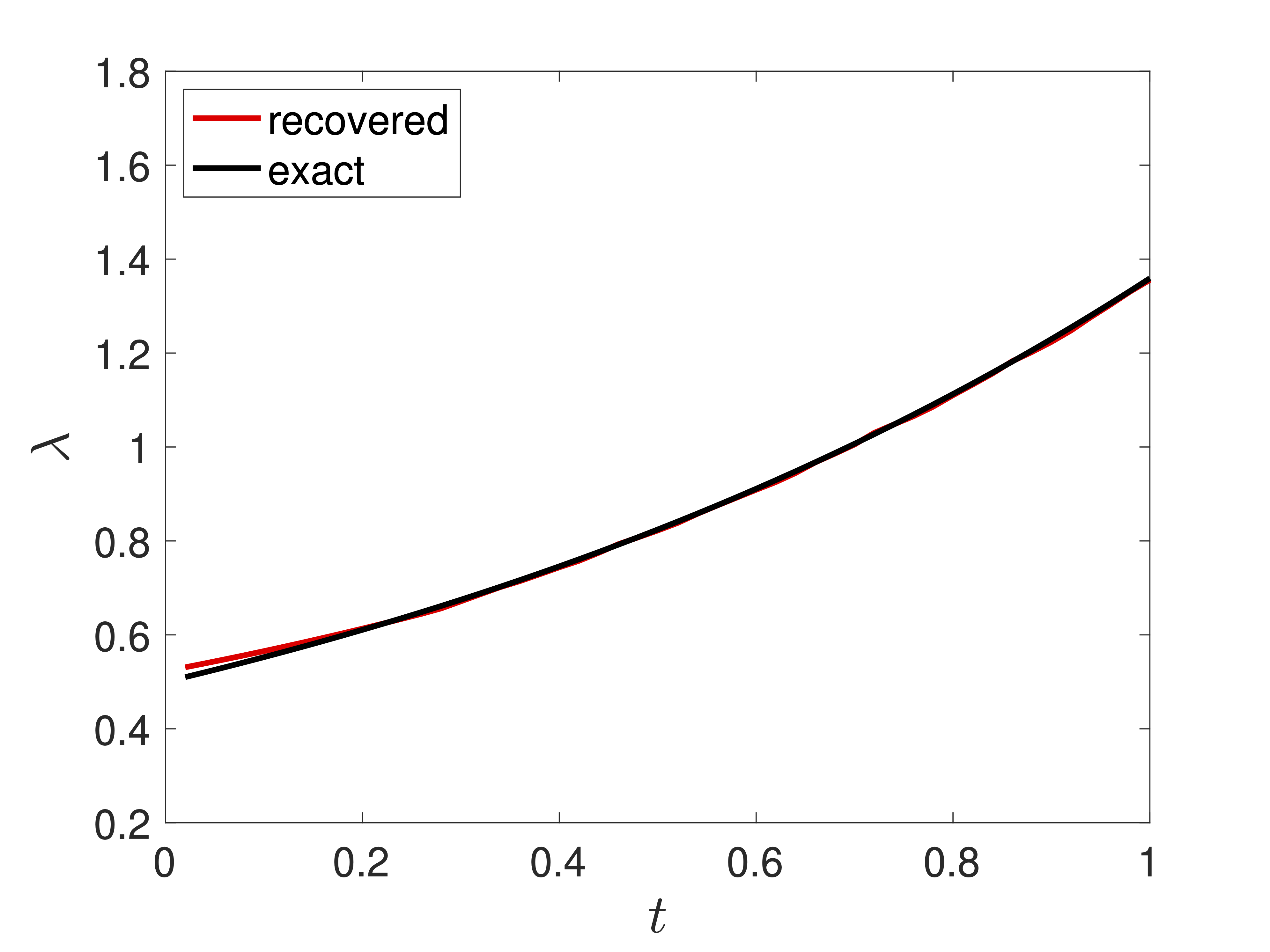} & \includegraphics[width=.32\textwidth]{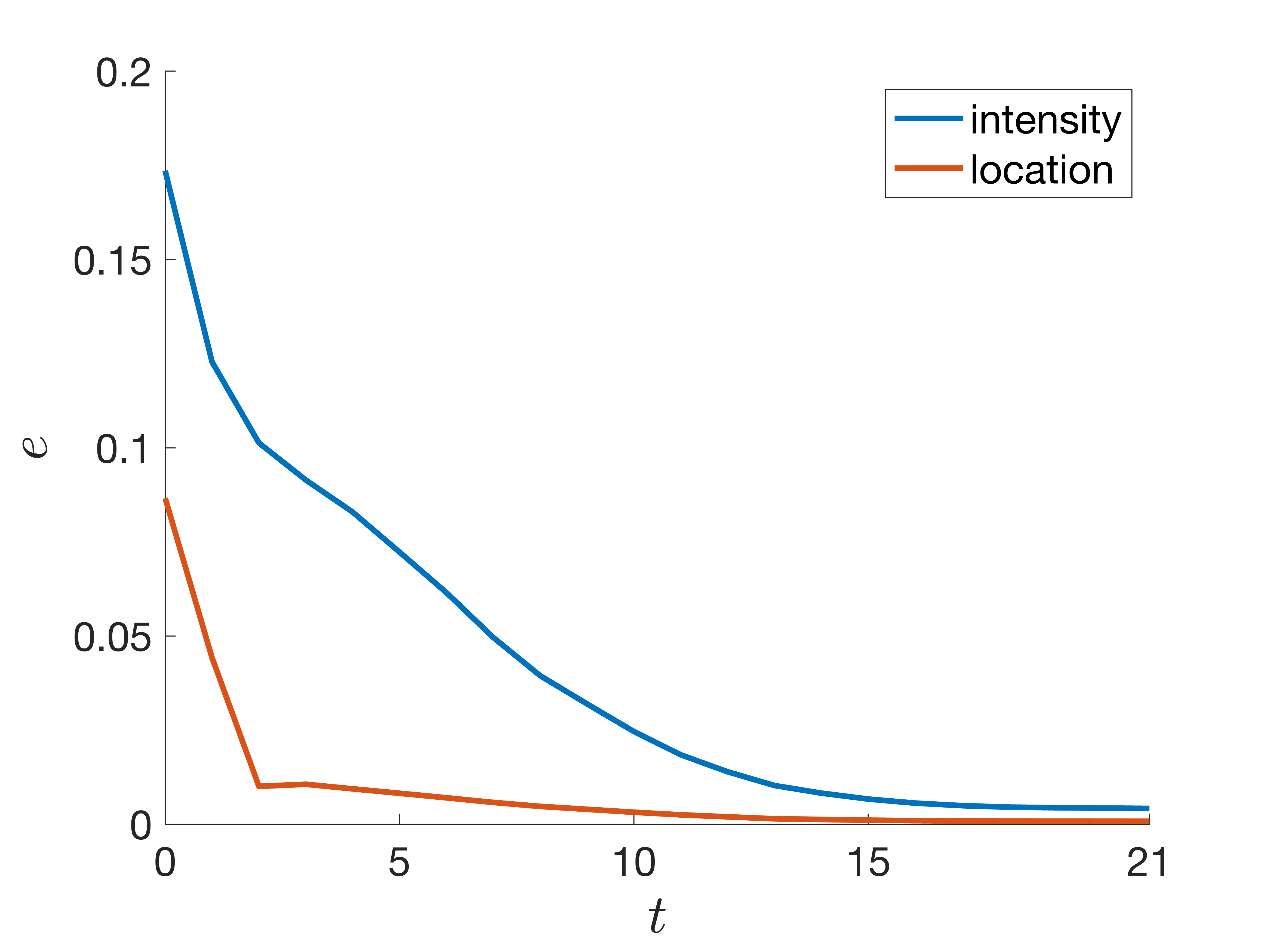}\\
(a) $\lambda$ & (b) convergence
\end{tabular}
\caption{The recovered $\lambda$ and convergence  for Example~\ref{exam6}(i) with $\epsilon=\frac34T$.}
\label{fig:exp_6}
\end{figure}

\begin{figure}[htbp!]
\centering\setlength{\tabcolsep}{0pt}
\begin{tabular}{ccc}
    \includegraphics[width=.32\textwidth]{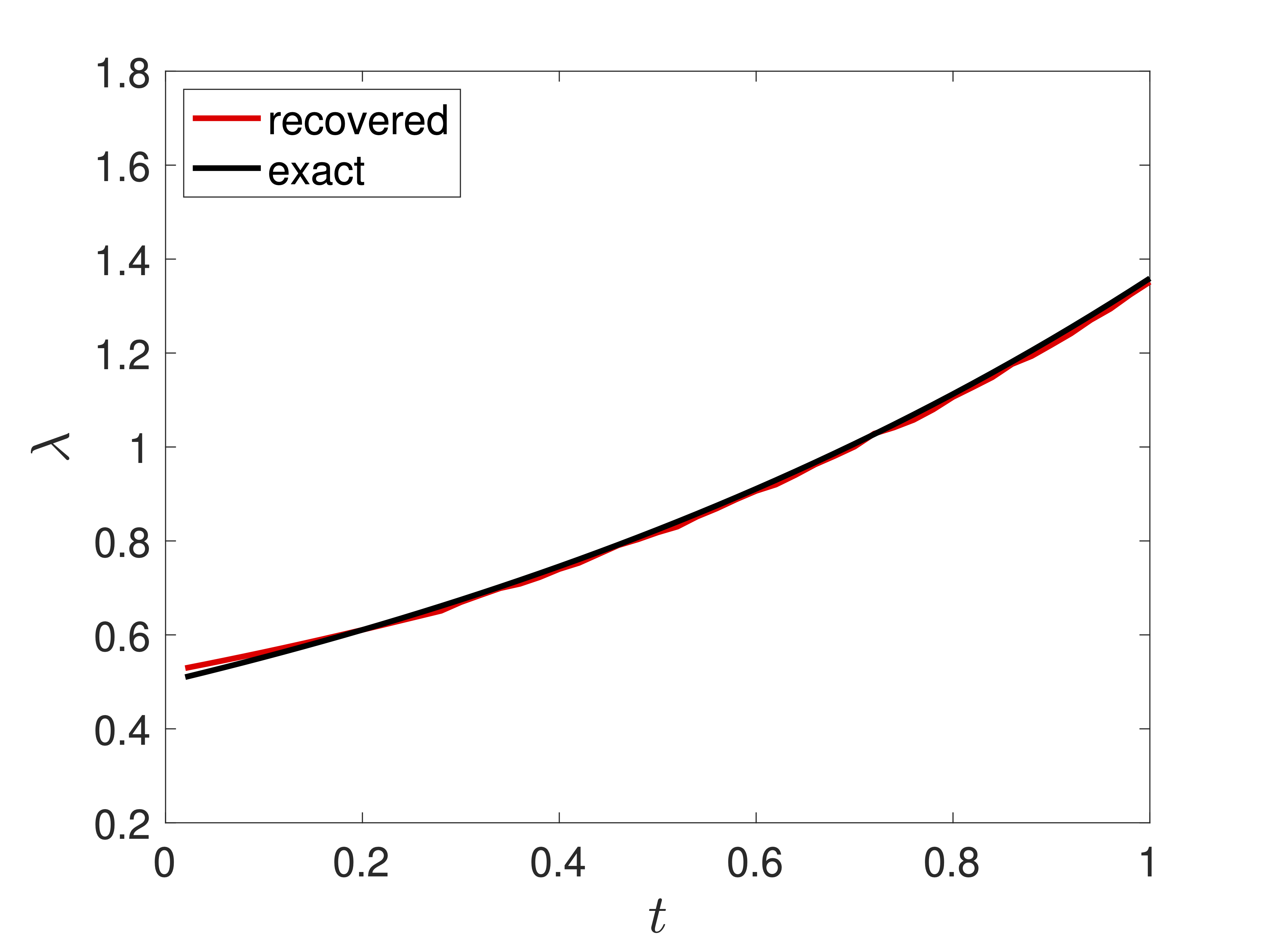} & \includegraphics[width=.32\textwidth]{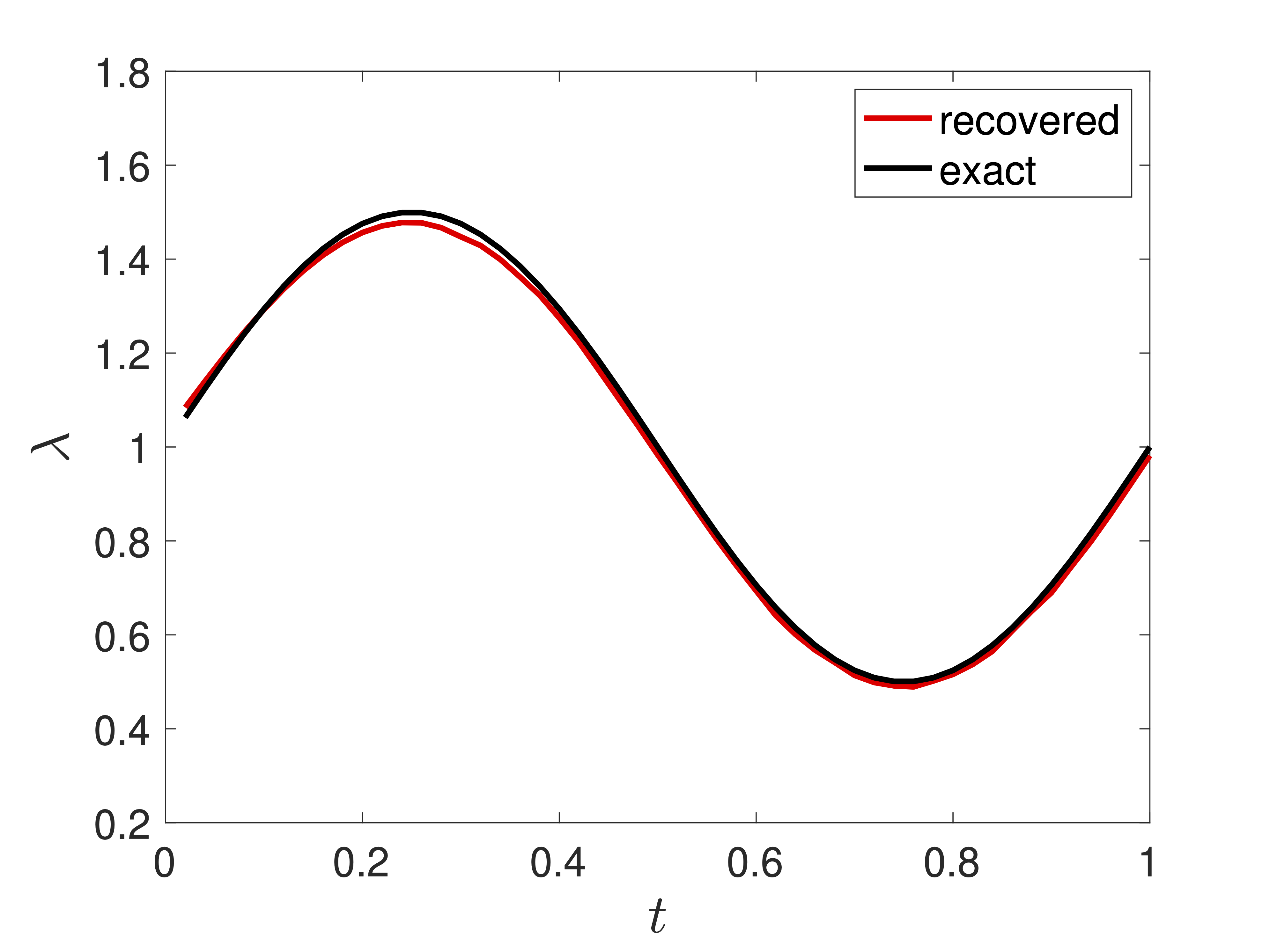} & \includegraphics[width=.32\textwidth]{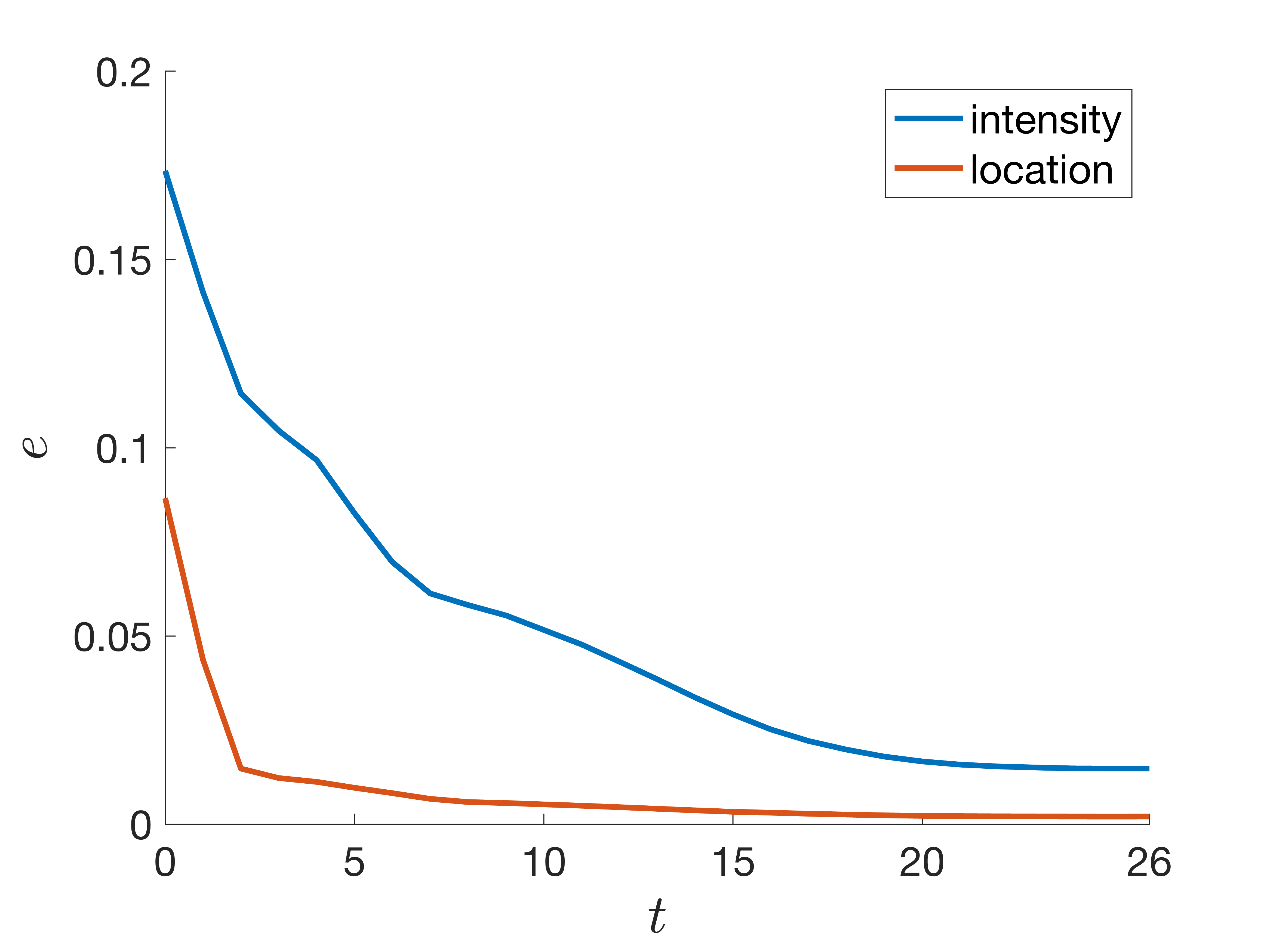}\\
    (a) $\lambda_1$ & (b) $\lambda_2$  & (c) convergence
\end{tabular}
\caption{The recovered $\lambda_i$ and convergence  for Example~\ref{exam6}(ii) with $\epsilon=\frac34T$.}
\label{fig:exp_7}
\end{figure}

In case (i), the algorithm is initialized with $x^0=(0.45, 0.45, 0.85)$ and $\lambda^0(t)=0.4e^t$. The recovered source location is $(0.5001,0.4995,0.8994)$, which is fairly accurate. See Fig.~\ref{fig:exp_6} for the recovered intensity, which has an error larger than that for the recovered location. In case (ii), the algorithm is initialized with $x_1^0=(0.45, 0.45, 0.85)$, $x_2^0=(0.45, 0.45, 0.15)$, $\lambda_1^0(t)=0.4e^t$ and $ \lambda_2^0(t)=0.45\sin(2\pi t)+0.9$. The recovered locations are $(0.4996,0.4997,0.1019)$ and   $(0.5000,0.4995,0.8988)$, which are again very accurate. See Fig. \ref{fig:exp_7} for the recovered intensities. For both cases (i) and (ii), the convergence of the algorithm is fairly stable.

\vspace{-0.8cm}
\section{Conclusion}
We have investigated the inverse problem of recovering point sources, their strengths and initial data from a posteriori partial internal observational data. We prove two unique identifiability results for the subdiffusion model involving general elliptic operators, including time-dependent coefficients in the elliptic operator in the one-dimensional case. The analysis relies on the classical unique continuation principle of the subdiffusion model, and improved local regularity of the solution. Numerical results obtained using the Levenberg-Marquardt algorithm also show the feasibility of the reconstruction. Future research topics include stability estimates of the inverse problem, and moving point sources etc.

\bibliographystyle{abbrv}
\bibliography{frac}
\end{document}